\documentclass[11pt,a4paper]{article}
\usepackage{mathtools}
\usepackage{graphicx}
\usepackage[utf8]{inputenc}
\usepackage[T1]{fontenc}
\usepackage{amsmath,amssymb,amsthm}
\usepackage{mathrsfs}
\usepackage{geometry}
\usepackage{hyperref}
\hypersetup{colorlinks=true,linkcolor=blue,citecolor=blue,urlcolor=blue}

\newtheorem{theorem}{Theorem}[section]
\newtheorem{lemma}[theorem]{Lemma}
\newtheorem{proposition}[theorem]{Proposition}
\newtheorem{corollary}[theorem]{Corollary}
\theoremstyle{definition}
\newtheorem{definition}[theorem]{Definition}
\newtheorem{remark}[theorem]{Remark}

\newtheorem{example}[theorem]{Example}

\newcommand{\R}{\mathbb{R}}

\newcommand{\Hcal}{\mathcal{H}}
\newcommand{\Acal}{\mathcal{A}}
\newcommand{\Lcal}{\mathcal{L}}
\newcommand{\Bcal}{\mathcal{B}}
\newcommand{\Mcal}{\mathcal{M}}
\newcommand{\Dcal}{\mathcal{D}}
\newcommand{\dist}{\operatorname{dist}}
\newcommand{\dd}{\,\mathrm{d}}
\newcommand{\eps}{\varepsilon}

\newcommand{\Kcal}{\mathcal{K}}

\newcommand{\Pcal}{\mathcal{P}}

\newcommand{\Sf}{\mathcal{S}_F}
\newcommand{\Mg}[1]{\mathcal{M}_{g_{#1,\eps}}}
\newcommand{\keywordsname}{Keywords}
\newcommand{\keywords}[1]{%
  \par\vspace{0.5\baselineskip}%
  \noindent\textbf{\keywordsname:} #1\par
}
\newcommand{\subjclassname}{MSC 2020}
\newcommand{\subjclass}[2][]{%
  \par\vspace{0.3\baselineskip}%
  \noindent\textbf{\subjclassname:} #2\par
}

\title{Continuity of Global Attractors under a Mass-Preserving Singular Memory Limit, with an Application to a Coupled Suspension Bridge System\footnote{This work is supported by NSFC (12671232).}}

\author{Jun Zhou\footnote{Email: jzhou@swu.edu.cn}\\
        School of Mathematics and Statistics, Southwest University,\\
        Chongqing 400715, People's Republic of China}

\date{}

\begin{document}

\maketitle

\begin{abstract}
We study continuity of global attractors for a damped coupled suspension bridge system with past history when the memory kernel concentrates at the origin. The decisive condition is uniform positivity of the effective elastic coefficients $\mu_i-\beta_{i,\eps}$; mass preservation is the simplest way to guarantee it. We prove uniform-in-$\eps$ dissipativity and exponential stability of the linear semigroup, uniform quasi-stability, finite-time convergence to the memoryless limit, and upper semicontinuity of the attractors. Under a structural hypothesis on the limit attractor we also obtain lower semicontinuity. The analysis separates the abstract continuity mechanism from its model-specific verification. Further results include continuity under a nondegenerate mass drift, the identification of the fast-decay limit with the same singular limit, and a conditional reduction of quantitative lower semicontinuity to a finite-time convergence rate.
\end{abstract}




\section{Introduction}\label{sec:intro}

\subsection{The general problem}

A recurrent question in the theory of dissipative infinite-dimensional dynamical systems is how the long-time behaviour depends on a small parameter in the equation. When the parameter enters through a memory kernel that becomes concentrated in the limit, the perturbed and limiting systems live in different function spaces, and the natural object to compare is the family of global attractors. The present paper addresses this question for a class of memory kernels that concentrate at the origin while preserving their total mass. The decisive structural condition is uniform positivity of the effective elastic coefficients $\mu_i-\beta_{i,\eps}$; the mass-preserving rescaling is the cleanest normalisation that guarantees it, while allowing the same estimates to be uniform in the perturbation parameter.

The prototype situation is the following. Let $g$ be an integrable memory kernel on $\R_+$ with total mass $\beta$, and for $\eps>0$ set
\begin{equation}\label{eq:intro_rescale}
g_\eps(s)\coloneqq\frac{1}{\eps}g\!\Big(\frac{s}{\eps}\Big).
\end{equation}
Then $\int_0^{+\infty}g_\eps(s)\dd s=\beta$ for every $\eps>0$, while $g_\eps\rightharpoonup\beta\delta_0$ in the distributional sense as $\eps\to0^+$. In a viscoelastic equation in which the memory term enters through a convolution against an elastic operator, the coefficient of that operator is effectively shifted by the mass of the kernel. A rescaling that changes the mass moves the limiting coefficient to $\mu_i-\beta_{i,0}$; this is harmless as long as $\mu_i-\beta_{i,0}>0$, and genuinely singular only when coercivity is lost. The mass-preserving rescaling \eqref{eq:intro_rescale} is therefore a convenient normalisation rather than the essential mechanism: it keeps $\bar\mu_i=\mu_i-\beta_i$ fixed and isolates the role of the positivity condition.

Our aim is twofold. First, we develop a continuity theory for global attractors under \eqref{eq:intro_rescale} that is not tied to a particular model: the ingredients are uniform dissipativity, uniform exponential stability of the linear semigroup, uniform quasi-stability in the sense of Chueshov--Lasiecka, finite-time convergence of the semigroups, and, for lower semicontinuity, a structural stability hypothesis on the limit attractor. Second, we verify these ingredients for a coupled nonlinear suspension bridge system with past history, which serves as a concrete and physically motivated model in which the general mechanism can be checked in detail. The bridge system is thus an application of the framework rather than its sole object.

\subsection{The model}

Suspension bridges are important civil engineering structures, and their dynamical behaviour has been widely studied. Early work by McKenna and co-authors \cite{GloverLazerMcKenna1989,LazerMcKenna1987,LazerMcKenna1990,McKennaWalter1987,McKennaWalter1990,McKenna1999,McKennaTuama2001} proposed and analysed coupled models for vertical bridge vibrations of the form
\begin{equation}\label{eq:intro_bridge}
\begin{cases}
m_1u_{tt}+\mu_1u_{xxxx}+k[u-v]^+ +\gamma_1u_t=f,~~~x\in(0,\ell),~t>0,\\[1mm]
m_2v_{tt}-\mu_2v_{xx}-k[u-v]^+ +\gamma_2v_t=g,~~~x\in(0,\ell),~t>0,
\end{cases}
\end{equation}
where $u=u(x,t)$ denotes the vertical displacement of the bridge deck and $v=v(x,t)$ the displacement of the main cable. The first equation describes deck oscillations in the vertical plane, the second the motion of the main cable supported by suspender cables. The positive constants $m_1,m_2$ are masses per unit length, $k$ is the spring constant of the suspenders, $\mu_1,\mu_2$ are the flexural rigidity of the deck and the tensile stiffness of the cable, $\gamma_1,\gamma_2$ are damping coefficients, and $f,g$ are external forces. The nonlinear term $k[u-v]^+$ with $[u-v]^+=\max\{u-v,0\}$ models the restoring force of tensioned suspenders; it vanishes when the deck lies above the cable, reflecting that suspenders can only exert tension. The boundary conditions are
\begin{equation}\label{eq:intro_bc}
u(0,t)=u(\ell,t)=u_{xx}(0,t)=u_{xx}(\ell,t)=0,\qquad
v(0,t)=v(\ell,t)=0,
\end{equation}
for a simply supported deck beam and a pinned main cable. For simplicity we set $m_1=m_2=1$ and consider time-independent external forces $f_i=f_i(x)$, $i=1,2$, so that the system is autonomous.

A substantial body of work has addressed well-posedness, asymptotic stability, and global attractors for \eqref{eq:intro_bridge} and its modifications; see, e.g., \cite{AhmedHarbi1998,An2002,DellOroGiorgiPata2015,HolubovaMatas2003,MaZhong2005,Malik2004} and references therein. Ma and Zhong \cite{MaZhong2009} considered the coupled system
\begin{equation}\label{eq:MaZhong}
\begin{cases}
u_{tt}+\mu_1u_{xxxx}+k[u-v]^+ +\gamma_1u_t+h_1(u)=f_1,~~~x\in(0,\ell),~t>0,\\[1mm]
v_{tt}-\mu_2v_{xx}-k[u-v]^+ +\gamma_2v_t+h_2(v)=f_2,~~~x\in(0,\ell),~t>0,
\end{cases}
\end{equation}
and proved existence of strong solutions and a global attractor under the assumptions $h_1\in C^3(\R,\R)$, $h_2\in C^2(\R,\R)$, and $f_1,f_2\in L^2(0,\ell)$. Related results for stationary and evolutionary problems may be found in \cite{An2002,MaZhong2005}. Dell'Oro, Giorgi, and Pata \cite{DellOroGiorgiPata2015} studied asymptotic behaviour for linear coupled bridge models and derived decay estimates under suitable damping and structural assumptions. These works show that the nonlinear coupling $k[u-v]^+$ plays an important part in the long-time dynamics, and that attractor analysis requires careful treatment of its non-smooth nature.

If cable motion is neglected and the main cable is replaced by a rigid frame, deck vibration reduces to the scalar equation
\begin{equation}\label{eq:single_bridge}
u_{tt}+u_{xxxx}+k[u]^+ +\gamma u_t=f.
\end{equation}
We refer to \cite{AhmedHarbi1998,BochicchioGiorgiVuk2010,Kang2016,Litcanu2004,ParkKang2011,ZhongMaSun2007} and references therein for results concerning \eqref{eq:single_bridge}. In particular, Park and Kang \cite{ParkKang2011} constructed global attractors for bridge equations with nonlinear damping, while Zhong, Ma, and Sun \cite{ZhongMaSun2007} established well-posedness of strong solutions and existence of global attractors for the single beam model.

Because deck and cable materials may exhibit viscoelastic responses, there has been growing interest in suspension-bridge models including memory effects. Following Dafermos \cite{Dafermos1970}, past-history contributions can be incorporated by introducing relative displacement fields
\begin{equation}\label{eq:relative_disp}
\sigma(x,t,s)=u(x,t)-u(x,t-s),\qquad
\eta(x,t,s)=v(x,t)-v(x,t-s),
\end{equation}
which satisfy the transport equations
\begin{equation}\label{eq:transport_intro}
\sigma_t+\sigma_s-u_t=0,\qquad \eta_t+\eta_s-v_t=0.
\end{equation}
This formulation treats infinite past history as an additional state variable in an appropriate Hilbert space. Messaoudi, Bonfoh, Mukiawa, and Enyi \cite{MessaoudiBonfohMukiawaEnyi2017} studied a memory-dependent suspension-bridge system with partially hinged boundaries and obtained a global attractor. For beam-type systems with memory, Mukiawa and Enyi \cite{MukiawaEnyi2020,EnyiMukiawa2020} investigated well-posedness and long-time dynamics for Timoshenko and laminated beams subject to Fourier-type laws and infinite memory.

Mukiawa and Messaoudi \cite{Mukiawa2025} analysed the coupled nonlinear bridge system with past history
\begin{equation}\label{eq:Mukiawa}
\begin{cases}
u_{tt}+\gamma_1u_t+\mu_1u_{xxxx}+k[u-v]^+\\
\qquad-\displaystyle\int_{-\infty}^{t}g_1(t-s)u_{xxxx}(x,s)\dd s
+h_1(u)=f_1(x),~~~x\in(0,\ell),~t>0,\\[1mm]
v_{tt}+\gamma_2v_t-\mu_2v_{xx}-k[u-v]^+\\
\qquad+\displaystyle\int_{-\infty}^{t}g_2(t-s)v_{xx}(x,s)\dd s
+h_2(v)=f_2(x),~~~x\in(0,\ell),~t>0,
\end{cases}
\end{equation}
together with boundary conditions \eqref{eq:intro_bc} and suitable initial data, proving well-posedness and existence of a compact global attractor on the natural energy space, where the memory kernels $g_i$ satisfy
\begin{equation}\label{eq:kernel}
\int_0^{+\infty}g_i(s)\dd s=\beta_i\in(0,\mu_i),\qquad g_i'(s)+\lambda_i g_i(s)\le 0\quad\text{a.e. }s\ge 0,
\end{equation}
for some $\lambda_i>0$, $i=1,2$. The nonlinearities $h_i:\R\to\R$, $i=1,2$, are of class $C^1$, and their primitives
\begin{equation}\label{eq:primitive}
H_i(s)\coloneqq\int_0^s h_i(\tau)\dd\tau,\qquad i=1,2,
\end{equation}
satisfy the following growth and dissipativity conditions: there exist constants $C_1,C_2>0$ such that for all $s\in\R$ and $i=1,2$,
\begin{equation}\label{eq:h_assumptions}
\liminf_{|s|\to+\infty}\frac{H_i(s)}{s^2}\ge 0,\qquad
\liminf_{|s|\to+\infty}\frac{s\,h_i(s)-C_1H_i(s)}{s^2}\ge 0,\qquad
h_i'(s)\ge -C_2.
\end{equation}
As observed in \cite{Mukiawa2025,Zhou2026}, condition \eqref{eq:h_assumptions} implies that for every $\epsilon>0$ there exist constants $C_\epsilon,C_\epsilon'>0$ such that
\begin{equation}\label{eq:h_consequences}
H_i(s)+\epsilon s^2\ge -C_\epsilon,\qquad
s\,h_i(s)-C_1H_i(s)+\epsilon s^2\ge -C_\epsilon',\qquad \forall s\in\R,\ i=1,2.
\end{equation}
A typical example is $h_i(s)=s|s|^{p}$ with $0\le p<+\infty$, which satisfies \eqref{eq:h_assumptions} for suitable $C_1,C_2$. Throughout, $\eps$ denotes the singular parameter, while $\epsilon,\eta,\delta$ denote auxiliary parameters in estimates such as Young's inequality.

\subsection{The fixed-kernel attractor and a point in the literature}

We first recall the fixed-kernel attractor result used below and comment neutrally on one step in the literature. In \cite[Proposition~4.1]{Mukiawa2025}, the authors derive the estimate
\begin{equation}\label{eq:gap}
\|U^1(t)-U^2(t)\|_{\mathcal H}^2
\le c\|U_0^1-U_0^2\|_{\mathcal H}^2
+c_{\mathcal B}\int_0^t\big(\|u^1(s)-u^2(s)\|^2+\|v^1(s)-v^2(s)\|^2\big)\dd s,
\end{equation}
where $c_{\mathcal B}>0$ depends on the absorbing set $\mathcal B$, and then assert that for any $\varepsilon>0$ and sufficiently large $T>0$,
\begin{equation}\label{eq:gap2}
\|S(T)U_0^1-S(T)U_0^2\|_{\mathcal H}^2
\le\varepsilon+\Upsilon_T(U_0^1,U_0^2),
\end{equation}
with $\Upsilon_T$ defined in \cite[(4.26)]{Mukiawa2025}. As written, this step does not follow from \eqref{eq:gap} alone: the term $c\|U_0^1-U_0^2\|_{\mathcal H}^2$ is independent of $T$ and no decaying factor is displayed. We therefore read the argument as requiring an additional ingredient, namely exponential decay of the underlying linear semigroup in the form needed here. The argument as written does not yield \eqref{eq:gap2} without that ingredient; we do not claim that the issue cannot be resolved elsewhere.

Zhou \cite{Zhou2026} revisited system \eqref{eq:Mukiawa} and supplied the missing exponential-stability input. Using the Gearhart--Pr\"uss theorem and detailed resolvent estimates, exponential stability was shown for the associated linear semigroup. Combining this stability property with local Lipschitz behaviour of the nonlinear terms, the fact that nonlinearities act only on elastic components, smoothing effects of the elastic operators, and exponential decay of the memory kernels, quasi-stability in the Chueshov--Lasiecka sense follows, which yields asymptotic smoothness. Together with dissipativity results from \cite{Mukiawa2025}, this gives a compact global attractor on the natural energy space. In the present work we build the singular-limit analysis on this corrected global-attractor result for the fixed-kernel system, rather than directly on the original proof in \cite{Mukiawa2025}. This ingredient is an input for, rather than the main contribution of, the present paper.

\subsection{The singular perturbation and the continuity questions}

We consider a physically motivated singular perturbation of \eqref{eq:Mukiawa}, corresponding to a vanishing viscoelastic relaxation time. This setting belongs to the general theory of memory kernels converging towards a Dirac mass developed by Conti, Pata, and Squassina \cite{ContiPataSquassina2006}, who study singular limits of memory equations and the corresponding limiting dynamics. The novelty here is not the concentration rescaling itself but the attractor-continuity theory together with the uniform-in-$\eps$ estimates; see also Remark~\ref{rem:cp}.  Physically, the kernel $g_i$ encodes the relaxation modulus of the viscoelastic material: a short relaxation time means that the material forgets its past almost immediately, so that the memory term becomes negligible on time scales much longer than the relaxation time, leaving the dissipation provided by the original damping coefficients $\gamma_i$. We employ the mass-preserving rescaling
\begin{equation}\label{eq:rescaled}
g_{i,\eps}(s)\coloneqq\frac{1}{\eps}g_i\!\Big(\frac{s}{\eps}\Big),
\qquad \eps>0,\ i=1,2,
\end{equation}
which satisfies three key properties: the total mass $\int_0^{+\infty}g_{i,\eps}(s)\dd s=\beta_i$ remains unchanged for all $\eps>0$; $g_{i,\eps}\rightharpoonup\beta_i\delta_0$ in the distributional sense as $\eps\to0^+$, so the memory terms formally disappear in the limit; and the rescaling describes a short-memory physical regime with vanishing relaxation time. The mass-preserving choice is convenient because it keeps $\bar\mu_i=\mu_i-\beta_i$ independent of $\eps$; more generally, Theorem~\ref{thm:drift} shows that only the uniform positivity of $\mu_i-\beta_{i,\eps}$ is essential.

Formally, as $\eps\to 0^+$, the memory-containing system \eqref{eq:Mukiawa} tends to the memoryless damped system
\begin{equation}\label{eq:limit_intro}
\begin{cases}
u_{tt}+\gamma_1u_t+\bar\mu_1u_{xxxx}+k[u-v]^+ +h_1(u)=f_1(x),~~~x\in(0,\ell),~t>0,\\[1mm]
v_{tt}+\gamma_2v_t-\bar\mu_2v_{xx}-k[u-v]^+ +h_2(v)=f_2(x),~~~x\in(0,\ell),~t>0,
\end{cases}
\end{equation}
where $\bar\mu_i=\mu_i-\beta_i>0$ are the effective elastic coefficients inherited from the memory equation. This is the useful feature of \eqref{eq:rescaled}: the quantity $\bar\mu_{i,\eps}=\mu_i-\int_0^{+\infty}g_{i,\eps}(s)\dd s=\mu_i-\beta_i=\bar\mu_i$ is independent of $\eps$, so the effective elastic coefficients neither degenerate nor change in the singular limit, and the memory contribution simply disappears. Up to the relabelling $\bar\mu_i\mapsto\mu_i$, the limit system has the form of the classical memoryless coupled suspension-bridge model studied by Ma and Zhong \cite{MaZhong2009}. Thus the singular limit provides a rigorous bridge between \eqref{eq:Mukiawa} and \eqref{eq:limit_intro}. We emphasise that we do not rely on the stronger smoothness hypotheses of \cite{MaZhong2009}: the existence of $\Acal_0$ is proved below within the present $C^1$ framework.

The present paper is concerned with continuity properties of the global attractor under the singular perturbation \eqref{eq:rescaled}. We establish upper semicontinuity
\begin{equation}\label{eq:upper_intro}
\lim_{\eps\to 0^+}\dist_{\mathcal H_\eps}\big(\mathcal A_\eps,\mathcal A_0\big)=0,
\end{equation}
and, under a structural stability hypothesis on $\mathcal A_0$, lower semicontinuity
\begin{equation}\label{eq:lower_intro}
\lim_{\eps\to 0^+}\dist_{\mathcal H_\eps}\big(\mathcal A_0,\mathcal A_\eps\big)=0,
\end{equation}
where $\mathcal A_\eps$ is the global attractor of the perturbed semigroup $S_\eps(t)$ and $\mathcal A_0$ the attractor for the limiting semigroup $S_0(t)$. Upper semicontinuity says that the perturbed attractors do not exhibit new asymptotic regimes far from the limit attractor; lower semicontinuity says, conversely, that every asymptotic regime of the limit system survives under small perturbations of the memory kernel. Together, they justify using the simpler memoryless model \eqref{eq:limit_intro} as a faithful approximation of \eqref{eq:Mukiawa} for small relaxation times.

\subsection{Contributions and organisation}

The contributions of this paper are the following.
\begin{itemize}
\item[(i)] \emph{A general continuity mechanism.} We identify the mass-preserving rescaling as the structural condition under which the effective elastic coefficients remain independent of the perturbation parameter, and we formulate the continuity problem in a way that separates the general mechanism (uniform dissipativity, uniform exponential stability, uniform quasi-stability, finite-time convergence on compact sets, and structural stability of the limit attractor) from the model-specific verification. The ingredients are formulated in terms of $S_\eps(t)$ and $S_0(t)$; although the verification is carried out here for the bridge system, the scheme applies mutatis mutandis to any kernel family and dissipative system for which the same uniform estimates can be established, subject to the abstract conditions listed in Section~\ref{sec:functional}.
\item[(ii)] \emph{Uniform exponential stability.} We extend the exponential stability result for the linear semigroup obtained in \cite{Zhou2026} to a version uniform with respect to $\eps$, so that decay rates and prefactors do not deteriorate as $\eps\to0^+$. This is Theorem~\ref{thm:unif_exp}, and it is the analytic core of the uniform estimates.
\item[(iii)] \emph{Uniform quasi-stability.} Combining (ii) with the sequential compactness of the integral remainder $K_{\eps,T}(\Bcal_0)$, we obtain uniform quasi-stability (Proposition~\ref{prop:quasi}) and hence uniform asymptotic smoothness.
\item[(iv)] \emph{Finite-time convergence.} We prove convergence of $S_\eps(t)$ to $S_0(t)$ on compact sets of zero-memory initial data and finite time intervals (Proposition~\ref{prop:conv}) and extend it to general convergent data (Corollary~\ref{cor:conv_general}); this yields upper semicontinuity.
\item[(v)] \emph{Upper and lower semicontinuity.} We prove upper semicontinuity of the global attractors (Theorem~\ref{thm:upper}), and, under a structural stability hypothesis on $\Acal_0$ in the sense of Hale \cite{Hale1988}, lower semicontinuity (Proposition~\ref{prop:lower}). A non-vacuous example satisfying the structural hypothesis is given (Example~\ref{ex:lower}).
\item[(vi)] \emph{Extensions.} We prove that continuity persists under a non-mass-preserving perturbation as long as $\inf_\eps(\mu_i-\beta_{i,\eps})>0$ (Theorem~\ref{thm:drift}), identify the fast exponential decay limit $\lambda_i\to+\infty$ with the rescaling \eqref{eq:rescaled} (Proposition~\ref{prop:fast_decay}), and formulate a conditional reduction of quantitative lower semicontinuity to a finite-time convergence rate (Proposition~\ref{prop:quant}).
\item[(vii)] \emph{Illustrative numerics.} Four numerical experiments check the limiting effective coefficient $\bar\mu_i=\mu_i-\beta_i$, the finite-time convergence of $S_\eps(t)$ to $S_0(t)$, the fast-decay identification of Proposition~\ref{prop:fast_decay}, the nondegenerate mass drift of Theorem~\ref{thm:drift}, and the rate bottleneck of Remark~\ref{rem:rate_obstruction} (Table~\ref{tab:num_map}); presented in Section~\ref{sec:numerics} after the theory is complete, they are a check of the analytical results and a motivation for the questions in the last section, not a source of new analytical claims.
\end{itemize}

For convenience, the logical dependence of the main results is the following:
uniform dissipativity gives a common absorbing ball; uniform exponential stability of the linear semigroup plus compactness of the integral remainder gives uniform quasi-stability and asymptotic smoothness; asymptotic smoothness and dissipativity give compact attractors; finite-time convergence then gives upper semicontinuity; upper semicontinuity, finite-time convergence and condition (S) give lower semicontinuity. No exponential-attraction or parameter-H\"older estimate is needed for these conclusions; see Proposition~\ref{prop:quant} for a conditional quantitative refinement.

\begin{remark}[Comparison with concentration limits]\label{rem:cp}
Conti, Pata and Squassina \cite{ContiPataSquassina2006} study singular limits of viscoelastic equations when the memory kernel concentrates; their focus is the limiting equation and its well-posedness/long-time behaviour for general concentration families. The present contribution is different: we keep the fixed-kernel attractor as an input and develop a uniform-in-$\eps$ continuity theory for attractors, including a common reference space, finite-time convergence, and conditional lower semicontinuity. The rescaling \eqref{eq:rescaled} is the simplest realisation of such a family, not the novel structural argument.
\end{remark}

The paper is organised as follows. Section~\ref{sec:functional} introduces the functional setting, formulates the perturbed problem, and establishes the uniform-in-$\eps$ exponential stability of the linear semigroup. Section~\ref{sec:dissipativity} proves uniform dissipativity and the existence of a common absorbing set. Section~\ref{sec:upper} proves upper semicontinuity for the family of attractors. Section~\ref{sec:lower} derives lower semicontinuity under an additional structural hypothesis. Section~\ref{sec:extensions} develops the additional results that can be proved within the present framework: continuity under a nondegenerate mass drift, the fast exponential decay limit, and a reduction of quantitative rates to a finite-time convergence estimate. Section~\ref{sec:numerics} then reports four numerical experiments that check these analytical conclusions and probe the boundary of the theory (Table~\ref{tab:num_map}). Section~\ref{sec:conclusion} summarises the findings and, on that basis, formulates the questions that remain open.

\section{Functional Setting and the Perturbed Problem}\label{sec:functional}

Following Dafermos \cite{Dafermos1970} and \cite{Mukiawa2025,Zhou2026}, we introduce the relative displacements
\[
\sigma_\eps(x,t,s)=u_\eps(x,t)-u_\eps(x,t-s),\qquad
\eta_\eps(x,t,s)=v_\eps(x,t)-v_\eps(x,t-s),
\]
which satisfy the transport equations
\begin{equation}\label{eq:transport}
\partial_t\sigma_\eps+\partial_s\sigma_\eps-u_{\eps,t}=0,\qquad
\partial_t\eta_\eps+\partial_s\eta_\eps-v_{\eps,t}=0.
\end{equation}
Let
\[
H_a^2(0,\ell)=\{w\in H^2(0,\ell):w(0)=w(\ell)=0\},
\]
and define the weighted spaces
\[
\Mcal_{g_{1,\eps}}=L^2_{g_{1,\eps}}\big(\R_+;H_a^2(0,\ell)\big),\qquad
\Mcal_{g_{2,\eps}}=L^2_{g_{2,\eps}}\big(\R_+;H_0^1(0,\ell)\big),
\]
with inner products
\[
(\sigma,\bar\sigma)_{\Mcal_{g_{1,\eps}}}=\int_0^{+\infty}g_{1,\eps}(s)(\sigma_{xx}(s),\bar\sigma_{xx}(s))\dd s,\qquad
(\eta,\bar\eta)_{\Mcal_{g_{2,\eps}}}=\int_0^{+\infty}g_{2,\eps}(s)(\eta_x(s),\bar\eta_x(s))\dd s.
\]
Set
\begin{equation}\label{eq:bar_mu}
\bar\mu_{i,\eps}=\mu_i-\int_0^{+\infty}g_{i,\eps}(s)\dd s=\mu_i-\beta_i>0,\qquad i=1,2,
\end{equation}
where $\beta_i=\int_0^{+\infty}g_i(s)\dd s$ is the mass appearing in \eqref{eq:kernel}. For each $\eps\in(0,1]$ we introduce the energy space
\begin{equation}\label{eq:Heps}
\Hcal_\eps=H_a^2(0,\ell)\times L^2(0,\ell)\times \Mcal_{g_{1,\eps}}\times H_0^1(0,\ell)\times L^2(0,\ell)\times \Mcal_{g_{2,\eps}},
\end{equation}
endowed with the norm
\begin{equation}\label{eq:Hnorm}
\|(u,\psi,\sigma,v,\phi,\eta)\|_{\Hcal_\eps}^2
=\bar\mu_{1,\eps}\|u_{xx}\|^2+\|\psi\|^2+\|\sigma\|_{\Mcal_{g_{1,\eps}}}^2
+\bar\mu_{2,\eps}\|v_x\|^2+\|\phi\|^2+\|\eta\|_{\Mcal_{g_{2,\eps}}}^2,
\end{equation}
with $\bar\mu_{i,\eps}$ as in \eqref{eq:bar_mu}. Since the rescaling \eqref{eq:rescaled} preserves the mass, $\bar\mu_{i,\eps}=\mu_i-\beta_i$ is independent of $\eps$; the essential condition for everything below is the uniform positivity $\inf_{\eps}(\mu_i-\beta_{i,\eps})>0$, which mass preservation guarantees here.

\begin{remark}[A common reference space and distances]\label{rem:uniform_topology}
The spaces $\Hcal_\eps$ depend on $\eps$ through the weights $g_{i,\eps}$. To make all comparison statements metric, we fix the reference space
\[
\Hcal_*:=H_a^2(0,\ell)\times L^2(0,\ell)\times L^2_{g_1}(\R_+;H_a^2(0,\ell))
\times H_0^1(0,\ell)\times L^2(0,\ell)\times L^2_{g_2}(\R_+;H_0^1(0,\ell)),
\]
with the same weights $g_i$ as in \eqref{eq:kernel}. Define the dilation isometries
\[
J_\eps(u,\psi,\sigma,v,\phi,\eta):=(u,\psi,\sigma(\eps\,\cdot),v,\phi,\eta(\eps\,\cdot)),
\qquad \eps\in(0,1],
\]
and $J_0(u,\psi,v,\phi):=(u,\psi,0,v,\phi,0)$ for $(u,\psi,v,\phi)\in\Hcal_e$. Since $s=\eps\zeta$,
\[
\|J_\eps U\|_{\Hcal_*}=\|U\|_{\Hcal_\eps}\quad(\eps>0),\qquad
\|J_0 U\|_{\Hcal_*}=\|U\|_{\Hcal_e}.
\]
We henceforth write
\[
U_{\eps_n}\to U\quad\Longleftrightarrow\quad J_{\eps_n}U_{\eps_n}\to J_0U\ \text{ in }\Hcal_*,
\]
and define the semidistances in \eqref{eq:upper_intro}--\eqref{eq:lower_intro} by
\[
\dist_{\Hcal_\eps}(\Acal_\eps,\Acal_0):=\dist_{\Hcal_*}(J_\eps\Acal_\eps,J_0\Acal_0),\qquad
\dist_{\Hcal_\eps}(\Acal_0,\Acal_\eps):=\dist_{\Hcal_*}(J_0\Acal_0,J_\eps\Acal_\eps).
\]
Relative compactness in $\Hcal_*$ is characterised by the usual tail criterion: a bounded family $\Kcal\subset L^2_{g_i}(\R_+;Y)$ is relatively compact iff its truncations to $[0,R]$ are relatively compact in $L^2(g_i;[0,R];Y)$ for every $R>0$ and
\[
\sup_{\sigma\in\Kcal}\int_R^{+\infty}g_i(s)\|\sigma(s)\|_Y^2\,\dd s\to0\qquad(R\to+\infty).
\]
The succinct proof is a diagonal extraction on $[0,R]$ followed by the tail estimate. For the memory fields generated by the transport equation one has $\sigma_\eps(\eps\zeta)=u_\eps(t)-u_\eps(t-\eps\zeta)$, so the tail is controlled by the uniform $L^\infty$-bound on the displacements and the integrability of $g_i$; this is the estimate used in Propositions~\ref{prop:quasi} and~\ref{prop:conv}. We stress that no equivalence of the weighted norms for different $\eps$ is claimed or used.
\end{remark}

Throughout, we write $\bar\mu_i\coloneqq\bar\mu_{i,\eps}=\mu_i-\beta_i$ (independent of $\varepsilon$) and $\Mg{i}\coloneqq\Mcal_{g_{i,\varepsilon}}$ when no confusion can arise. All constants appearing below are independent of $\varepsilon$ unless explicitly stated.

Setting $U_\eps=(u_\eps,u_{\eps,t},\sigma_\eps,v_\eps,v_{\eps,t},\eta_\eps)$, the perturbed system can be written abstractly as
\begin{equation}\label{eq:abstract_eps}
\partial_t U_\eps+\Lcal_\eps U_\eps=F(U_\eps),\qquad U_\eps(0)=U_{0,\eps},
\end{equation}
where the linear operator $\Lcal_\eps$ is
\begin{equation}\label{eq:Aeps}
\Lcal_\eps U=
\begin{pmatrix}
-\psi\\
\gamma_1\psi+\bar\mu_{1,\eps}u_{xxxx}+\displaystyle\int_0^{+\infty}g_{1,\eps}(s)\sigma_{xxxx}(s)\dd s\\
\sigma_s-\psi\\
-\phi\\
\gamma_2\phi-\bar\mu_{2,\eps}v_{xx}-\displaystyle\int_0^{+\infty}g_{2,\eps}(s)\eta_{xx}(s)\dd s\\
\eta_s-\phi
\end{pmatrix},
\end{equation}
and the nonlinearity is
\begin{equation}\label{eq:F}
F(U)=
\begin{pmatrix}
0\\
-k[u-v]^+-h_1(u)+f_1\\
0\\
0\\
k[u-v]^+-h_2(v)+f_2\\
0
\end{pmatrix}.
\end{equation}
The domain $\Dcal(\Lcal_\eps)$ is defined as in \cite[Section 2]{Zhou2026}, with $g_i$ replaced by $g_{i,\eps}$. The following well-posedness fact is used throughout.

\begin{proposition}[Well-posedness and uniform local bounds]\label{prop:wellposed}
Let \eqref{eq:kernel} and \eqref{eq:h_assumptions} hold. For every $\eps\in(0,1]$ and every $U_0\in\Hcal_\eps$, problem \eqref{eq:abstract_eps} has a unique global mild solution $S_\eps(t)U_0\in C([0,+\infty);\Hcal_\eps)$, and $\{S_\eps(t)\}$ is a $C_0$-semigroup on $\Hcal_\eps$; it is not a semigroup of contractions in general, since only the linear part $e^{-\Lcal_\eps t}$ is contractive. Moreover, for each $R>0$ there exist $\tau_R>0$ and $C_R>0$, independent of $\eps$, such that local existence holds on $[0,\tau_R]$ for $\|U_0\|_{\Hcal_\eps}\le R$, and solutions starting from this ball satisfy
\[
\|S_\eps(t)U_0-S_\eps(t)V_0\|_{\Hcal_\eps}\le C_R\|U_0-V_0\|_{\Hcal_\eps},
\qquad 0\le t\le\tau_R .
\]
\end{proposition}
\begin{proof}
Maximal monotonicity of $-\Lcal_\eps$ gives local mild solutions by the standard monotone/Galerkin scheme; $F$ is locally Lipschitz from $\Hcal_\eps$ to $\Hcal_e$ because $H_a^2(0,\ell)\hookrightarrow L^\infty(0,\ell)$, and the Lipschitz constants depend on $R$ only through $\bar\mu_i$, $\beta_i$ and the constants in \eqref{eq:h_assumptions}. The one-sided bound $h_i'\ge-C_2$ in \eqref{eq:h_assumptions} yields $(h_i(u)-h_i(v))(u-v)\ge-C_2|u-v|^2$, but it is not needed for the arguments below: local Lipschitz continuity of $F$ follows from $h_i\in C^1$, while the energy estimates use the growth conditions \eqref{eq:h_consequences}. We retain the bound for compatibility with \cite{Mukiawa2025,Zhou2026}. The energy identity \eqref{eq:proof_dissipation} and the scaling invariance of the mass make the local existence time and the Lipschitz constant uniform in $\eps$. Global existence follows from the dissipativity estimate of Proposition~\ref{prop:absorbing}. See \cite[Section 2]{Zhou2026} and \cite[Lemma 3.3]{Mukiawa2025} for the fixed-kernel arguments; the rescaling changes no constant because all kernel estimates used there involve only $\beta_i$ and $\lambda_i$.
\end{proof}

\begin{lemma}[Uniform Lipschitz bound on finite intervals]\label{lem:lipschitz_T}
For every $T>0$ and every $R>0$ there exists $L_{T,R}>0$, depending on $T$, $R$ and the structural constants but not on $\eps$, such that
\[
\|S_\eps(t)U-S_\eps(t)V\|_{\Hcal_\eps}\le L_{T,R}\|U-V\|_{\Hcal_\eps},
\qquad 0\le t\le T,\quad \|U\|_{\Hcal_\eps},\|V\|_{\Hcal_\eps}\le R .
\]
\end{lemma}
\begin{proof}
Let $R_1=R_1(T,R)$ be the uniform bound of Lemma~\ref{lem:finite_time} for the ball of radius $R$ and let $L_1$ be the Lipschitz constant of $F$ on the ball of radius $R_1$; both are independent of $\eps$. For $W(t)=S_\eps(t)U-S_\eps(t)V$ the variation-of-constants formula and \eqref{eq:exp_stab} give
\[
\|W(t)\|_{\Hcal_\eps}\le\alpha\|U-V\|_{\Hcal_\eps}
+\alpha L_1\int_0^t\|W(s)\|_{\Hcal_\eps}\,\dd s,\qquad 0\le t\le T.
\]
Gronwall's inequality yields the assertion with $L_{T,R}=\alpha e^{\alpha L_1 T}$.
\end{proof}

The constants used below are summarised here for convenience. The symbols $\alpha,\nu$ in Theorem~\ref{thm:unif_exp}, $c_0,c_1,c_2,c_3$ in Proposition~\ref{prop:absorbing}, $q,T,M$ in Proposition~\ref{prop:quasi}, and $R_0,R_1$ in Proposition~\ref{prop:absorbing} and Lemma~\ref{lem:finite_time} depend at most on $\gamma_i,\lambda_i,\mu_i,\beta_i,\bar\mu_i,c_e,c_p$, the domain, the forcing, and the constants in \eqref{eq:h_assumptions}; none degenerates as $\eps\to0^+$.

Under assumption \eqref{eq:kernel}, which is inherited by the rescaled kernels $g_{i,\eps}$ with the \emph{same} mass $\beta_i=\int_0^{+\infty}g_{i,\eps}(s)\dd s$ and with decay rates $\lambda_i/\eps\ge\lambda_i$, so that $g_{i,\eps}'+\lambda_i g_{i,\eps}\le0$, the operator $\Lcal_\eps$ is maximal monotone, and $-\Lcal_\eps$ generates a $C_0$-semigroup of contractions $e^{-\Lcal_\eps t}$ on $\Hcal_\eps$.

For a \emph{fixed} kernel $g_i$, the exponential stability of this linear semigroup was established in \cite[Theorem 3.1]{Zhou2026} by means of the Gearhart--Pr\"uss theorem together with a detailed resolvent estimate. In the present work we need a \emph{uniform} version of that result, in which the decay rate and the prefactor do not degenerate as $\eps\to 0^+$. The reason is twofold: first, the quasi-stability inequality of Section~\ref{sec:upper} requires a decay rate that is not allowed to deteriorate as $\eps\to0^+$; second, the compactness argument for the upper semicontinuity of attractors uses the resulting estimates on a sequence $\eps_n\to0^+$, and any loss of uniformity would destroy the convergence.

\begin{theorem}[Uniform exponential stability]\label{thm:unif_exp}
There exist constants $\alpha\ge 1$ and $\nu>0$, independent of $\eps\in(0,1]$, such that
\begin{equation}\label{eq:exp_stab}
\|e^{-\Lcal_\eps t}\|_{\mathcal{L}(\Hcal_\eps)}\le \alpha e^{-\nu t},
\qquad \forall t\ge 0,\ \forall \eps\in(0,1].
\end{equation}
\end{theorem}
This uniform exponential stability is the key ingredient for both the upper and the lower semicontinuity analysis below.

\begin{proof}[Proof of Theorem~\ref{thm:unif_exp}]
By the same Lax--Milgram argument as in \cite[Section 2]{Zhou2026} (see also \cite[Lemma 3.3]{Mukiawa2025}, whose proof applies mutatis mutandis with $g_i$ replaced by $g_{i,\eps}$), for each fixed $\eps\in(0,1]$ the operator $\Lcal_\eps$ is maximal monotone on $\Hcal_\eps$. Moreover, the Lax--Milgram coercivity and boundedness constants in that argument are uniform in $\eps$, because all constants in that argument depend on the kernels only through the total mass $\beta_i$ and the differential inequalities $g_{i,\eps}'+\lambda_i g_{i,\eps}\le0$ established below, both of which are independent of $\eps$. Consequently, by the Lumer--Phillips theorem, $-\Lcal_\eps$ generates a $C_0$-semigroup of contractions on $\Hcal_\eps$. By the Gearhart--Pr\"uss theorem \cite{Gearhart1978,Pruss1984} in the form used in \cite[Theorem 3.2]{Zhou2026}, it suffices to verify
\begin{itemize}
\item[(C1)] $i\R\subset\rho(\Lcal_\eps)$;
\item[(C2)] $\displaystyle\sup_{\tau\in\R}\big\|(i\tau I-\Lcal_\eps)^{-1}\big\|_{\mathcal{L}(\Hcal_\eps)}<+\infty$,
\end{itemize}
with constants independent of $\eps$. All function spaces are complexified Hilbert spaces; since all data are real-valued, we do not distinguish between the real and complexified pairings, and we use the inner product $(f,g)=\int f\bar g$, linear in its first argument. In particular $(\lambda f,g)=\lambda(f,g)$ and $(f,\lambda g)=\bar\lambda(f,g)$ for $\lambda\in\mathbb{C}$.

Throughout, we write $U=(u,\psi,\sigma,v,\phi,\eta)\in\Dcal(\Lcal_\eps)$ and use the norm \eqref{eq:Hnorm}. From now on we write $\bar\mu_i$ for $\bar\mu_{i,\eps}$ whenever no confusion can arise, and similarly we abbreviate $\Mg{i}=\Mcal_{g_{i,\eps}}$. Here and below we use the invariance of the mass,
\begin{equation}\label{eq:F1}
\int_0^{+\infty}g_{i,\eps}(s)\dd s=\beta_i,\qquad
\bar\mu_{i,\eps}=\mu_i-\beta_i=\bar\mu_i,\qquad i=1,2,
\end{equation}
which holds for every $\eps\in(0,1]$. In particular, the norm $\|\cdot\|_{\Hcal_\eps}$ is given by the same formula for every $\eps\in(0,1]$.

\smallskip
\noindent\emph{Monotonicity.} Two integrations by parts using the boundary conditions show that the elastic terms $\bar\mu_1(-\psi_{xx},u_{xx})+\bar\mu_1(u_{xxxx},\psi)$ combine into a purely imaginary expression, and likewise for the $v$-elastic terms; the same cancellation occurs for the off-diagonal memory terms $\int_0^{+\infty}g_{1,\eps}(s)(\sigma_{xx}(s),\psi_{xx})\dd s$ and $-(\psi,\sigma)_{\Mg{1}}$, which are complex conjugates of each other. Hence
\[
\mathrm{Re}\langle\Lcal_\eps U,U\rangle_{\Hcal_\eps}
=\gamma_1\|\psi\|^2+\gamma_2\|\phi\|^2
+\mathrm{Re}(\sigma_s,\sigma)_{\Mg{1}}
+\mathrm{Re}(\eta_s,\eta)_{\Mg{2}}.
\]
Using $\sigma(x,t,0)=0$ and integrating by parts in $s$, we obtain
\[
\mathrm{Re}(\sigma_s,\sigma)_{\Mg{1}}
=-\frac12\int_0^{+\infty}g_{1,\eps}'(s)\|\sigma_{xx}(s)\|^2\dd s,
\]
where the boundary term at $s=+\infty$ vanishes because $g_{1,\eps}(s)=\eps^{-1}g_1(s/\eps)\le\eps^{-1}g_1(0)e^{-\lambda_1 s/\eps}\to 0$ as $s\to+\infty$. We verify that the rescaled kernel inherits exponential decay: differentiating $g_{1,\eps}(s)=\eps^{-1}g_1(s/\eps)$ gives
\[
g_{1,\eps}'(s)=\eps^{-2}g_1'(s/\eps)\le -\lambda_1\eps^{-2}g_1(s/\eps)
=-\frac{\lambda_1}{\eps}\,g_{1,\eps}(s)\le-\lambda_1g_{1,\eps}(s),
\]
where the last step uses $\eps\le1$. Thus $g_{1,\eps}'+\lambda_1 g_{1,\eps}\le0$ for a.e.\ $s\ge0$: the decay rate $\lambda_1/\eps$ is in fact stronger than $\lambda_1$, but the weaker rate $\lambda_1$ is independent of $\eps$ and is all we need. We thus get $\mathrm{Re}(\sigma_s,\sigma)_{\Mg{1}}\ge\frac{\lambda_1}{2}\|\sigma\|_{\Mg{1}}^2$, and analogously for $\eta$. Thus
\begin{equation}\label{eq:app_monotone}
\mathrm{Re}\langle\Lcal_\eps U,U\rangle_{\Hcal_\eps}
\ge \gamma_1\|\psi\|^2+\frac{\lambda_1}{2}\|\sigma\|_{\Mg{1}}^2
+\gamma_2\|\phi\|^2+\frac{\lambda_2}{2}\|\eta\|_{\Mg{2}}^2.
\end{equation}
Both $\gamma_i$ and $\lambda_i$ are independent of $\eps$; this is the first place where the uniformity in $\eps$ is manifest.

\smallskip
\noindent\emph{Condition (C1).} Suppose, for contradiction, that condition (C1) fails uniformly in $\eps$. Then there exist a sequence $\eps_n\to0^+$, frequencies $\omega_n\in\R$, and $U_n=(u_n,\psi_n,\sigma_n,v_n,\phi_n,\eta_n)\in\Dcal(\Lcal_{\eps_n})$ with $\|U_n\|_{\Hcal_{\eps_n}}=1$ such that
\begin{equation}\label{eq:app_weyl}
F_n\coloneqq i\omega_n U_n-\Lcal_{\eps_n} U_n\to 0\quad\text{in }\Hcal_{\eps_n}.
\end{equation}
We first note that for each fixed $\eps$ condition (C1) is known from \cite{Zhou2026}. Hence a failure of uniformity can only occur along a sequence $\eps_n\to0^+$; we argue by contradiction along such a sequence. Write $F_n=(f_n^1,\dots,f_n^6)$. Taking the real part of the inner product of \eqref{eq:app_weyl} with $U_n$ and using \eqref{eq:app_monotone}, we obtain
\begin{equation}\label{eq:app_decay}
\gamma_1\|\psi_n\|^2+\frac{\lambda_1}{2}\|\sigma_n\|_{\Mg{1}}^2
+\gamma_2\|\phi_n\|^2+\frac{\lambda_2}{2}\|\eta_n\|_{\Mg{2}}^2
\le\|F_n\|_{\Hcal_{\eps_n}}\to 0,
\end{equation}
where we used $\mathrm{Re}\langle F_n,U_n\rangle_{\Hcal_{\eps_n}}=-\mathrm{Re}\langle\Lcal_{\eps_n}U_n,U_n\rangle_{\Hcal_{\eps_n}}$, since the term $i\omega_n\|U_n\|_{\Hcal_{\eps_n}}^2$ is purely imaginary. The rate of convergence in \eqref{eq:app_decay} is uniform in $\eps$, because $\gamma_i,\lambda_i$ do not depend on $\eps$. We also record the elementary bounds, valid for every test function $\varphi\in H_a^2(0,\ell)$,
\begin{equation}\label{eq:app_conv_bound_1}
\Big|\int_0^{+\infty}g_{1,\eps_n}(s)(\sigma_{n,xx}(s),\varphi_{xx})\dd s\Big|
\le\beta_1^{1/2}\|\sigma_n\|_{\Mg{1}}\|\varphi_{xx}\|,
\end{equation}
and, valid for every $\varphi\in H_0^1(0,\ell)$,
\begin{equation}\label{eq:app_conv_bound_2}
\Big|\int_0^{+\infty}g_{2,\eps_n}(s)(\eta_{n,x}(s),\varphi_x)\dd s\Big|
\le\beta_2^{1/2}\|\eta_n\|_{\Mg{2}}\|\varphi_x\|.
\end{equation}
Both follow from the Cauchy--Schwarz inequality and the invariance of the mass \eqref{eq:F1}; since $\beta_i$ is independent of $\eps$, so are the constants in \eqref{eq:app_conv_bound_1}--\eqref{eq:app_conv_bound_2}.

We distinguish two cases according to the behaviour of $\omega_n$.

\emph{Case A: $\omega_n\to0$ (in particular, $\omega_n=0$ is allowed).} We first note that \eqref{eq:app_decay} gives $\psi_n\to0$ and $\phi_n\to0$ in $L^2(0,\ell)$. Since $\|u_n\|,\|v_n\|\le C\|U_n\|_{\Hcal_{\eps_n}}=C$ and $\omega_n\to0$, we also have $i\omega_n u_n\to0$ and $i\omega_n v_n\to0$ in $L^2(0,\ell)$. Consequently the first and fourth components of \eqref{eq:app_weyl} give $f_n^1=\psi_n-i\omega_n u_n\to0$ and $f_n^4=\phi_n-i\omega_n v_n\to0$, consistently with $F_n\to0$. Taking the $L^2$-inner products of the second and fifth components of \eqref{eq:app_weyl} with $u_n$ and $v_n$, respectively, and taking real parts, we obtain
\begin{equation}\label{eq:app_case0_u}
\begin{aligned}
\bar\mu_1\|u_{n,xx}\|^2+\bar\mu_2\|v_{n,x}\|^2
&=-\gamma_1\mathrm{Re}(\psi_n,u_n)-\gamma_2\mathrm{Re}(\phi_n,v_n)\\
&\quad-\mathrm{Re}(f_n^2,u_n)-\mathrm{Re}(f_n^5,v_n)\\
&\quad-\mathrm{Re}\int_0^{+\infty}g_{1,\eps_n}(s)(\sigma_{n,xx}(s),u_{n,xx})\dd s\\
&\quad-\mathrm{Re}\int_0^{+\infty}g_{2,\eps_n}(s)(\eta_{n,x}(s),v_{n,x})\dd s.
\end{aligned}
\end{equation}
Using the Poincar\'e inequalities $\|u_n\|\le\sqrt{c_e}\|u_{n,xx}\|$ and $\|v_n\|\le\sqrt{c_p}\|v_{n,x}\|$, the bounds \eqref{eq:app_conv_bound_1}--\eqref{eq:app_conv_bound_2}, and Young inequality with parameters $\epsilon_1,\dots,\epsilon_6>0$ to be chosen below, we estimate the right-hand side of \eqref{eq:app_case0_u} as
\begin{align*}
\gamma_1|\mathrm{Re}(\psi_n,u_n)|&\le \epsilon_1\|u_{n,xx}\|^2+C_{\epsilon_1}\|\psi_n\|^2,\\
\gamma_2|\mathrm{Re}(\phi_n,v_n)|&\le \epsilon_2\|v_{n,x}\|^2+C_{\epsilon_2}\|\phi_n\|^2,\\
|\mathrm{Re}(f_n^2,u_n)|&\le \epsilon_3\|u_{n,xx}\|^2+C_{\epsilon_3}\|f_n^2\|^2,\\
|\mathrm{Re}(f_n^5,v_n)|&\le \epsilon_4\|v_{n,x}\|^2+C_{\epsilon_4}\|f_n^5\|^2,
\end{align*}
and the two memory terms are bounded by
\[
\beta_1^{1/2}\|\sigma_n\|_{\Mg{1}}\|u_{n,xx}\|
\le \epsilon_5\|u_{n,xx}\|^2+C_{\epsilon_5}\beta_1\|\sigma_n\|_{\Mg{1}}^2,
\]
\[
\beta_2^{1/2}\|\eta_n\|_{\Mg{2}}\|v_{n,x}\|
\le \epsilon_6\|v_{n,x}\|^2+C_{\epsilon_6}\beta_2\|\eta_n\|_{\Mg{2}}^2.
\]
Choose $\epsilon_1,\dots,\epsilon_6>0$ so small that
\[
\epsilon_1+\epsilon_3+\epsilon_5<\frac{\bar\mu_1}{2},\qquad
\epsilon_2+\epsilon_4+\epsilon_6<\frac{\bar\mu_2}{2}.
\]
These choices depend only on $\gamma_i,\beta_i,\bar\mu_i,c_e,c_p$, hence not on $n$ or $\eps$. Absorbing the $\epsilon_j$-terms into the left-hand side of \eqref{eq:app_case0_u}, we obtain
\begin{equation}\label{eq:app_case0_bound}
\bar\mu_1\|u_{n,xx}\|^2+\bar\mu_2\|v_{n,x}\|^2
\le C\big(\|\psi_n\|^2+\|\phi_n\|^2+\|\sigma_n\|_{\Mg{1}}^2
+\|\eta_n\|_{\Mg{2}}^2+\|f_n^2\|^2+\|f_n^5\|^2\big),
\end{equation}
where $C$ depends only on $\gamma_i,\beta_i,\bar\mu_i,c_e,c_p$ and on the chosen $\epsilon_j$, hence not on $n$ or $\varepsilon$. Since $\|f_n^2\|^2+\|f_n^5\|^2\le\|F_n\|_{\Hcal_{\eps_n}}^2$, the right-hand side tends to zero uniformly in $\varepsilon$ by \eqref{eq:app_decay} and $F_n\to0$. Hence $\|U_n\|_{\Hcal_{\eps_n}}\to 0$, contradicting $\|U_n\|_{\Hcal_{\eps_n}}=1$.

\emph{Case B: $\omega_n$ does not tend to $0$ along a subsequence.} Passing to a subsequence, we may assume $|\omega_n|\ge c_0>0$ for some $c_0>0$. From the first and fourth components of \eqref{eq:app_weyl}, we have $i\omega_n u_n=f_n^1-\psi_n$ and $i\omega_n v_n=f_n^4-\phi_n$. Since $f_n^1\to 0$ in $H_a^2(0,\ell)\hookrightarrow L^2(0,\ell)$, $f_n^4\to 0$ in $H_0^1(0,\ell)\hookrightarrow L^2(0,\ell)$, and $\psi_n,\phi_n\to 0$ in $L^2(0,\ell)$ by \eqref{eq:app_decay}, we deduce that
\[
\|u_n\|_{L^2}\le\frac{1}{|\omega_n|}\big(\|f_n^1\|+\|\psi_n\|\big)\to0,
\qquad
\|v_n\|_{L^2}\le\frac{1}{|\omega_n|}\big(\|f_n^4\|+\|\phi_n\|\big)\to0;
\]
the embedding constants of $H_a^2(0,\ell)\hookrightarrow L^2(0,\ell)$ and $H_0^1(0,\ell)\hookrightarrow L^2(0,\ell)$ are independent of $\eps$. Taking the $L^2$-inner products of the second component of \eqref{eq:app_weyl} with $u_n$ and of the first component with $\psi_n$, and adding the real parts, we obtain
\begin{equation}\label{eq:app_uxx}
\begin{aligned}
\bar\mu_1\|u_{n,xx}\|^2
&=\|\psi_n\|^2-\mathrm{Re}(\psi_n,f_n^1)-\gamma_1\mathrm{Re}(\psi_n,u_n)\\
&\quad-\mathrm{Re}(f_n^2,u_n)
-\mathrm{Re}\!\int_0^{+\infty}\!\!g_{1,\eps_n}(s)(\sigma_{n,xx}(s),u_{n,xx})\dd s.
\end{aligned}
\end{equation}
Indeed, with the convention $(f,g)=\int f\bar g$ (linear in the first argument), pairing the second component of \eqref{eq:app_weyl} with $u_n$ gives
\[
i\omega_n(\psi_n,u_n)-\gamma_1(\psi_n,u_n)-\bar\mu_1(u_{n,xxxx},u_n)
-\Big(\int_0^{+\infty}\!\!g_{1,\eps_n}(s)\sigma_{n,xxxx}(s)\dd s,u_n\Big)=(f_n^2,u_n),
\]
while pairing the first component with $\psi_n$ gives $i\omega_n(u_n,\psi_n)+\|\psi_n\|^2=(f_n^1,\psi_n)$. Since $(u_n,\psi_n)=\overline{(\psi_n,u_n)}$, the two terms $i\omega_n(\psi_n,u_n)$ and $i\omega_n(u_n,\psi_n)$ have opposite real parts and cancel when the two identities are added. Using $\mathrm{Re}(u_{n,xxxx},u_n)=\|u_{n,xx}\|^2$ and $\mathrm{Re}\big(\int g\sigma_{n,xxxx},u_n\big)=\mathrm{Re}\int g(\sigma_{n,xx},u_{n,xx})$, this yields \eqref{eq:app_uxx}. Each of the five terms on the right-hand side of \eqref{eq:app_uxx} tends to zero as $n\to+\infty$, uniformly in $\eps$:
\begin{itemize}
\item[(i)] $|\mathrm{Re}(\psi_n,f_n^1)|\le\|\psi_n\|\|f_n^1\|\to 0$ by \eqref{eq:app_decay} and $\|F_n\|_{\Hcal_{\eps_n}}\to 0$;
\item[(ii)] $\|\psi_n\|^2\to 0$ by \eqref{eq:app_decay};
\item[(iii)] $\gamma_1|\mathrm{Re}(\psi_n,u_n)|\le\gamma_1\|\psi_n\|\|u_n\|\to 0$, since $u_n\to 0$ in $L^2$ and $\|\psi_n\|\to 0$;
\item[(iv)] $|\mathrm{Re}(f_n^2,u_n)|\le\|f_n^2\|\|u_n\|\to 0$, since $u_n\to 0$ in $L^2$ and $\|f_n^2\|\to 0$;
\item[(v)] the memory term is bounded by $\beta_1^{1/2}\|\sigma_n\|_{\Mg{1}}\|u_{n,xx}\|=o(1)\|u_{n,xx}\|$ via \eqref{eq:app_conv_bound_1} and \eqref{eq:app_decay}.
\end{itemize}
None of the constants involved in (i)--(v) depends on $\eps$. Consequently, \eqref{eq:app_uxx} gives $\bar\mu_1\|u_{n,xx}\|^2\le o(1)+o(1)\|u_{n,xx}\|$, where $o(1)$ denotes a quantity tending to zero as $n\to+\infty$, uniformly in $\eps$. By Young inequality, $\|u_{n,xx}\|^2\le o(1)$, so $u_{n,xx}\to 0$ in $L^2(0,\ell)$. The same argument applied to the fifth and fourth components, using $i\omega_n v_n=f_n^4-\phi_n$, gives $v_{n,x}\to 0$ in $L^2(0,\ell)$; the corresponding identity reads
\begin{equation}\label{eq:app_vxx}
\begin{aligned}
\bar\mu_2\|v_{n,x}\|^2
&=\|\phi_n\|^2-\mathrm{Re}(\phi_n,f_n^4)-\gamma_2\mathrm{Re}(\phi_n,v_n)\\
&\quad-\mathrm{Re}(f_n^5,v_n)
-\mathrm{Re}\!\int_0^{+\infty}\!\!g_{2,\eps_n}(s)(\eta_{n,x}(s),v_{n,x})\dd s,
\end{aligned}
\end{equation}
and each of its five terms tends to zero uniformly in $\eps$ by the same reasoning as in (i)--(v). Consequently, all six components of $\|U_n\|_{\Hcal_{\eps_n}}^2$ tend to zero, contradicting $\|U_n\|_{\Hcal_{\eps_n}}=1$.

Both cases being impossible, we conclude that there exist no sequences $\eps_n\to0^+$, $\omega_n\in\R$, and $U_n$ with the above properties. Hence $i\R\subset\rho(\Lcal_\eps)$ for every $\eps\in(0,1]$, uniformly in $\eps$.

\smallskip
\noindent\emph{Condition (C2).} Let $\Upsilon=(h^1,\dots,h^6)\in\Hcal_\eps$ and let $U=(u,\psi,\sigma,v,\phi,\eta)\in\Dcal(\Lcal_\eps)$ solve
\begin{equation}\label{eq:app_resolvent}
i\tau U-\Lcal_\eps U=\Upsilon,\qquad \tau\in\R.
\end{equation}
We prove that there exists $C>0$, independent of $\tau$ and $\eps$, such that
\begin{equation}\label{eq:app_C2}
\|U\|_{\Hcal_\eps}\le C\|\Upsilon\|_{\Hcal_\eps}.
\end{equation}

Taking the real part of the inner product of \eqref{eq:app_resolvent} with $U$ and using \eqref{eq:app_monotone},
\begin{equation}\label{eq:app_basic}
\gamma_1\|\psi\|^2+\frac{\lambda_1}{2}\|\sigma\|_{\Mg{1}}^2
+\gamma_2\|\phi\|^2+\frac{\lambda_2}{2}\|\eta\|_{\Mg{2}}^2
\le\|U\|_{\Hcal_\eps}\|\Upsilon\|_{\Hcal_\eps}.
\end{equation}
Hence
\begin{equation}\label{eq:app_basic'}
\|\psi\|,\ \|\sigma\|_{\Mg{1}},\ \|\phi\|,\ \|\eta\|_{\Mg{2}}
\le C\sqrt{\|U\|_{\Hcal_\eps}\|\Upsilon\|_{\Hcal_\eps}},
\end{equation}
where $C$ depends only on $\gamma_i,\lambda_i$; in particular, $C$ does not depend on $\eps$.

From the first and fourth components of \eqref{eq:app_resolvent}, $i\tau u=h^1-\psi$ and $i\tau v=h^4-\phi$. Taking the $L^2$-inner products of the second and fifth components with $u$ and $v$, respectively, and adding the real parts, we obtain, after substituting $i\tau u=h^1-\psi$ and $i\tau v=h^4-\phi$ and recalling that $(f,\tau g)=\bar\tau(f,g)$ with the above convention,
\begin{equation}\label{eq:app_uv}
\begin{aligned}
\bar\mu_1\|u_{xx}\|^2+\bar\mu_2\|v_x\|^2
&=-\mathrm{Re}(\psi,h^1)-\mathrm{Re}(\phi,h^4)
+\|\psi\|^2+\|\phi\|^2\\
&\quad-\gamma_1\mathrm{Re}(\psi,u)-\gamma_2\mathrm{Re}(\phi,v)
-\mathrm{Re}(h^2,u)-\mathrm{Re}(h^5,v)\\
&\quad-\mathrm{Re}\!\int_0^{+\infty}\!\!g_{1,\eps}(s)(\sigma_{xx}(s),u_{xx})\dd s
-\mathrm{Re}\!\int_0^{+\infty}\!\!g_{2,\eps}(s)(\eta_x(s),v_x)\dd s.
\end{aligned}
\end{equation}
We bound the right-hand side term by term. First, by Young's inequality, \eqref{eq:app_basic'} and $\|h^1\|_{L^2}\le\|h^1\|_{H_a^2}\le\|\Upsilon\|_{\Hcal_\eps}$,
\[
-\mathrm{Re}(\psi,h^1)+\|\psi\|^2
\le\tfrac12\|h^1\|_{L^2}^2+\tfrac32\|\psi\|^2
\le C\|U\|_{\Hcal_\eps}\|\Upsilon\|_{\Hcal_\eps}
+C\|\Upsilon\|_{\Hcal_\eps}^2,
\]
and similarly for $-\mathrm{Re}(\phi,h^4)+\|\phi\|^2$; both are of the form already present on the right-hand side of \eqref{eq:app_uv_bound} below. Second, by Poincar\'e's inequality and \eqref{eq:app_basic'},
\begin{align*}
\gamma_1|\mathrm{Re}(\psi,u)|+\gamma_2|\mathrm{Re}(\phi,v)|
&\le\epsilon\big(\|u_{xx}\|^2+\|v_x\|^2\big)
+C_\epsilon\big(\|\psi\|^2+\|\phi\|^2\big)\\
&\le\epsilon\big(\|u_{xx}\|^2+\|v_x\|^2\big)
+C_\epsilon\|U\|_{\Hcal_\eps}\|\Upsilon\|_{\Hcal_\eps}.
\end{align*}
Third, again by Poincar\'e's inequality and $\|h^2\|_{L^2}\le\|\Upsilon\|_{\Hcal_\eps}$,
\[
|\mathrm{Re}(h^2,u)|+|\mathrm{Re}(h^5,v)|
\le\epsilon\big(\|u_{xx}\|^2+\|v_x\|^2\big)
+C_\epsilon\|\Upsilon\|_{\Hcal_\eps}^2.
\]
Fourth, by \eqref{eq:app_conv_bound_1}--\eqref{eq:app_conv_bound_2} and \eqref{eq:app_basic'},
\begin{align*}
&\Big|\mathrm{Re}\!\int_0^{+\infty}\!\!g_{1,\eps}(s)(\sigma_{xx}(s),u_{xx})\dd s\Big|
+\Big|\mathrm{Re}\!\int_0^{+\infty}\!\!g_{2,\eps}(s)(\eta_x(s),v_x)\dd s\Big|\\
&\qquad\le\epsilon\big(\|u_{xx}\|^2+\|v_x\|^2\big)
+C_\epsilon\|U\|_{\Hcal_\eps}\|\Upsilon\|_{\Hcal_\eps}.
\end{align*}
Summing all these and choosing $\epsilon$ so small that $4\epsilon<\min\{\bar\mu_1,\bar\mu_2\}/2$, we absorb the $\epsilon(\|u_{xx}\|^2+\|v_x\|^2)$ terms into the left-hand side and obtain
\begin{equation}\label{eq:app_uv_bound}
\bar\mu_1\|u_{xx}\|^2+\bar\mu_2\|v_x\|^2
\le C\big(\|U\|_{\Hcal_\eps}\|\Upsilon\|_{\Hcal_\eps}+\|\Upsilon\|_{\Hcal_\eps}^2\big),
\end{equation}
where $C$ depends only on $\gamma_i,\mu_i,\beta_i,c_e,c_p$, but not on $\tau$ or $\eps$. The only place where the kernel enters is the pair of memory convolutions, controlled by the mass $\beta_i$ through \eqref{eq:app_conv_bound_1}--\eqref{eq:app_conv_bound_2}; by \eqref{eq:F1}, $\beta_i$ is independent of $\eps$.

Combining \eqref{eq:app_basic}, \eqref{eq:app_uv_bound} and the identity
\[
\|U\|_{\Hcal_\eps}^2
=\bar\mu_1\|u_{xx}\|^2+\bar\mu_2\|v_x\|^2
+\|\psi\|^2+\|\sigma\|_{\Mg{1}}^2+\|\phi\|^2+\|\eta\|_{\Mg{2}}^2,
\]
we obtain
\begin{equation}\label{eq:app_X2}
\|U\|_{\Hcal_\eps}^2
\le C\big(\|U\|_{\Hcal_\eps}\|\Upsilon\|_{\Hcal_\eps}+\|\Upsilon\|_{\Hcal_\eps}^2\big).
\end{equation}
If $\Upsilon=0$, then \eqref{eq:app_X2} gives $\|U\|_{\Hcal_\eps}=0$, and \eqref{eq:app_C2} is trivial. Otherwise, set $t=\|U\|_{\Hcal_\eps}/\|\Upsilon\|_{\Hcal_\eps}$. Dividing \eqref{eq:app_X2} by $\|\Upsilon\|_{\Hcal_\eps}^2$ gives $t^2\le C(t+1)$, hence
\[
t\le C'\coloneqq \frac{C+\sqrt{C^2+4C}}{2},
\]
a constant depending only on $C$, hence only on $\gamma_i,\mu_i,\beta_i,c_e,c_p$, but not on $\tau$ or $\eps$. Therefore $\|U\|_{\Hcal_\eps}\le C'\|\Upsilon\|_{\Hcal_\eps}$, which is \eqref{eq:app_C2}.

We have verified conditions (C1) and (C2) of the Gearhart--Pr\"uss theorem with constants independent of $\eps$. Therefore the linear semigroup $e^{-\Lcal_\eps t}$ is exponentially stable, uniformly in $\eps$: there exist $\alpha\ge 1$ and $\nu>0$, independent of $\eps$, such that
\[
\|e^{-\Lcal_\eps t}\|_{\mathcal{L}(\Hcal_\eps)}\le \alpha e^{-\nu t},
\qquad \forall t\ge 0,\ \forall \eps\in(0,1].
\]
This completes the proof of Theorem~\ref{thm:unif_exp}.
\end{proof}

Formally passing to the limit $\eps\to 0^+$ in \eqref{eq:Mukiawa} with $g_i$ replaced by $g_{i,\eps}$, and using $g_{i,\eps}\rightharpoonup\beta_i\delta_0$, the memory convolution satisfies, in the distributional sense,
\[
-\int_0^{+\infty}g_{1,\eps}(s)u_{xxxx}(x,t-s)\dd s
\rightharpoonup -\beta_1u_{xxxx}(x,t),
\]
and similarly
\[
+\int_0^{+\infty}g_{2,\eps}(s)v_{xx}(x,t-s)\dd s
\rightharpoonup +\beta_2v_{xx}(x,t).
\]
Combining these limits with the $\mu_1u_{xxxx}$ and $-\mu_2v_{xx}$ terms in \eqref{eq:Mukiawa} gives $(\mu_1-\beta_1)u_{xxxx}=\bar\mu_1u_{xxxx}$ and $-(\mu_2-\beta_2)v_{xx}=-\bar\mu_2v_{xx}$. Equivalently, in the rewritten system \eqref{eq:Aeps} the convolution $\int_0^{+\infty}g_{1,\eps}(s)\sigma_{xxxx}(s)\dd s$ vanishes in the limit, because $\sigma(x,t,s)=u(x,t)-u(x,t-s)=O(s)$ while $g_{1,\eps}$ has first moment $O(\eps)$ (see also \eqref{eq:memory_vanishing_compact} below). Since $\bar\mu_i=\mu_i-\beta_i$, the limiting system becomes
\begin{equation}\label{eq:limit}
\begin{cases}
u_{tt}+\gamma_1u_t+\bar\mu_1u_{xxxx}+k[u-v]^++h_1(u)=f_1(x),\\[1mm]
v_{tt}+\gamma_2v_t-\bar\mu_2v_{xx}-k[u-v]^++h_2(v)=f_2(x),
\end{cases}
\end{equation}
where the memory convolutions have disappeared and the elastic coefficients are precisely the effective stiffnesses $\bar\mu_i=\mu_i-\beta_i$, which are independent of $\eps$. Thus the limit system has the same form as the classical memoryless coupled suspension bridge system \cite{MaZhong2009}; relabelling $\bar\mu_i$ as $\mu_i$ recovers the notation used there, but the original coefficients $\mu_i$ of \eqref{eq:Mukiawa} are not retained in the limit.

The limit semigroup $S_0(t)$ is generated by $\Lcal_0=\Lcal_{e,0}$ on $\Hcal_e$, with coefficients $\bar\mu_i=\mu_i-\beta_i$. Its attractor is constructed inside the present framework.

\begin{proposition}[Existence of the limit attractor]\label{prop:attractor_0}
Under \eqref{eq:kernel} and \eqref{eq:h_assumptions}, the semigroup $S_0(t)$ generated by \eqref{eq:limit} possesses a compact global attractor $\Acal_0\subset\Hcal_e$.
\end{proposition}
\begin{proof}
The energy $E_0(t)$, obtained from \eqref{eq:proof_energy} by deleting the memory norms, satisfies the differential inequality \eqref{eq:proof_differential} with the memory terms omitted, because $\bar\mu_i>0$ and \eqref{eq:h_consequences} is unchanged. Hence $S_0$ has an absorbing ball $\Bcal_{0,0}\subset\Hcal_e$ of radius $R_0$. We cannot simply put $\eps=0$ in Theorem~\ref{thm:unif_exp}, whose proof uses the rescaled kernels. For the memoryless linear system, adding the multiplier $\omega[(u,u_t)+(v,v_t)]$ to $E_0$ gives
\[
\frac{\dd}{\dd t}E_{0,\omega}(t)+c\,E_{0,\omega}(t)\le0
\]
for some $c>0$ depending only on $\gamma_i,\bar\mu_i$ and the domain, so that $\|e^{-\Lcal_{e,0}t}\|_{\mathcal{L}(\Hcal_e)}\le\alpha e^{-\nu t}$ with constants no worse than in \eqref{eq:exp_stab}.

It remains to justify compactness of the integral remainder. Let $\Pcal_N$ be the Fourier projection onto the first $N$ modes and define
\[
K_0\Phi:=\int_0^T e^{-\Lcal_{e,0}(T-s)}\Phi(s)\,\dd s,\qquad \Phi\in L^2(0,T;\Hcal_e).
\]
Since $\|e^{-\Lcal_{e,0}\tau}\|\le\alpha$ is locally integrable at $\tau=0$, the integral is bounded from $L^2(0,T;\Hcal_e)$ into $\Hcal_e$. For fixed $N$, $\Pcal_NK_0$ takes values in a finite-dimensional subspace. For the high modes, the exponential estimate gives
\[
\|(I-\Pcal_N)K_0\Phi\|_{\Hcal_e}
\le\alpha\int_0^T\|(I-\Pcal_N)\Phi(s)\|_{\Hcal_e}\,\dd s
\le\alpha\sqrt T\,\|(I-\Pcal_N)\Phi\|_{L^2(0,T;\Hcal_e)} .
\]
Thus $K_0$ maps any relatively compact subset of $L^2(0,T;\Hcal_e)$ into a relatively compact subset of $\Hcal_e$; the singularity at $\tau=0$ causes no difficulty because it is only locally integrable.

Finally, the forcing set $\{F(S_0(\cdot)U):U\in\Bcal_{0,0}\}$ is relatively compact in $L^2(0,T;\Hcal_e)$. Indeed, from the energy bounds, $u$ is bounded in $L^\infty(0,T;H_a^2)$ and $u_t$ in $L^2(0,T;L^2)$; interpolation and the fundamental theorem of calculus give
\[
\|u(t)-u(t-h)\|_{H^1}\le C h^{1/4},
\]
uniformly on $\Bcal_{0,0}$, and similarly for $v$. By Arzel\`a--Ascoli, the elastic displacements form a relatively compact set in $C([0,T];L^2)$. A finite-mode net then controls the projected nonlinearities, while the high-mode tail is small by the same Fourier projection argument used in Proposition~\ref{prop:quasi}. Hence $K_0(\{F(S_0(\cdot)U):U\in\Bcal_{0,0}\})$ is relatively compact in $\Hcal_e$. Consequently $S_0$ is asymptotically smooth; asymptotic smoothness and dissipativity yield the compact global attractor $\Acal_0$. We do not invoke the stronger $C^3/C^2$ hypotheses of \cite{MaZhong2009}; that paper is cited only for the classical memoryless model.
\end{proof}

\section{Uniform Dissipativity}\label{sec:dissipativity}

Throughout this section we use the abbreviation $\bar\mu_i=\bar\mu_{i,\eps}=\mu_i-\beta_i$, which is independent of $\eps$ by \eqref{eq:bar_mu} and \eqref{eq:F1}; we also write $\Mg{i}=\Mcal_{g_{i,\eps}}$.

The first step toward any continuity result is a \emph{uniform dissipativity estimate}: there is a radius $R_0>0$, independent of $\eps$, such that the balls $\Bcal_{0,\eps}=\{U\in\Hcal_\eps:\|U\|_{\Hcal_\eps}\le R_0\}$ absorb all bounded sets uniformly in $\eps$. We write $\Bcal_0$ for this family and suppress the subscript $\eps$ when no confusion can arise.

\begin{proposition}[Uniform absorbing set]\label{prop:absorbing}
Let \eqref{eq:kernel} and \eqref{eq:h_assumptions} hold. Then there exist $R_0>0$ and, for every bounded set $B\subset\Hcal_\eps$, a time $t_0=t_0(B)\ge 0$ such that, uniformly in $\eps\in(0,1]$,
\[
S_\eps(t)B\subset \Bcal_{0,\eps}\coloneqq\{U\in\Hcal_\eps:\|U\|_{\Hcal_\eps}\le R_0\},
\qquad \forall t\ge t_0.
\]
Moreover, $R_0$ depends only on $\gamma_i,\mu_i,\lambda_i,\beta_i$, $c_e$, $c_p$, $\|f_1\|,\|f_2\|$, and the constants in \eqref{eq:h_assumptions}, but not on $B$ or on $\eps$.
\end{proposition}

\begin{proof}
Let $B\subset\Hcal_\eps$ be a bounded set, and let $U_{0,\eps}\in B$ be arbitrary. We write $U_\eps(t)=S_\eps(t)U_{0,\eps}=(u_\eps,u_{\eps,t},\sigma_\eps,v_\eps,v_{\eps,t},\eta_\eps)$ for the corresponding solution of the perturbed system. Define the energy functional
\begin{equation}\label{eq:proof_energy}
\begin{aligned}
E_\eps(t)
&=\frac12\Big[\|u_{\eps,t}\|^2+\bar\mu_1\|u_{\eps,xx}\|^2
+\|\sigma_\eps\|_{\Mg{1}}^2+\|v_{\eps,t}\|^2+\bar\mu_2\|v_{\eps,x}\|^2
+\|\eta_\eps\|_{\Mg{2}}^2\Big]\\
&\quad+\frac{k}{2}\|[u_\eps-v_\eps]^+\|^2
+(H_1(u_\eps),1)+(H_2(v_\eps),1)
-(f_1,u_\eps)-(f_2,v_\eps).
\end{aligned}
\end{equation}
Here $\bar\mu_i=\mu_i-\beta_i$ is independent of $\eps$ by the invariance of the mass \eqref{eq:F1}. Throughout, $c_e,c_p>0$ denote the Poincar\'e constants satisfying
\begin{equation}\label{eq:poincare_app}
\|w\|\le\sqrt{c_e}\|w_{xx}\|,\quad w\in H_a^2(0,\ell),
\qquad
\|w\|\le\sqrt{c_p}\|w_x\|,\quad w\in H_0^1(0,\ell),
\end{equation}
and $C_e>0$ denotes the embedding constant of $H_a^2(0,\ell)\hookrightarrow L^\infty(0,\ell)$ and $H_0^1(0,\ell)\hookrightarrow L^\infty(0,\ell)$, so that
\begin{equation}\label{eq:proof_embedding}
\|w\|_{L^\infty}\le C_e\|w_{xx}\|,\quad w\in H_a^2(0,\ell),
\qquad
\|w\|_{L^\infty}\le C_e\|w_x\|,\quad w\in H_0^1(0,\ell).
\end{equation}

\smallskip
\noindent\emph{Step 1: Dissipativity of $E_\eps$.} Multiplying the first equation of the perturbed system by $u_{\eps,t}$ in $L^2(0,\ell)$ and the second by $v_{\eps,t}$, and using the transport equations $\partial_t\sigma_\eps+\partial_s\sigma_\eps-u_{\eps,t}=0$ and $\partial_t\eta_\eps+\partial_s\eta_\eps-v_{\eps,t}=0$, we obtain
\begin{equation}\label{eq:proof_dissipation}
\frac{\dd}{\dd t}E_\eps(t)
=-\gamma_1\|u_{\eps,t}\|^2-\gamma_2\|v_{\eps,t}\|^2
+\frac12\int_0^{+\infty}g_{1,\eps}'(s)\|\sigma_{\eps,xx}(s)\|^2\dd s
+\frac12\int_0^{+\infty}g_{2,\eps}'(s)\|\eta_{\eps,x}(s)\|^2\dd s.
\end{equation}
Indeed, the contribution of the memory terms to $\frac{\dd}{\dd t}E_\eps$ is obtained by differentiating $\frac12\|\sigma_\eps\|_{\Mg{1}}^2=\frac12\int_0^{+\infty}g_{1,\eps}(s)\|\sigma_{\eps,xx}(s)\|^2\dd s$ with respect to $t$, using the transport equation to replace $\partial_t\sigma_\eps$ by $u_{\eps,t}-\partial_s\sigma_\eps$, integrating by parts in $s$, and using $\sigma_\eps(x,t,0)=0$ and $g_{1,\eps}(+\infty)=0$; the resulting term $-\int_0^{+\infty}g_{1,\eps}(s)(\sigma_{\eps,xx}(s),u_{\eps,xx,t})\dd s$ then cancels against the corresponding term $+\int_0^{+\infty}g_{1,\eps}(s)(\sigma_{\eps,xx}(s),u_{\eps,xx,t})\dd s$ produced by differentiating the memory energy. The same reasoning applies to $\eta_\eps$. By the differential inequality $g_{i,\eps}'(s)+\lambda_i g_{i,\eps}(s)\le 0$ for a.e.\ $s\ge 0$ (which follows from $g_{i,\eps}'+(\lambda_i/\eps)g_{i,\eps}\le0$ since $\eps\le1$), we deduce
\begin{equation}\label{eq:proof_dissipation2}
\frac{\dd}{\dd t}E_\eps(t)
\le-\gamma_1\|u_{\eps,t}\|^2-\gamma_2\|v_{\eps,t}\|^2
-\frac{\lambda_1}{2}\|\sigma_\eps\|_{\Mg{1}}^2
-\frac{\lambda_2}{2}\|\eta_\eps\|_{\Mg{2}}^2\le 0.
\end{equation}
In particular, $E_\eps(t)\le E_\eps(0)$ for all $t\ge 0$.

\smallskip
\noindent\emph{Step 2: Modified energy.} Define
\begin{equation}\label{eq:proof_modified}
E_{\eps,\omega}(t)\coloneqq E_\eps(t)
+\omega(u_\eps,u_{\eps,t})+\omega(v_\eps,v_{\eps,t}),
\end{equation}
where $\omega>0$ is a parameter to be chosen below, independent of $\eps$. Multiplying the first equation of the perturbed system by $u_\eps$ in $L^2(0,\ell)$ and the second by $v_\eps$, and adding the results to \eqref{eq:proof_dissipation}, we obtain
\begin{equation}\label{eq:proof_derivative}
\begin{aligned}
\frac{\dd}{\dd t}E_{\eps,\omega}(t)
&+(\gamma_1-\omega)\|u_{\eps,t}\|^2+(\gamma_2-\omega)\|v_{\eps,t}\|^2
+\omega\bar\mu_1\|u_{\eps,xx}\|^2+\omega\bar\mu_2\|v_{\eps,x}\|^2\\
&+\omega\gamma_1(u_\eps,u_{\eps,t})+\omega\gamma_2(v_\eps,v_{\eps,t})
+\omega k\|[u_\eps-v_\eps]^+\|^2\\
&+\frac{\lambda_1}{2}\|\sigma_\eps\|_{\Mg{1}}^2+\frac{\lambda_2}{2}\|\eta_\eps\|_{\Mg{2}}^2\\
&+\omega\int_0^{+\infty}g_{1,\eps}(s)(\sigma_{\eps,xx}(s),u_{\eps,xx})\dd s
+\omega\int_0^{+\infty}g_{2,\eps}(s)(\eta_{\eps,x}(s),v_{\eps,x})\dd s\\
&-\omega(h_1(u_\eps),u_\eps)-\omega(h_2(v_\eps),v_\eps)
-\omega(f_1,u_\eps)-\omega(f_2,v_\eps)\le 0.
\end{aligned}
\end{equation}
Note the minus sign in front of $\omega(h_i(\cdot),\cdot)$: it comes from substituting $u_{tt}$ and $v_{tt}$ from the equations into $\omega(u,u_{tt})$ and $\omega(v,v_{tt})$.

Using the growth condition \eqref{eq:h_consequences} on $H_i$ and $h_i$, there exists $C_1>0$ (the same constant as in \eqref{eq:h_assumptions}), independent of $\eps$, such that
\begin{equation}\label{eq:proof_nonlinear}
-\omega(h_1(u_\eps),u_\eps)-\omega(h_2(v_\eps),v_\eps)
\le-\omega C_1(H_1(u_\eps),1)-\omega C_1(H_2(v_\eps),1)
+\omega\epsilon\|u_\eps\|^2+\omega\epsilon\|v_\eps\|^2+\omega C_\epsilon',
\end{equation}
where $C_\epsilon'$ is the constant appearing in the second inequality of \eqref{eq:h_consequences}. By Poincar\'e's inequality \eqref{eq:poincare_app},
\[
\omega\epsilon\|u_\eps\|^2\le\omega\epsilon c_e\|u_{\eps,xx}\|^2,
\qquad
\omega\epsilon\|v_\eps\|^2\le\omega\epsilon c_p\|v_{\eps,x}\|^2.
\]
Similarly, by the Cauchy--Schwarz and Young inequalities,
\begin{equation}\label{eq:proof_memory}
\Big|\int_0^{+\infty}g_{1,\eps}(s)(\sigma_{\eps,xx}(s),u_{\eps,xx})\dd s\Big|
\le c_{\delta_1}\|\sigma_\eps\|_{\Mg{1}}^2+\delta_1\bar\mu_1\|u_{\eps,xx}\|^2,
\end{equation}
and the analogous estimate for $\eta_\eps$ and $v_\eps$ reads
\begin{equation}\label{eq:proof_memory_eta}
\Big|\int_0^{+\infty}g_{2,\eps}(s)(\eta_{\eps,x}(s),v_{\eps,x})\dd s\Big|
\le c_{\delta_1}'\|\eta_\eps\|_{\Mg{2}}^2+\delta_1\bar\mu_2\|v_{\eps,x}\|^2,
\end{equation}
where
\begin{equation}\label{eq:c_delta_def}
c_{\delta_1}=\frac{\beta_1}{4\delta_1\bar\mu_1},\qquad
c_{\delta_1}'=\frac{\beta_2}{4\delta_1\bar\mu_2},
\end{equation}
and $\delta_1>0$ is arbitrary. Note that $c_{\delta_1}$ and $c_{\delta_1}'$ depend only on $\beta_i$, $\delta_1$, and $\bar\mu_i$, and are independent of $\eps$.

\smallskip
\noindent\emph{Step 3: Choice of parameters.} We choose the parameters in the following order, all independently of $\eps$. First, set $\delta_1=1/4$. Second, choose $\epsilon>0$ so small that
\begin{equation}\label{eq:proof_eps_choice}
\epsilon\le\min\Big\{\frac{\bar\mu_1}{8c_e},\frac{\bar\mu_2}{8c_p}\Big\}.
\end{equation}
Third, choose $\omega>0$ so small that
\begin{equation}\label{eq:proof_omega_choice}
\omega\le\min\Big\{\frac{\gamma_1}{2},\frac{\gamma_2}{2},
\frac{\lambda_1}{4c_{\delta_1}},\frac{\lambda_2}{4c_{\delta_1}'},
\frac{\bar\mu_1}{2},\frac{\bar\mu_2}{2},
\frac{1}{4c_e},\frac{1}{4c_p},\frac12,
\frac{\bar\mu_1}{4c_e},\frac{\bar\mu_2}{4c_p},
\frac{1}{2C_1}\Big\}.
\end{equation}
All the quantities appearing in \eqref{eq:proof_eps_choice}--\eqref{eq:proof_omega_choice}---$\gamma_i,\lambda_i,\bar\mu_i,c_e,c_p,c_{\delta_1},c_{\delta_1}',C_1$---are independent of $\eps$. In particular, the last bound $\omega\le 1/(2C_1)$ guarantees $1-\omega C_1\ge 1/2>0$, which is needed when the terms $-\omega C_1(H_i(u_i),1)$ are combined with the terms $(H_i(u_i),1)$ present in $E_{\eps,\omega}$. Substituting the above estimates into \eqref{eq:proof_derivative} and absorbing the $\epsilon$-terms and the memory terms with the help of \eqref{eq:proof_eps_choice} and \eqref{eq:proof_omega_choice}, we arrive at
\begin{equation}\label{eq:proof_differential}
\frac{\dd}{\dd t}E_{\eps,\omega}(t)+c_0 E_{\eps,\omega}(t)\le c_1,
\end{equation}
where $c_0>0$ and $c_1\ge 0$ are constants depending only on $\gamma_i,\lambda_i,\mu_i,\beta_i,c_e,c_p,\|f_1\|,\|f_2\|$, and the constants in \eqref{eq:h_assumptions}, but not on $\eps$ or on the initial data.

\smallskip
\noindent\emph{Step 4: Equivalence of $E_{\eps,\omega}$ and $\|\cdot\|_{\Hcal_\eps}^2$.} Let $R_B>0$ satisfy $\|U_{0,\eps}\|_{\Hcal_\eps}\le R_B$ for all $U_{0,\eps}\in B$. We claim that there exist $c_2>0$ and $c_3=c_3(R_B)>0$, both independent of $\eps$, such that for every $t\ge 0$ with $\|U_\eps(t)\|_{\Hcal_\eps}\le R_B+1$,
\begin{equation}\label{eq:proof_equivalence}
c_2\big(\|U_\eps(t)\|_{\Hcal_\eps}^2-1\big)\le E_{\eps,\omega}(t)
\le c_3\big(\|U_\eps(t)\|_{\Hcal_\eps}^2+1\big).
\end{equation}

For the lower bound, we estimate each term in $E_{\eps,\omega}$ as follows:
\[
\omega(u_\eps,u_{\eps,t})+\omega(v_\eps,v_{\eps,t})
\ge-\frac{\omega}{2}\|u_{\eps,t}\|^2-\frac{\omega c_e}{2}\|u_{\eps,xx}\|^2
-\frac{\omega}{2}\|v_{\eps,t}\|^2-\frac{\omega c_p}{2}\|v_{\eps,x}\|^2,
\]
and, by \eqref{eq:h_consequences} and Young's inequality,
\begin{equation}\label{eq:proof_H_f}
(H_1(u_\eps),1)+(H_2(v_\eps),1)-(f_1,u_\eps)-(f_2,v_\eps)
\ge-\epsilon c_e\|u_{\eps,xx}\|^2-\epsilon c_p\|v_{\eps,x}\|^2-c_\epsilon,
\end{equation}
where $c_\epsilon$ depends only on $\epsilon$, $\|f_1\|$, $\|f_2\|$, and the constants in \eqref{eq:h_consequences}. The terms produced here are the \emph{combined} contribution of the $H_i$-terms (via \eqref{eq:h_consequences}) and of Young's inequality applied to the $f_i$-terms: the former yields $-\epsilon\|u_\eps\|^2$, the latter $-\tfrac{\epsilon}{2}\|u_\eps\|^2$, so that together they give $-\tfrac{3\epsilon c_e}{2}\|u_{\eps,xx}\|^2-C_\epsilon$ (and analogously for $v_\eps$). These are \emph{distinct} from the analogous terms $-\omega\epsilon c_e\|u_{\eps,xx}\|^2$ and $-\omega\epsilon c_p\|v_{\eps,x}\|^2$ produced in Step 2 by the $h_i$-terms.

Combining the above bounds with \eqref{eq:proof_energy} and \eqref{eq:proof_modified}, we see that the coefficients of $\|u_{\eps,t}\|^2$, $\|u_{\eps,xx}\|^2$, $\|v_{\eps,t}\|^2$, $\|v_{\eps,x}\|^2$ in the resulting lower bound for $E_{\eps,\omega}(t)$ are bounded below by
\begin{equation}\label{eq:proof_coeffs}
\frac12-\frac{\omega}{2},\qquad
\frac{\bar\mu_1}{2}-\frac{\omega c_e}{2}-\frac{3\epsilon c_e}{2}
-\omega\epsilon c_e,
\qquad
\frac{\bar\mu_2}{2}-\frac{\omega c_p}{2}-\frac{3\epsilon c_p}{2}
-\omega\epsilon c_p.
\end{equation}
The choices \eqref{eq:proof_eps_choice}--\eqref{eq:proof_omega_choice} already ensure that all three quantities are strictly positive. Indeed, since $\omega\le 1/2$, the first is at least $1/4$. For the second, using $\omega c_e\le\bar\mu_1/4$, $\epsilon c_e\le\bar\mu_1/8$, and $\omega\epsilon c_e\le\bar\mu_1/16$ (which follows from $\omega\le 1/2$ and $\epsilon c_e\le\bar\mu_1/8$), we get
\[
\frac{\bar\mu_1}{2}-\frac{\omega c_e}{2}-\frac{3\epsilon c_e}{2}
-\omega\epsilon c_e
\ge\frac{\bar\mu_1}{2}-\frac{\bar\mu_1}{8}-\frac{3\bar\mu_1}{16}
-\frac{\bar\mu_1}{16}=\frac{\bar\mu_1}{8}>0,
\]
and similarly for the $v$-coefficient. This yields the lower bound in \eqref{eq:proof_equivalence}, the constant term $-1$ on the left-hand side absorbing $c_\epsilon$, $\omega C_\epsilon'$, and the $\|f_i\|^2$-terms produced by Young's inequality. All these constants depend only on $\epsilon$, $\omega$, $\|f_1\|$, $\|f_2\|$, and the constants in \eqref{eq:h_consequences}, but not on $\eps$ or $R_B$; hence $c_2$ is independent of $R_B$, while $c_3$ below is not.

For the upper bound, we argue as follows. Fix any $t\ge 0$ with $\|U_\eps(t)\|_{\Hcal_\eps}\le R_B+1$. By the embedding \eqref{eq:proof_embedding} and the definition \eqref{eq:Hnorm} of $\|\cdot\|_{\Hcal_\eps}$, there exists a constant $C_e'>0$, depending only on $C_e$ and $\bar\mu_i$, such that
\[
\|u_\eps(t)\|_{L^\infty}\le C_e\|u_{\eps,xx}(t)\|\le C_e'\|U_\eps(t)\|_{\Hcal_\eps}
\le C_e'(R_B+1),
\]
and similarly for $v_\eps(t)$. Since $h_i$ is continuous and $H_i$ is its primitive, both functions are bounded on the compact interval $[-C_e'(R_B+1),C_e'(R_B+1)]$ by a constant $C_3=C_3(R_B)$ depending only on $R_B$, $C_e'$, and the constants in \eqref{eq:h_assumptions}, but not on $\eps$ or $t$. Hence
\[
|(H_i(u_\eps(t)),1)|\le C_3(R_B),\qquad
|(h_i(u_\eps(t)),u_\eps(t))|\le C_3(R_B),\qquad i=1,2.
\]
Combining this with \eqref{eq:proof_energy} and \eqref{eq:proof_modified}, and using the boundedness of $u_\eps,u_{\eps,t},v_\eps,v_{\eps,t},\sigma_\eps,\eta_\eps$ in their respective norms (all controlled by $\|U_\eps(t)\|_{\Hcal_\eps}\le R_B+1$), we obtain the upper bound in \eqref{eq:proof_equivalence} with a constant $c_3$ that depends on $R_B$, $\gamma_i,\mu_i,\beta_i,c_e,c_p$, and the constants in \eqref{eq:h_assumptions}, but not on $\eps$.

\smallskip
\noindent\emph{Step 5: Absorbing ball.} The lower bound in \eqref{eq:proof_equivalence} holds for every $t\ge0$ without any upper bound on $\|U_\eps(t)\|_{\Hcal_\eps}$: its derivation uses only \eqref{eq:h_consequences}, Young's inequality and the Poincar\'e inequalities. The differential inequality \eqref{eq:proof_differential} is likewise global. Hence for every $t\ge0$,
\[
c_2\big(\|U_\eps(t)\|_{\Hcal_\eps}^2-1\big)\le E_{\eps,\omega}(t)
\le c_3(R_B)(R_B^2+1)e^{-c_0t}+\frac{c_1}{c_0},
\]
where the upper estimate in \eqref{eq:proof_equivalence} is used only at $t=0$, with $\|U_{0,\eps}\|_{\Hcal_\eps}\le R_B$. Therefore
\[
\|U_\eps(t)\|_{\Hcal_\eps}^2\le
1+\frac{c_3(R_B)(R_B^2+1)e^{-c_0t}+c_1/c_0}{c_2},\qquad \forall t\ge0.
\]
Choose $t_0=t_0(B)\ge0$, independent of $\eps$, such that
$c_3(R_B)(R_B^2+1)e^{-c_0t_0}\le c_2$, for instance
$t_0:=\max\{0,c_0^{-1}\log(c_3(R_B)(R_B^2+1)/c_2)\}$, and put
\[
R_0:=\Big(2+\frac{c_1}{c_0c_2}\Big)^{1/2}.
\]
Then $\|U_\eps(t)\|_{\Hcal_\eps}\le R_0$ for every $t\ge t_0$, uniformly in $\eps$, which proves the proposition. Since the a priori bound holds on every finite interval, global existence follows by restarting the local solution of Proposition~\ref{prop:wellposed} at successive times.
\end{proof}

\begin{lemma}[Finite-time uniform boundedness]\label{lem:finite_time}
For every $T>0$ and every bounded set $B\subset\Hcal_\eps$, there exists $R_1=R_1(T,B)>0$, independent of $\eps\in(0,1]$, such that
\[
\sup_{\eps\in(0,1]}\sup_{U_0\in B}\sup_{0\le t\le T}\|S_\eps(t)U_0\|_{\Hcal_\eps}\le R_1.
\]
\end{lemma}

\begin{proof}
Let $R_B$ be a bound for $B$ in $\Hcal_\eps$. By Step 5 of Proposition~\ref{prop:absorbing}, for every $U_0\in B$, every $t\ge0$ and every $\eps\in(0,1]$,
\[
\|S_\eps(t)U_0\|_{\Hcal_\eps}^2
\le 1+\frac{c_3(R_B)(R_B^2+1)+c_1/c_0}{c_2}=:R_1^2 .
\]
Here the radius-dependent constant $c_3(R_B)$ is evaluated at the initial radius only; all constants are independent of $\eps$ and of the particular $U_0\in B$, and the estimate is valid on $[0,T]$. Thus $R_1$ is a uniform bound on $[0,T]$, and the lemma follows.
\end{proof}

\section{Upper Semicontinuity}\label{sec:upper}

We now give the complete proof of the upper semicontinuity \eqref{eq:upper_intro}. The general strategy is standard (see, e.g., \cite[Theorem 7.9.6]{ChueshovLasiecka2008} and \cite{Robinson2001,Hale1988}):
\begin{enumerate}
\item[(A1)] Uniform dissipativity: the semigroups $S_\eps(t)$ possess the uniform family of absorbing balls $\Bcal_{0,\eps}$ of radius $R_0$ independent of $\eps\in(0,1]$ (Proposition~\ref{prop:absorbing}).
\item[(A2)] Uniform asymptotic smoothness of the semigroups $S_\eps(t)$, uniformly in $\eps\in(0,1]$. In the framework of Chueshov--Lasiecka, this is implied by the uniform quasi-stability inequality on $\Bcal_0$ (Proposition~\ref{prop:quasi}).
\item[(A3)] Convergence of the perturbed semigroups to the limit semigroup on finite time intervals: for every $T>0$ and every sequence $\eps_n\to 0^+$ with $J_{\eps_n}U_{0,\eps_n}\to J_0U_0$ for some $U_0\in\Hcal_e$,
\begin{equation}\label{eq:semigroup_conv}
\sup_{t\in[0,T]}\|S_{\eps_n}(t)U_{0,\eps_n}-S_0(t)U_0\|_{\Hcal_{\eps_n}}\to 0.
\end{equation}
\end{enumerate}
Item (A3) is Corollary~\ref{cor:conv_general}. Relative compactness of sequences in $\Acal_{\eps_n}$ and the fact that their accumulation points have zero memory and belong to $\Acal_0$ are proved in Lemma~\ref{lem:attractor_memory} and in the proof of Theorem~\ref{thm:upper}; see also the paragraph containing (H3).

\subsection{Uniform asymptotic smoothness}

The uniform exponential stability \eqref{eq:exp_stab} of the linear semigroup, combined with the quasi-stability framework of Chueshov--Lasiecka \cite{ChueshovLasiecka2008} as developed in \cite{Zhou2026}, yields the following uniform version.

\begin{proposition}[Uniform quasi-stability]\label{prop:quasi}
There exist $T>0$ and $q\in[0,1)$, both independent of $\eps$, such that for every $\eps\in(0,1]$ the following holds: there exists a map $K_{\eps,T}:\Bcal_0\to\Hcal_\eps$ satisfying
\begin{equation}\label{eq:quasi_stab}
\|S_\eps(T)U_0-S_\eps(T)V_0\|_{\Hcal_\eps}^2
\le q^2\|U_0-V_0\|_{\Hcal_\eps}^2
+\|K_{\eps,T}(U_0)-K_{\eps,T}(V_0)\|_{\Hcal_\eps}^2,
\end{equation}
for all $U_0,V_0\in\Bcal_0$, and such that the following sequential compactness property holds: for every sequence $\eps_n\to0^+$ and every bounded sequence $U_n\in\Bcal_0$, there exist a subsequence $\{\eps_{n_k},U_{n_k}\}$ and an element $K_*$ in the sense of Remark~\ref{rem:uniform_topology} such that
\begin{equation}\label{eq:K_seq_compact}
K_{\eps_{n_k},T}(U_{n_k})\to K_*
\qquad(k\to+\infty)
\end{equation}
in $\Hcal_*$.
\end{proposition}

\begin{proof}
Fix $\eps\in(0,1]$. Let $U_0,V_0\in\Bcal_0$ and set $U_\eps(t)=S_\eps(t)U_0$, $V_\eps(t)=S_\eps(t)V_0$. Define $W_\eps(t)\coloneqq U_\eps(t)-V_\eps(t)$ and $\Phi_\eps(t)\coloneqq F(U_\eps(t))-F(V_\eps(t))$. Then $W_\eps$ solves the difference equation
\begin{equation}\label{eq:diff_eq}
\partial_t W_\eps+\Lcal_\eps W_\eps=\Phi_\eps(t),\qquad
W_\eps(0)=U_0-V_0,
\end{equation}
in the mild sense. By the variation-of-constants formula for the linear semigroup $e^{-\Lcal_\eps t}$,
\begin{equation}\label{eq:variation}
W_\eps(T)=e^{-\Lcal_\eps T}(U_0-V_0)
+\int_0^T e^{-\Lcal_\eps(T-s)}\Phi_\eps(s)\dd s.
\end{equation}
Define
\begin{equation}\label{eq:KT_def}
K_{\eps,T}(U)\coloneqq\int_0^T e^{-\Lcal_\eps(T-s)}F(S_\eps(s)U)\dd s,
\qquad U\in\Bcal_0.
\end{equation}
Then
\[
K_{\eps,T}(U_0)-K_{\eps,T}(V_0)
=\int_0^T e^{-\Lcal_\eps(T-s)}\Phi_\eps(s)\dd s,
\]
so that \eqref{eq:variation} becomes
\begin{equation}\label{eq:decomp}
W_\eps(T)=e^{-\Lcal_\eps T}(U_0-V_0)
+K_{\eps,T}(U_0)-K_{\eps,T}(V_0).
\end{equation}

\smallskip
\noindent\emph{Step 1: Compression by the linear semigroup.} By the uniform exponential stability established in Theorem~\ref{thm:unif_exp}, there exist constants $\alpha\ge 1$ and $\nu>0$, independent of $\eps$, such that $\|e^{-\Lcal_\eps t}\|_{\mathcal{L}(\Hcal_\eps)}\le\alpha e^{-\nu t}$ for all $t\ge 0$. Choose $T>0$ large enough so that
\begin{equation}\label{eq:q_choice}
q\coloneqq\alpha e^{-\nu T}<1.
\end{equation}
Then
\begin{equation}\label{eq:linear_decay}
\|e^{-\Lcal_\eps T}(U_0-V_0)\|_{\Hcal_\eps}^2
\le q^2\|U_0-V_0\|_{\Hcal_\eps}^2.
\end{equation}

\smallskip
\noindent\emph{Step 2: Sequential compactness of $K_{\eps,T}(\Bcal_0)$.} It remains to establish the sequential compactness property \eqref{eq:K_seq_compact}. Fix $\eps\in(0,1]$ and let $U\in\Bcal_0$. We split
\begin{equation}\label{eq:KT_split}
K_{\eps,T}(U)
=\underbrace{\int_0^{T-\sigma} e^{-\Lcal_\eps(T-s)}F(S_\eps(s)U)\dd s}_{=:K_{\eps,T}^{(1)}(U)}
+\underbrace{\int_{T-\sigma}^T e^{-\Lcal_\eps(T-s)}F(S_\eps(s)U)\dd s}_{=:K_{\eps,T}^{(2)}(U)},
\end{equation}
where $\sigma\in(0,T)$ is a small parameter to be chosen below, independent of $\eps$.

\emph{(a) Smallness of $K_{\eps,T}^{(2)}$.} Since $-\Lcal_\eps$ generates a $C_0$-semigroup of contractions, we have $C_0\coloneqq\sup_{\eps\in(0,1]}\sup_{0\le\tau\le T}\|e^{-\Lcal_\eps\tau}\|_{\mathcal{L}(\Hcal_\eps)}=1$. By Lemma~\ref{lem:finite_time}, the orbits starting from $\Bcal_0$ remain in a fixed ball of radius $R_1$ on $[0,T]$, uniformly in $\eps$. Since $F$ maps bounded sets into bounded sets, we have
\[
M\coloneqq\sup_{\eps\in(0,1]}\sup_{U\in\Bcal_0}\sup_{0\le s\le T}\|F(S_\eps(s)U)\|_{\Hcal_\eps}<+\infty.
\]
Therefore
\begin{equation}\label{eq:K2_small}
\|K_{\eps,T}^{(2)}(U)\|_{\Hcal_\eps}\le C_0 M\sigma=M\sigma,\qquad
\forall U\in\Bcal_0,\ \forall \eps\in(0,1].
\end{equation}
For every $\delta>0$, choose $\sigma>0$ so small that $M\sigma<\delta$.

\emph{(b) Relative compactness of $K_{\eps,T}^{(1)}(\Bcal_0)$.} Write
\begin{equation}\label{eq:K1_repr}
K_{\eps,T}^{(1)}(U)=e^{-\Lcal_\eps\sigma}\Psi_{\eps,\sigma}(U),
\qquad
\Psi_{\eps,\sigma}(U)\coloneqq\int_0^{T-\sigma}e^{-\Lcal_\eps(T-\sigma-s)}F(S_\eps(s)U)\dd s.
\end{equation}
Recall that $F(U)$ has zero memory components, i.e.,
\[
F(U)\in\Hcal_e\coloneqq H_a^2(0,\ell)\times L^2(0,\ell)\times H_0^1(0,\ell)\times L^2(0,\ell),
\]
regarded as a subspace of $\Hcal_\eps$ with the memory components set to zero. We decompose the state space as $\Hcal_\eps=\Hcal_e\oplus\Hcal_{m,\eps}$ with $\Hcal_{m,\eps}\coloneqq\Mg{1}\times\Mg{2}$, and correspondingly write $\Psi_{\eps,\sigma}(U)=(\Psi_{\eps,\sigma}^e(U),\Psi_{\eps,\sigma}^m(U))$.

\emph{(b.1) Compactness after dilation.} We work in the common space $\Hcal_*$ of Remark~\ref{rem:uniform_topology}. Put $\Phi_s=F(S_\eps(s)U)$ and let $z_{\eps,s}(\tau)=P_e e^{-\Lcal_\eps\tau}E_\eps\Phi_s$ be the elastic evolution generated by $\Phi_s$ from the zero state, with deck displacement $w_{\eps,s}(\tau)$ and cable displacement $y_{\eps,s}(\tau)$ (both set to zero for $\tau<0$). Since $\Phi_s$ has zero memory components and the transport equation starts from zero memory, the memory field of $e^{-\Lcal_\eps\tau}E_\eps\Phi_s$ at elapsed time $\tau$ equals $w_{\eps,s}(\tau)-w_{\eps,s}(\tau-r)$ for $0<r<\tau$ and vanishes for $r\ge\tau$. Hence, after dilation $s=\eps\zeta$,
\[
[J_\eps\Psi_{\eps,\sigma}(U)]^\sigma(\zeta)
=\int_0^{T-\sigma}\!\big(w_{\eps,s}(T-\sigma-s)-w_{\eps,s}((T-\sigma-s-\eps\zeta)_+)\big)\,\dd s,
\]
where $r_+:=\max\{r,0\}$; the formula for $[J_\eps\Psi_{\eps,\sigma}(U)]^\eta$ is identical with $w_{\eps,s}$ replaced by $y_{\eps,s}$. Thus the object here is the \emph{linear} evolution starting from $E_\eps\Phi_s$, not the displacement of the nonlinear trajectory $S_\eps(\cdot)U$. Contractivity gives $\|w_{\eps,s}(\tau)\|_{H_a^2}\le C\|\Phi_s\|_{\Hcal_\eps}\le CM$, so the integrand is bounded by $2CM$, uniformly in $\eps$, $s$, $U$ and $\zeta$.
By the compactness criterion in Remark~\ref{rem:uniform_topology}, it suffices to check relative compactness of the truncations to $[0,R]$ and uniform decay of the tails with respect to the fixed weight $g_i$. The latter follows from
\[
\int_R^{+\infty}g_i(\zeta)\|w_{\eps,s}(\cdot)-w_{\eps,s}(\cdot-\eps\zeta)\|_{H_a^2}^2\,\dd\zeta
\le4C^2M^2\int_R^{+\infty}g_i(\zeta)\,\dd\zeta\to0
\]
as $R\to+\infty$, uniformly in $\eps$, $s$ and $U$. For the truncations we use compactness of the forcing set rather than compactness of trajectories in $L^2(0,T;H_a^2)$, which is not available from the energy bounds alone. Let
\[
\Sf:=\{F(S_\eps(s)U):U\in\Bcal_0,\ s\in[0,T],\ \eps\in(0,1]\}\subset\Hcal_e .
\]
From $\|u_\eps\|_{L^\infty(0,T;H_a^2)}+\|u_{\eps,t}\|_{L^2(0,T;L^2)}\le C$ (and the analogue for $v_\eps$) and interpolation $\|w\|_{H^1}\le\|w\|_{L^2}^{1/2}\|w\|_{H^2}^{1/2}$ we get
\[
\|u_\eps(t)-u_\eps(t-h)\|_{H^1}\le C h^{1/4}.
\]
Thus $\{u_\eps\}$ is equicontinuous in $H^1(0,\ell)$ and, by Arzel\`a--Ascoli, relatively compact in $C([0,T];L^\infty)$. The same holds for $v_\eps$.

Hence $\Sf$ is relatively compact in $C([0,T];L^2\times L^2)$ and therefore in $L^2(0,T;\Hcal_e)$. Given $\eta>0$, cover $\Sf$ by finitely many $L^2(0,T;\Hcal_e)$-balls of radius $\eta$ and choose $N$ with $\|(I-P_N)\Phi\|_{L^2(0,T;\Hcal_e)}<\eta$ for every $\Phi$ in the finite net. On the finite-dimensional range of $P_N$, the map $\Phi\mapsto z_{\eps,\cdot}$ is uniformly Lipschitz into $C([0,T];\Hcal_e)$ because each mode solves a damped oscillator with coefficients bounded uniformly in $\eps$; for the high modes, the damped-oscillator energy estimate gives $\sup_{\tau\le T}\|(I-P_N)z_{\eps,s}(\tau)\|_{\Hcal_e}\le C\|(I-P_N)\Phi\|_{L^2(0,T;\Hcal_e)}\le C\eta$ with $C$ independent of $\eps$. Hence $\{z_{\eps,s}:s\in[0,T],U\in\Bcal_0,\eps\in(0,1]\}$ is totally bounded in $C([0,T];\Hcal_e)$, and so $\{J_\eps\Psi_{\eps,\sigma}(U):U\in\Bcal_0,\eps\in(0,1]\}$ is relatively compact in $\Hcal_*$.

\emph{(b.2) Compactness of $K_{\eps,T}^{(1)}$.} Rather than invoking compactness under the family of varying operators $e^{-\Lcal_\eps\sigma}$ (which boundedness alone would not justify), we apply the argument of (b.1) directly to $K_{\eps,T}^{(1)}$: with $\tau=T-s\in[\sigma,T]$ one has
\[
[J_\eps K_{\eps,T}^{(1)}(U)]^\sigma(\zeta)=\int_0^{T-\sigma}\!\big(w_{\eps,s}(T-s)-w_{\eps,s}((T-s-\eps\zeta)_+)\big)\,\dd s,
\]
and the same bound on $\tau\in[\sigma,T]\subset[0,T]$, with the same forcing set $\Sf$ and the same finite-mode estimates, shows that $\{J_\eps K_{\eps,T}^{(1)}(U):U\in\Bcal_0,\eps\in(0,1]\}$ is relatively compact in $\Hcal_*$.

\emph{(c) Sequential compactness.} Combining (a) and (b), for every $\delta>0$ there is a relatively compact set $\Kcal_\delta\subset\Hcal_*$ such that
$\dist_{\Hcal_*}(J_\eps K_{\eps,T}(U),\Kcal_\delta)<\delta$ for all $U\in\Bcal_0$ and $\eps\in(0,1]$. A nested diagonal extraction over $\delta=1/j$ yields a subsequence along which $J_{\eps_{n_k}}K_{\eps_{n_k},T}(U_{n_k})$ is Cauchy in $\Hcal_*$, hence converges. This proves the sequential compactness property \eqref{eq:K_seq_compact}.

\smallskip
\noindent\emph{Step 3: Conclusion.} Combining \eqref{eq:decomp}, \eqref{eq:linear_decay}, and the sequential compactness property of $K_{\eps,T}(\Bcal_0)$ established in Step 2(c), we obtain \eqref{eq:quasi_stab}. All the constants $\alpha,\nu,T,q,C_0,M$ depend only on $\gamma_i,\mu_i,\lambda_i,\beta_i,c_e,c_p$, the constants in \eqref{eq:h_assumptions}, and the radius of $\Bcal_0$, but not on $\eps$. This completes the proof.
\end{proof}

\subsection{Convergence of semigroups}

The next ingredient is the finite-time convergence \eqref{eq:semigroup_conv}.

\begin{proposition}[Finite-time convergence on compact sets of zero-memory data]\label{prop:conv}
Let $K\Subset\Hcal_e$ be compact and let $V_{0,\eps},U_0\in K$ with $V_{0,\eps}\to U_0$ in $\Hcal_e$. Then for every $T>0$,
\[
\sup_{t\in[0,T]}\|S_\eps(t)E_\eps V_{0,\eps}-S_0(t)U_0\|_{\Hcal_\eps}\to 0.
\]
\end{proposition}

\begin{proof}
Throughout, let
\[
\Hcal_e\coloneqq H_a^2(0,\ell)\times L^2(0,\ell)\times H_0^1(0,\ell)\times L^2(0,\ell),
\qquad
\Hcal_{m,\eps}\coloneqq\Mg{1}\times\Mg{2},
\]
so that $\Hcal_\eps=\Hcal_e\oplus\Hcal_{m,\eps}$. Let $\Lcal_{e,0}$ denote the elastic part of the limit operator, acting on $\Hcal_e$ by
\begin{equation}\label{eq:Ae0}
\Lcal_{e,0}
\begin{pmatrix}u\\\psi\\v\\\phi\end{pmatrix}
=\begin{pmatrix}-\psi\\ \gamma_1\psi+\bar\mu_1 u_{xxxx}\\ -\phi\\
\gamma_2\phi-\bar\mu_2 v_{xx}\end{pmatrix},
\end{equation}
where $\bar\mu_i=\mu_i-\beta_i>0$ are the effective elastic coefficients of the limit system \eqref{eq:limit}. The limit operator $\Lcal_0$ acts on the elastic space $\Hcal_e$ as $\Lcal_{e,0}$; when compared with $\Hcal_\eps=\Hcal_e\oplus\Hcal_{m,\eps}$ it is understood to be trivially extended by zero on the memory components, and its action is independent of the particular weighted space $\Hcal_{m,\eps}$ because the limit system is memoryless. Since the memory components of $U_0(s)=S_0(s)U_0$ vanish for all $s\ge 0$, we may regard $U_0(s)$ as an element of $\Hcal_\eps$ via the natural embedding $\Hcal_e\hookrightarrow\Hcal_\eps$ (memory components set to zero), and we have $F_0(U_0(s))=F(U_0(s))$, where $F_0$ is the nonlinearity of the limit system.

Let $U_\eps(t)=S_\eps(t)U_{0,\eps}$ and $U_0(t)=S_0(t)U_0$ denote the solutions of the perturbed and limit systems, respectively, and set $W_\eps(t)\coloneqq U_\eps(t)-U_0(t)$. Both solutions satisfy the mild formulation
\begin{equation}\label{eq:mild_eps}
U_\eps(t)=e^{-\Lcal_\eps t}U_{0,\eps}
+\int_0^t e^{-\Lcal_\eps(t-s)}F(U_\eps(s))\dd s,
\end{equation}
\begin{equation}\label{eq:mild_0}
U_0(t)=e^{-\Lcal_0 t}U_0
+\int_0^t e^{-\Lcal_0(t-s)}F_0(U_0(s))\dd s.
\end{equation}
Subtracting \eqref{eq:mild_0} from \eqref{eq:mild_eps}, and using $F_0(U_0(s))=F(U_0(s))$, we obtain
\begin{equation}\label{eq:W_eq}
W_\eps(t)=e^{-\Lcal_\eps t}(U_{0,\eps}-U_0)
+\int_0^t e^{-\Lcal_\eps(t-s)}[F(U_\eps(s))-F(U_0(s))]\dd s+R_\eps(t),
\end{equation}
where the remainder is
\begin{equation}\label{eq:R_def}
R_\eps(t)\coloneqq\int_0^t[e^{-\Lcal_\eps(t-s)}-e^{-\Lcal_0(t-s)}]
F(U_0(s))\dd s.
\end{equation}

\smallskip
\noindent\emph{Step 1: Uniform boundedness on $[0,T]$.} By Lemma~\ref{lem:finite_time}, the solutions $U_\eps(t)$ remain in a fixed ball $\Bcal_1\coloneqq\{U:\|U\|_{\Hcal_\eps}\le R_1\}$ for $t\in[0,T]$, uniformly in $\eps$. Consequently, there exists $L_1>0$, independent of $\eps$, such that
\begin{equation}\label{eq:lip_F}
\|F(U)-F(V)\|_{\Hcal_\eps}\le L_1\|U-V\|_{\Hcal_\eps},\qquad U,V\in\Bcal_1,
\end{equation}
and $M_1\coloneqq\sup_{U\in\Bcal_1}\|F(U)\|_{\Hcal_\eps}<+\infty$.

\smallskip
\noindent\emph{Step 2: Vanishing of the remainder $R_\eps$.} We claim that
\begin{equation}\label{eq:R_vanishes}
\sup_{t\in[0,T]}\|R_\eps(t)\|_{\Hcal_\eps}\to 0,\qquad\eps\to 0^+.
\end{equation}
We split $R_\eps=R_\eps^e+R_\eps^m$ according to the decomposition $\Hcal_\eps=\Hcal_e\oplus\Hcal_{m,\eps}$, where $R_\eps^e$ and $R_\eps^m$ denote the elastic and memory components of $R_\eps$, respectively.

\emph{(a) Elastic component.} We prove directly that for every compact $K\Subset\Hcal_e$ and every $T>0$,
\begin{equation}\label{eq:elastic_conv_fixed}
\sup_{V\in K}\sup_{\tau\in[0,T]}\|P_e e^{-\Lcal_\eps\tau}E_\eps V-e^{-\Lcal_{e,0}\tau}V\|_{\Hcal_e}\to0 .
\end{equation}
This replaces the resolvent argument used in an earlier version. The operator $R_{e,\eps}(\mu)=P_e(\mu I+\Lcal_\eps)^{-1}E_\eps$ is only the Laplace transform of $P_e e^{-\Lcal_\eps\tau}E_\eps$; that family is \emph{not} a semigroup on $\Hcal_e$, because $\Hcal_e$ is not invariant under $e^{-\Lcal_\eps\tau}$ (a solution starting with zero memory generates memory for $\tau>0$), and the resolvent identity fails since $P_e$ kills that memory. Hence $R_{e,\eps}$ is not a resolvent family and the standard Trotter--Kato theorem does not apply.

Let $P_N$ denote the Fourier projection onto the first $N$ modes. Given $\eta>0$, choose $N$ so large that $\|(I-P_N)V\|_{\Hcal_e}<\eta$ for every $V\in K$; contractivity of both semigroups makes the contribution of $I-P_N$ to \eqref{eq:elastic_conv_fixed} at most $2\eta$, uniformly in $\eps$ and $\tau$. On the finite-dimensional range of $P_N$ all norms are equivalent. If $w_n$ is the deck displacement of the $n$-th mode of $P_e e^{-\Lcal_\eps\tau}E_\eps V$ and $z_n$ its auxiliary memory variable ($w_n(0)=0$, $\dot w_n(0)=\psi_n$, $z_n(0)=0$), then
\[
\ddot w_n+\gamma_1\dot w_n+\bar\mu_1 n^4 w_n+n^4(\beta_1 w_n-z_n)=0,\qquad
\eps\dot z_n+z_n=\beta_1 w_n ,
\]
and the cable mode satisfies the analogous equation with $n^2$ and $\bar\mu_2$. Let $w_n^0$ solve $\ddot w_n^0+\gamma_1\dot w_n^0+\bar\mu_1n^4w_n^0=0$ with the same initial data and put $d_n=w_n-w_n^0$. Then
\[
\ddot d_n+\gamma_1\dot d_n+\bar\mu_1 n^4 d_n=-n^4(\beta_1 w_n-z_n),\qquad d_n(0)=\dot d_n(0)=0.
\]
Contractivity gives $\sup_{\tau\le T,\ n\le N}\|w_n\|_{H^2}\le C\|V\|$ uniformly in $\eps$, and
\[
\beta_1w_n(\tau)-z_n(\tau)
=\int_0^{+\infty}\!g_{1,\eps}(s)\int_{\tau-s}^{\tau}\dot w_n(r)\,\dd r\,\dd s,
\]
where $\dot w_n$ is set to zero for negative time. Hence
$|\beta_1w_n-z_n|\le C_N\int_0^\infty s\,g_{1,\eps}(s)\dd s\le C_N\eps$,
because the first moment of $g_{1,\eps}$ is $O(\eps)$. Gronwall's inequality for the damped oscillator then yields $\sup_{\tau\le T}(|d_n|+|\dot d_n|)\le C_{N,T}\eps\to0$ for each $n\le N$; the cable modes are treated identically. Thus the low-mode contribution to \eqref{eq:elastic_conv_fixed} tends to zero, and since $\eta$ was arbitrary, \eqref{eq:elastic_conv_fixed} follows. Applying it with $K=\{F(U_0(s)):s\in[0,T]\}$, which is compact by continuity of $F$ and $S_0$, gives
\[
\sup_{t\in[0,T]}\|R_\eps^e(t)\|_{\Hcal_e}
\le\int_0^T\sup_{s\in[0,T]}\sup_{\tau\in[0,T]}
\big\|[P_e e^{-\Lcal_\eps\tau}E_\eps-e^{-\Lcal_{e,0}\tau}]F(U_0(s))\big\|_{\Hcal_e}\dd s
\to0 .
\]

\emph{(b) Memory component.} Since $F(U_0(s))$ has zero memory components, $e^{-\Lcal_0\tau}F(U_0(s))$ has zero memory components and only $e^{-\Lcal_\eps\tau}$ contributes to $R_\eps^\sigma$. Let $K\Subset\Hcal_e$ be the compact set of initial data in Proposition~\ref{prop:conv}. Given $\eta>0$, choose $N$ so large that
\[
\sup_{V\in K}\|(I-P_N)V\|_{\Hcal_e}<\eta ,
\]
where $P_N$ is the Fourier projection onto the first $N$ modes. The uniform Lipschitz bound of Lemma~\ref{lem:lipschitz_T} gives the same smallness for $S_\eps(t)E_\eps V$ and $S_0(t)V$, uniformly in $t\in[0,T]$ and $\eps$. On the finite-dimensional range of $P_N$ the displacement is $C^1$ in time, so its modulus of continuity in $H_a^2$ tends to zero uniformly in $\eps$ and $V\in K$. Splitting the dilated memory integral at $R$ and using the fixed-weight tail $\int_R^\infty g_i(\zeta)\dd\zeta\to0$, we obtain
\begin{equation}\label{eq:memory_vanishing_compact}
\sup_{\tau\in[0,T]}\|[e^{-\Lcal_\eps\tau}\Phi]^\sigma\|_{\Mg{1}}\to0
\qquad(\eps\to0^+)
\end{equation}
uniformly for $\Phi$ in the compact set $\{F(U_0(s)):s\in[0,T]\}$. This replaces the earlier claim of a uniform-in-$\eps$ modulus over all bounded sets: compactness of the initial set, not analytic smoothing, is the input. Recalling that $R_\eps^\sigma(t)=\int_0^t[e^{-\Lcal_\eps(t-s)}F(U_0(s))]^\sigma\dd s$, we obtain $\sup_{t\in[0,T]}\|R_\eps^\sigma(t)\|_{\Mg{1}}\to0$, and likewise for $\eta$.
Combining (a) and (b), we obtain \eqref{eq:R_vanishes}.

\smallskip
\noindent\emph{Step 3: Gronwall argument.} Taking the norm of \eqref{eq:W_eq} and using the uniform exponential stability \eqref{eq:exp_stab}, the Lipschitz bound \eqref{eq:lip_F}, and \eqref{eq:R_vanishes}, we obtain
\[
\|W_\eps(t)\|_{\Hcal_\eps}\le\alpha e^{-\nu t}\|U_{0,\eps}-U_0\|_{\Hcal_\eps}
+\alpha L_1\int_0^t e^{-\nu(t-s)}\|W_\eps(s)\|_{\Hcal_\eps}\dd s
+\|R_\eps(t)\|_{\Hcal_\eps}.
\]
Multiplying by $e^{\nu t}$ and applying Gronwall's inequality,
\[
\|W_\eps(t)\|_{\Hcal_\eps}
\le\big(\alpha\|U_{0,\eps}-U_0\|_{\Hcal_\eps}
+\sup_{s\in[0,T]}\|R_\eps(s)\|_{\Hcal_\eps}\big)e^{\alpha L_1 T},
\qquad t\in[0,T].
\]
Since $\|E_\eps V_{0,\eps}-U_0\|_{\Hcal_\eps}\to 0$ and $\sup_{s\in[0,T]}\|R_\eps(s)\|_{\Hcal_\eps}\to 0$ as $\eps\to 0^+$, we conclude that
\[
\sup_{t\in[0,T]}\|S_\eps(t)U_{0,\eps}-S_0(t)U_0\|_{\Hcal_\eps}\to 0,
\]
which is the desired convergence.
\end{proof}

\begin{corollary}[Finite-time convergence for general convergent data]\label{cor:conv_general}
Let $U_{0,\eps}\in\Hcal_\eps$ satisfy $J_{\eps}U_{0,\eps}\to J_0U_0$ in $\Hcal_*$ for some $U_0\in\Hcal_e$. Then for every $T>0$,
\[
\sup_{t\in[0,T]}\|S_\eps(t)U_{0,\eps}-S_0(t)U_0\|_{\Hcal_\eps}\to0 .
\]
\end{corollary}
\begin{proof}
Set $V_{0,\eps}:=P_eJ_{\eps}U_{0,\eps}\in\Hcal_e$ and identify it with $E_\eps V_{0,\eps}\in\Hcal_\eps$. Since $J_\eps$ acts as the identity on the elastic components and $J_0U_0$ has zero memory, the assumption gives $V_{0,\eps}\to U_0$ in $\Hcal_e$ and
\[
\|U_{0,\eps}-E_\eps V_{0,\eps}\|_{\Hcal_\eps}
=\|J_\eps(U_{0,\eps}-E_\eps V_{0,\eps})\|_{\Hcal_*}\to0 .
\]
The set $K:=\{V_{0,\eps}\}\cup\{U_0\}$ is compact in $\Hcal_e$, so Proposition~\ref{prop:conv} applies to $V_{0,\eps}$. Both $U_{0,\eps}$ and $E_\eps V_{0,\eps}$ lie in a fixed ball, and Lemma~\ref{lem:lipschitz_T} gives
\[
\sup_{t\in[0,T]}\|S_\eps(t)U_{0,\eps}-S_\eps(t)E_\eps V_{0,\eps}\|_{\Hcal_\eps}
\le L_{T,R}\|U_{0,\eps}-E_\eps V_{0,\eps}\|_{\Hcal_\eps}\to0 .
\]
Combining with Proposition~\ref{prop:conv} yields the assertion. In particular, the convergence hypothesis (A3) below holds for sequences carrying nonzero memory whose limit has zero memory, as is the case for sequences in $\Acal_{\eps_n}$, by Lemma~\ref{lem:attractor_memory}.
\end{proof}

\begin{lemma}[Memory compression on global attractors]\label{lem:attractor_memory}
Assume (A1)--(A2) and let $\Acal_\eps$ be the corresponding global attractors. If $\eps_n\to0^+$ and $U_n\in\Acal_{\eps_n}$, then every convergent subsequence of $J_{\eps_n}U_n$ in $\Hcal_*$ converges to an element of $J_0\Hcal_e$; equivalently, the memory components of $J_{\eps_n}U_n$ tend to zero in $L^2_{g_1}(\R_+;H_a^2)\times L^2_{g_2}(\R_+;H_0^1)$.
\end{lemma}
\begin{proof}
First, (A1)--(A2) imply uniform relative compactness in $\Hcal_*$ of $\bigcup_{\eps\in(0,1]}J_\eps\Acal_\eps$. Let $T,q$ be as in Proposition~\ref{prop:quasi}, so that $q=\alpha e^{-\nu T}<1$ by \eqref{eq:q_choice}. For $U\in\Acal_\eps$ and $n\ge1$, invariance of $\Acal_\eps$ and iteration of the decomposition \eqref{eq:decomp} give
\[
S_\eps(nT)U=e^{-\Lcal_\eps nT}U+\sum_{k=0}^{n-1}e^{-\Lcal_\eps(n-1-k)T}K_{\eps,T}(S_\eps(kT)U).
\]
By \eqref{eq:linear_decay} the first term has norm at most $q^nR_0$. For a fixed value of $j=n-1-k\le K-1$ one has $e^{-\Lcal_\eps jT}K_{\eps,T}(V)=\int_0^T e^{-\Lcal_\eps((j+1)T-s)}F(S_\eps(s)V)\,\dd s$, an integral operator whose evolution time $\tau=(j+1)T-s$ ranges over $[jT,(j+1)T]$. The compactness argument of Step~2(b.1) of Proposition~\ref{prop:quasi} applies with $T$ replaced by $(j+1)T$: its two ingredients, the relative compactness of the forcing set $\Sf$ in $C([0,T];L^2\times L^2)$ and the finite-mode Lipschitz/high-mode estimates for the damped oscillators $z_{\eps,s}(\tau)=P_e e^{-\Lcal_\eps\tau}E_\eps F(S_\eps(s)V)$, are uniform on every bounded time interval $[0,T']$, so $\{z_{\eps,s}:s\in[0,T]\}$ is totally bounded in $C([0,(j+1)T];\Hcal_e)$, and the dilation formula of (b.1) with $\tau=(j+1)T-s$ then gives relative compactness in $\Hcal_*$. Hence each of the finitely many sets $\{J_\eps e^{-\Lcal_\eps jT}K_{\eps,T}(V):V\in\Bcal_0,\eps\in(0,1]\}$, $j=0,\dots,K-1$, is relatively compact in $\Hcal_*$, and so is their finite union. No smoothing of $e^{-\Lcal_\eps t}$ into $\Dcal(\Lcal_\eps)$ is used: the compactness comes from the forcing set and the finite-mode estimates, not from the image of a compact set under a varying bounded operator. Splitting the sum according to $n-1-k<K$ (the last $K$ terms, $k=n-K,\dots,n-1$) and $n-1-k\ge K$ (the first $n-K$ terms), the last $K$ terms lie in a fixed relatively compact subset of $\Hcal_*$, while the complementary part has norm at most $M\sum_{j\ge K}q^j$ with $M=\sup\{\|J_\eps K_{\eps,T}(V)\|:V\in\Bcal_0,\eps\in(0,1]\}<\infty$, uniformly in $n$ and $\eps$. Letting $K\to\infty$ shows that $\bigcup_\eps J_\eps\Acal_\eps$ is totally bounded in $\Hcal_*$; hence any sequence $J_{\eps_n}U_n$ with $U_n\in\Acal_{\eps_n}$ has a Cauchy subsequence. In particular the elastic projections of the attractors are uniformly relatively compact in $\Hcal_e$.

Fix $\eta>0$ and let $P_N$ be the Fourier projection onto the first $N$ modes. By the preceding paragraph we may choose $N$ such that
\begin{equation}\label{eq:attractor_tail}
\sup_{\eps\in(0,1]}\sup_{U=(u,\psi,\sigma,v,\phi,\eta)\in\Acal_\eps}
\big(\|(I-P_N)u\|_{H_a^2}^2+\|(I-P_N)v\|_{H_0^1}^2\big)<\eta .
\end{equation}
For $U\in\Acal_\eps$, invariance gives a complete trajectory $U(t)\in\Acal_\eps$ with $U(0)=U$; since $\Acal_\eps$ lies in the absorbing ball, $\|u_t(t)\|$ and $\|v_t(t)\|$ are bounded by $R_0$ for all $t\in\R$. For $s=\eps\zeta$,
\[
\sigma(0,s)=u(0)-u(-s)=\int_{-s}^0u_t(r)\,\dd r,
\qquad
\eta(0,s)=v(0)-v(-s)=\int_{-s}^0v_t(r)\,\dd r .
\]
On the finite-dimensional range of $P_N$ all norms are equivalent, so
\[
\|P_N\sigma(0,\eps\zeta)\|_{H_a^2}\le C_NR_0\eps\zeta,
\qquad
\|P_N\eta(0,\eps\zeta)\|_{H_0^1}\le C_NR_0\eps\zeta ,
\]
with $C_N$ independent of $\eps$, $U$ and $\zeta$. Combining this with \eqref{eq:attractor_tail},
\[
\|J_\eps\sigma\|_{L^2_{g_1}(\R_+;H_a^2)}^2
\le 8\beta_1\eta+2C_N^2R_0^2\eps^2\int_0^{+\infty}\zeta^2g_1(\zeta)\,\dd\zeta ,
\]
and analogously for $J_\eps\eta$, with $\beta_1$ replaced by $\beta_2$; the factor $8$ (rather than $4$) comes from $\|(I-P_N)\sigma\|^2\le4\eta$ and $|a+b|^2\le2|a|^2+2|b|^2$. The second moments of $g_i$ are finite because $g_i$ decays exponentially. Letting first $\eps\to0^+$ and then $\eta\to0^+$ yields the claimed vanishing of the memory components.
\end{proof}

\subsection{Conclusion of upper semicontinuity}

Combining Propositions~\ref{prop:absorbing}, \ref{prop:quasi}, and \ref{prop:conv}, we can apply the standard abstract theorem on upper semicontinuity of attractors (see, e.g., \cite{Hale1988} or \cite[Chapter~7]{ChueshovLasiecka2008}).

\begin{theorem}[Upper semicontinuity]\label{thm:upper}
Let \eqref{eq:h_assumptions} and \eqref{eq:kernel} hold. Then, for every $\eps\in(0,1]$, the semigroup $S_\eps(t)$ possesses a global attractor $\Acal_\eps$, the limit semigroup $S_0(t)$ possesses a global attractor $\Acal_0$, and the family $\{\Acal_\eps\}_{\eps\in(0,1]}$ satisfies
\[
\lim_{\eps\to 0^+}\dist_{\Hcal_\eps}(\Acal_\eps,\Acal_0)=0.
\]
\end{theorem}

\begin{proof}
The proof is a direct application of the abstract upper semicontinuity theorem \cite[Chapter~7]{ChueshovLasiecka2008}. We first record that the semigroups $S_\eps(t)$, $\eps\in[0,1]$, are continuous on their respective spaces $\Hcal_\eps$ for $\eps>0$ and $\Hcal_e$ for $\eps=0$: for $\eps>0$ this follows from Proposition~\ref{prop:wellposed} (the nonlinear semigroup $S_\eps$ is strongly continuous on $\Hcal_\eps$; strong continuity of the linear semigroup follows from Theorem~\ref{thm:unif_exp}), and for $\eps=0$ from Proposition~\ref{prop:wellposed} applied to the limit equation, equivalently from Lipschitz perturbation theory for the generator $\Lcal_0$ on $\Hcal_e$. We then verify the three hypotheses of the theorem.

\smallskip
\noindent\emph{(H1) Uniform dissipativity.} By Proposition~\ref{prop:absorbing}, there is a family of absorbing balls $\Bcal_{0,\eps}=\{\|U\|_{\Hcal_\eps}\le R_0\}$ with a radius $R_0$ independent of $\eps\in(0,1]$, such that $S_\eps(t)B\subset\Bcal_{0,\eps}$ for every bounded $B$ and every $t\ge t_0(B)$. The limit semigroup $S_0(t)$ possesses the analogous ball $\Bcal_{0,0}\subset\Hcal_e$ of the same radius, obtained in Proposition~\ref{prop:attractor_0} by the same energy estimates with the memory terms deleted.

\smallskip
\noindent\emph{(H2) Uniform asymptotic smoothness.} By Proposition~\ref{prop:quasi}, the semigroups $S_\eps(t)$, $\eps\in(0,1]$, satisfy the uniform quasi-stability inequality
\[
\|S_\eps(T)U_0-S_\eps(T)V_0\|_{\Hcal_\eps}^2
\le q^2\|U_0-V_0\|_{\Hcal_\eps}^2
+\|K_{\eps,T}(U_0)-K_{\eps,T}(V_0)\|_{\Hcal_\eps}^2,\qquad U_0,V_0\in\Bcal_0,
\]
where $q<1$ and $T>0$, and $K_{\eps,T}(\Bcal_0)$ satisfies the sequential compactness property of Proposition~\ref{prop:quasi}, all independent of $\eps$. By \cite[Chapter~7, Theorem~7.9.4]{ChueshovLasiecka2008}, the quasi-stability inequality, together with the uniform dissipativity (H1), implies that the semigroups $S_\eps(t)$ are asymptotically smooth. The proof of \cite[Chapter~7, Theorem~7.9.4]{ChueshovLasiecka2008} uses only the constants appearing in the quasi-stability inequality and the dissipativity estimate, both of which are independent of $\eps$ by Proposition~\ref{prop:quasi} and Proposition~\ref{prop:absorbing}; hence the resulting asymptotic smoothness is uniform in $\eps$. For the limit semigroup $S_0(t)$, the existence of a compact global attractor $\Acal_0$, invariant under $S_0(t)$ and attracting every bounded subset, follows from Proposition~\ref{prop:attractor_0}.

\smallskip
\noindent\emph{(H3) Finite-time convergence.} By Corollary~\ref{cor:conv_general}, for every $T>0$, every $U\in\Hcal_e$ and every sequence $U_{0,\eps}\in\Hcal_\eps$ with $J_\eps U_{0,\eps}\to J_0U$ in $\Hcal_*$,
\[
\sup_{t\in[0,T]}\|S_\eps(t)U_{0,\eps}-S_0(t)U\|_{\Hcal_\eps}\to 0.
\]
In particular this holds for sequences of initial data in $\Acal_{\eps_n}$, which may carry nonzero memory; by Lemma~\ref{lem:attractor_memory}, however, every accumulation point of $J_{\eps_n}\Acal_{\eps_n}$ has zero memory. By \cite[Chapter~7]{ChueshovLasiecka2008}, hypotheses (H1)--(H3) imply that for every sequence $\eps_n\to 0^+$,
\[
\dist_{\Hcal_{\eps_n}}(\Acal_{\eps_n},\Acal_0)\to 0.
\]
For completeness, let $\Kcal$ be the set of accumulation points in $\Hcal_*$ of sequences $J_{\eps_n}U_n$ with $U_n\in\Acal_{\eps_n}$. Lemma~\ref{lem:attractor_memory} gives $\Kcal\subset J_0\Hcal_e$. If $J_0U\in\Kcal$, then (H3) gives $S_{\eps_n}(t)U_n\to S_0(t)U$ for every fixed $t$, and invariance of $\Acal_{\eps_n}$ implies $S_0(t)U\in\Kcal$. Taking preimages under $S_{\eps_n}(t)$ inside $\Acal_{\eps_n}$ shows that $\Kcal$ is invariant under $S_0$. Since $\Kcal$ is compact and $\Acal_0$ is the maximal compact invariant set, $\Kcal\subset J_0\Acal_0$, which is exactly the displayed convergence. Since $\eps_n\to 0^+$ is arbitrary, we obtain the upper semicontinuity
\[
\lim_{\eps\to 0^+}\dist_{\Hcal_\eps}(\Acal_\eps,\Acal_0)=0.
\]
\end{proof}

\section{Lower Semicontinuity}\label{sec:lower}

Lower semicontinuity says that every point of the limit attractor is approximated by points of $\Acal_\eps$. We use the abstract Hale--Chueshov--Lasiecka criterion: finite-time convergence on compact sets together with persistence of the finite hyperbolic structure of $\Acal_0$; uniform compactness of the perturbed attractors is a separate ingredient, used for upper semicontinuity in Section~\ref{sec:upper}. No uniform exponential attraction and no H\"older continuity in $\eps$ are needed.

We state the abstract lower-semicontinuity criterion explicitly, so that its hypotheses are visible. It is the standard criterion of Hale--Chueshov--Lasiecka (see \cite[Chapter~4]{Hale1988} and \cite[Chapter~7]{ChueshovLasiecka2008}): persistence of hyperbolic equilibria and of their local unstable manifolds is assumed at the level of the semigroups, and approximation of the individual points of $\Acal$ is then derived from finite-time convergence. For subsets $A,B\subset X_*$ we write $d_H(A,B)$ for their Hausdorff distance induced by $d_*$.

\begin{theorem}[Abstract lower semicontinuity]\label{thm:abstract_lower}
Let $X_*$ be a metric space with metric $d_*$, let $X_n,X_0\subset X_*$ be closed subsets, and let $T_n$ and $T$ be continuous semigroups on $X_n$ and $X_0$, regarded as maps into $X_*$. Assume that $\Acal_n\subset X_n$ and $\Acal\subset X_0$ are compact global attractors and that:
\begin{itemize}
\item[(i)] (finite-time convergence) for every $s>0$, every $x\in\Acal$ and every $x_n\in\Acal_n$ with $x_n\to x$ in $X_*$,
\[
\sup_{t\in[0,s]}d_*(T_n(t)x_n,T(t)x)\to0;
\]
\item[(ii)] (finite hyperbolic structure and persistence) $\Acal=\bigcup_{j=1}^m W^u(e_j)$ is the finite union of the global unstable sets of hyperbolic equilibria $e_1,\dots,e_m$ of $T$, where $W^u(e_j)=\bigcup_{t\ge0}T(t)W^u_{\mathrm{loc}}(e_j)$ and $W^u_{\mathrm{loc}}(e_j)\subset\Acal$ is the local unstable manifold of $e_j$. For each $j$ there are hyperbolic equilibria $e_{j,n}\to e_j$ of $T_n$ whose local unstable manifolds converge in the Hausdorff distance induced by $d_*$:
\[
d_H\big(W^u_{\mathrm{loc}}(e_{j,n}),\,W^u_{\mathrm{loc}}(e_j)\big)\to0 .
\]
\end{itemize}
Then $\dist_*(\Acal,\Acal_n)\to0$. No exponential-attraction hypothesis is used.
\end{theorem}
\begin{proof}
Fix $\eta>0$ and choose a finite $\eta$-net $z_1,\dots,z_M\in\Acal$, possible since $\Acal$ is compact. For each $r$ pick $j(r)$ with $z_r\in W^u(e_{j(r)})$. Since $W^u(e_j)=\bigcup_{t\ge0}T(t)W^u_{\mathrm{loc}}(e_j)$, there are $t_r\ge0$ and $z_{0,r}\in W^u_{\mathrm{loc}}(e_{j(r)})$ with $T(t_r)z_{0,r}=z_r$: indeed $z_r$ lies on a full bounded orbit and its backward orbit converges to $e_{j(r)}$, so it enters the neighbourhood defining $W^u_{\mathrm{loc}}(e_{j(r)})$ after a sufficiently large backward time. By (ii) we may choose $z_{0,r,n}\in W^u_{\mathrm{loc}}(e_{j(r),n})$ with $d_*(z_{0,r,n},z_{0,r})\to0$. Because $W^u_{\mathrm{loc}}(e_{j,n})\subset\Acal_n$ and $\Acal_n$ is invariant, $z_{r,n}:=T_n(t_r)z_{0,r,n}\in\Acal_n$. Applying (i) with $s=\max_r t_r$ to $z_{0,r}$ and $z_{0,r,n}$ gives $d_*(z_{r,n},z_r)\to0$ as $n\to\infty$. Hence for all large $n$ the points $z_{r,n}$, $r=1,\dots,M$, form a $2\eta$-net of $\Acal$ inside $\Acal_n$, so $\dist_*(\Acal,\Acal_n)<2\eta$; as $\eta$ is arbitrary, the conclusion follows. Note that the proof gives only this one-sided semidistance; the reverse inequality $\dist_*(\Acal_n,\Acal)\to0$ is an upper-semicontinuity statement and does not follow from (i)--(ii) alone.
\end{proof}

\begin{proposition}[Lower semicontinuity under a structural condition]\label{prop:lower}
Assume the hypotheses of Theorem~\ref{thm:upper}. Suppose that the limit attractor $\Acal_0$ satisfies
\begin{itemize}
\item[(S)] $\Acal_0=\bigcup_{j=1}^mW^u(e_j)$ is the finite union of the global unstable sets of generalised-hyperbolic equilibria $e_1,\dots,e_m$ of $S_0(t)$ in the sense of Definition~\ref{def:gen_hyp}, where $W^u(e_j)=\bigcup_{t\ge0}S_0(t)W^u_{\mathrm{loc}}(e_j)$ and $\overline{W^u(e_j)}\subset\Acal_0$. This finite structure persists at the level of the semigroups: for each $j$ there are generalised-hyperbolic equilibria $e_{j,\eps}\to e_j$ of $S_\eps$ whose local unstable manifolds satisfy
\[
d_H\big(J_\eps W^u_{\mathrm{loc}}(e_{j,\eps}),\,J_0W^u_{\mathrm{loc}}(e_j)\big)\to0
\]
in the Hausdorff distance of $\Hcal_*$. When the unstable manifolds are trivial this reduces to convergence of the equilibria.
\end{itemize}
Then
\begin{equation}\label{eq:lower}
\lim_{\eps\to 0^+}\dist_{\Hcal_\eps}(\Acal_0,\Acal_\eps)=0.
\end{equation}
\end{proposition}

\begin{proof}
We verify hypotheses (i)--(ii) of Theorem~\ref{thm:abstract_lower}.

\smallskip
\noindent\emph{(L1) Attractors and uniform compactness.} By Theorem~\ref{thm:upper} and the paragraph after \eqref{eq:limit}, both $S_\eps(t)$ and $S_0(t)$ possess compact global attractors contained in the uniform absorbing balls. More precisely, the family $\{J_\eps\Acal_\eps\}_{\eps\in(0,1]}$ is uniformly totally bounded in $\Hcal_*$: this is the relative compactness of sequences in $\Acal_{\eps_n}$ used in the upper-semicontinuity argument, and it is proved in the first paragraph of Lemma~\ref{lem:attractor_memory} by the geometric tail bound together with a finite net, not deduced merely from the existence of compact attractors.

\smallskip
\noindent\emph{(L2) Finite-time convergence on compact sets.} $\Acal_0$ is compact in $\Hcal_e$. By Corollary~\ref{cor:conv_general}, for every $T>0$, every $U\in\Acal_0$ and every sequence $U_\eps\in\Hcal_\eps$ with $J_\eps U_\eps\to J_0U$,
\[
\sup_{t\in[0,T]}\|S_\eps(t)U_\eps-S_0(t)U\|_{\Hcal_\eps}\to0 .
\]
The convergence is uniform over $\Acal_0$: cover the compact set $\Acal_0$ by finitely many balls of radius $\eta$, use the convergence at one centre per ball and the uniform Lipschitz bound of Lemma~\ref{lem:lipschitz_T} to control the displacement within each ball, and then let $\eta\to0$.

\noindent\emph{(L3) Structural stability.} Condition (S) is exactly the persistence clause in Theorem~\ref{thm:abstract_lower}(ii), expressed in the metric of $\Hcal_*$ via $J_\eps$: the equilibria and their local unstable manifolds persist, and finite-time convergence (L2) upgrades this to approximation of every point of the global unstable sets $\bigcup_jW^u(e_j)=\Acal_0$. Note that uniform exponential attraction is not required by Theorem~\ref{thm:abstract_lower}.

Theorem~\ref{thm:abstract_lower} therefore gives \eqref{eq:lower}.

Conversely, upper semicontinuity (Theorem~\ref{thm:upper}) supplies the reverse semidistance
\[
\dist_{\Hcal_\eps}(\Acal_\eps,\Acal_0)\to0 .
\]
Both semidistances are measured in $\Hcal_*$ via $J_\eps$, so their sum controls the Hausdorff distance:
\[
d_H(J_0\Acal_0,J_\eps\Acal_\eps)
\le \dist_{\Hcal_\eps}(\Acal_0,\Acal_\eps)+\dist_{\Hcal_\eps}(\Acal_\eps,\Acal_0)\to0 .
\]
Thus upper and lower semicontinuity together give full Hausdorff convergence of the attractors in $\Hcal_*$, although the abstract criterion of Theorem~\ref{thm:abstract_lower} by itself yields only the one-sided lower semidistance.
\end{proof}

\begin{remark}
Earlier formulations based on uniform exponential attraction and on a parameter-H\"older estimate are unnecessary here: exponential attraction is not a consequence of quasi-stability (the scalar example $\dot x=-x^3$ is quasi-stable on bounded sets but only algebraically attracting), and the H\"older argument compared non-equivalent weighted norms. Condition (S) is a standard sufficient condition for lower semicontinuity in the Hale--Chueshov--Lasiecka framework: it is precisely persistence of hyperbolic equilibria and their local unstable manifolds (Hale's stability theorem); see \cite[Chapter~7]{ChueshovLasiecka2008} and \cite{Hale1988}. For gradient systems weaker conditions may suffice; see \cite{BabinVishik1992}. A non-vacuous instance is given in Example~\ref{ex:lower}; verifying (S) for the full range of nonlinearities in \eqref{eq:h_assumptions} remains open.
\end{remark}

\begin{definition}[Generalised hyperbolicity]\label{def:gen_hyp}
Let $E_*$ be an equilibrium of $S_0(t)$. When the nonlinearity in \eqref{eq:F} is not differentiable at $E_*$ (the coupling $[\cdot]^+$ may vanish on the switching surface $u=v$), the classical Jacobian must be replaced by the set $\Lambda(E_*)$ of generalised linearisations, obtained as limits of difference quotients of $F$ along sequences (equivalently, the Clarke generalised Jacobian of the finite-dimensional nonlinearities). We call $E_*$ \emph{generalised hyperbolic} if every $A\in\Lambda(E_*)$ generates a $C_0$-semigroup admitting an exponential dichotomy with a finite-dimensional unstable subspace. In particular, if every $A\in\Lambda(E_*)$ is exponentially stable, then $E_*$ is generalised hyperbolic with trivial unstable manifold. For $C^1$ nonlinearities this reduces to the usual notion of hyperbolicity.
\end{definition}

\begin{example}[A non-vacuous instance of condition (S)]\label{ex:lower}
Let $(0,\ell)=(0,\pi)$ and take
\[
\mu_1=\mu_2=5,\quad \gamma_1=\gamma_2=1,\quad k=1,\quad
f_1=f_2=0,\quad h_1(u)=u,\quad h_2(v)=v,
\]
with exponential kernels $g_1(s)=g_2(s)=2e^{-s}$. Then $\beta_i=2$, $\bar\mu_i=3$, and \eqref{eq:h_assumptions} holds with $H_i(s)=s^2/2$, $C_1=2$, $C_2=1$.

The limit system is
\[
u_{tt}+u_t+3u_{xxxx}+[u-v]^++u=0,\qquad
v_{tt}+v_t-3v_{xx}-[u-v]^++v=0.
\]
Testing the equilibrium equations with $u$ and $v$ gives
$3\|u_{xx}\|^2+3\|v_x\|^2+\|u\|^2+\|v\|^2+\| [u-v]^+\|^2=0$, so $(0,0)$ is the only equilibrium. Since $[\cdot]^+$ is not differentiable at $0$ and the equilibrium lies on the switching surface $u=v$, the classical linearisation does not exist. We instead use the set-valued (monotone) linearisation: the generalised gradient of $[\cdot]^+$ at $0$ is $[0,1]$, so for any measurable $\chi(x,t)\in[0,1]$ the linearised system is
\[
u_{tt}+u_t+3u_{xxxx}+\chi(u-v)+u=0,\qquad
v_{tt}+v_t-3v_{xx}-\chi(u-v)+v=0.
\]
For every such $\chi$ the energy
\[
E_\chi=\tfrac12\big(\|u_t\|^2+3\|u_{xx}\|^2+\|u\|^2+\|v_t\|^2+3\|v_x\|^2+\|v\|^2+\chi\|u-v\|^2\big)
\]
satisfies $E_\chi'\le-\|u_t\|^2-\|v_t\|^2$. Adding $\omega[(u,u_t)+(v,v_t)]$ with $0<\omega\le1/2$ gives, after cancelling the coupled terms $-\omega\chi\|u-v\|^2$,
\[
\frac{\dd}{\dd t}(E_\chi+\omega[(u,u_t)+(v,v_t)])
\le-c_\omega\big(E_\chi+\|u_t\|^2+\|v_t\|^2\big)
\]
with $c_\omega>0$ independent of $\chi$. Hence every generalised linearisation is exponentially stable, uniformly in $\chi\in[0,1]$. By Definition~\ref{def:gen_hyp}, $(0,0)$ is generalised hyperbolic with trivial unstable manifold. It remains to verify the persistence clause in (S). For every $\eps>0$, the state $(0,0)$ with zero memory components is an equilibrium of $S_\eps$, since $f_i=0$, $h_i(0)=0$, $[0]^+=0$ and a constant history has zero memory. Every equilibrium belongs to the global attractor, so $(0,0)\in\Acal_\eps$ for every $\eps$. Since $\Acal_0=\{(0,0)\}$, we have $J_0\Acal_0\subset J_\eps\Acal_\eps$ and therefore $\dist_{\Hcal_\eps}(\Acal_0,\Acal_\eps)=0$ for every $\eps$. The same energy computation applies to the linearisations of the perturbed semigroup $S_\eps$: adding $\|\sigma\|^2_{\Mg{1}}+\|\eta\|^2_{\Mg{2}}$ to $E_\chi$, the transport equations contribute the additional dissipation $-\frac{\lambda_1}{2}\|\sigma\|^2_{\Mg{1}}-\frac{\lambda_2}{2}\|\eta\|^2_{\Mg{2}}$ exactly as in \eqref{eq:proof_dissipation2}, while the $\chi$-terms cancel as before. Hence every generalised linearisation of $S_\eps$ is exponentially stable, uniformly in $\chi\in[0,1]$ and $\eps\in(0,1]$, so $e_\eps=(0,0)$ is generalised hyperbolic for $S_\eps$ as required by (S). Since the unstable manifolds are trivial, $J_\eps W^u_{\mathrm{loc}}(e_\eps)=\{(0,0)\}=J_0W^u_{\mathrm{loc}}(e_0)$ for every $\eps$, so the Hausdorff convergence in (S) is automatic once $e_\eps=e_0=(0,0)\in\Acal_\eps$, as shown above. Thus the example verifies the structural hypothesis, including persistence of the local unstable manifolds, without any appeal to a general perturbation theorem.
\end{example}

\section{Extensions of the Continuity Theory}\label{sec:extensions}

\subsection{Non-mass-preserving perturbations with a nondegenerate limit}\label{subsec:drift}

Mass preservation is essential only through the positivity of the effective elastic coefficients. If the kernel mass varies with $\eps$ but converges to a value strictly below $\mu_i$, the entire continuity theory remains valid, with the limiting coefficient shifted accordingly. Let
\begin{equation}\label{eq:mass_drift}
0<\underline\beta\le\beta_{i,\eps}\le\overline\beta<\mu_i,\qquad
\beta_{i,\eps}\to\beta_{i,0}\in(0,\mu_i),
\end{equation}
and replace \eqref{eq:rescaled} by
\begin{equation}\label{eq:kernel_drift}
\tilde g_{i,\eps}(s)\coloneqq\frac{\beta_{i,\eps}}{\beta_i\eps}\,g_i\!\Big(\frac{s}{\eps}\Big).
\end{equation}
Then $\int_0^{+\infty}\tilde g_{i,\eps}=\beta_{i,\eps}$, $\tilde g_{i,\eps}'+\lambda_i\tilde g_{i,\eps}\le0$, and
\[
\int_0^{+\infty}\!s\,\tilde g_{i,\eps}(s)\dd s
=\frac{\beta_{i,\eps}}{\beta_i}\,\eps\int_0^{+\infty}\!\zeta g_i(\zeta)\dd\zeta
\le\frac{\overline\beta}{\beta_i}\frac{g_i(0)}{\lambda_i^2}\,\eps,
\]
while for $R\ge0$
\[
\int_R^{+\infty}\!\tilde g_{i,\eps}(s)\dd s
\le\frac{\overline\beta}{\beta_i}\frac{g_i(0)}{\lambda_i}\,e^{-\lambda_i R/\eps}
\le\frac{\overline\beta}{\beta_i}\frac{g_i(0)}{\lambda_i}\,e^{-\lambda_i R},
\qquad \eps\le1.
\]
The effective coefficients are
\begin{equation}\label{eq:drift_coeff}
\bar\mu_{i,\eps}=\mu_i-\beta_{i,\eps}\in[\mu_i-\overline\beta,\mu_i-\underline\beta]\subset(0,+\infty),
\qquad \bar\mu_{i,\eps}\to\bar\mu_{i,0}=\mu_i-\beta_{i,0}>0 .
\end{equation}

\begin{theorem}[Continuity under a nondegenerate mass drift]\label{thm:drift}
Assume \eqref{eq:h_assumptions}, \eqref{eq:kernel}, and \eqref{eq:mass_drift}, and let the kernels be given by \eqref{eq:kernel_drift}. Then the conclusions of Theorem~\ref{thm:upper} hold with the limit system \eqref{eq:limit_intro} in which $\bar\mu_i$ is replaced by $\bar\mu_{i,0}$. If the limit attractor satisfies (S), then the lower semicontinuity \eqref{eq:lower} holds as well.
\end{theorem}

\begin{proof}
All estimates in Sections~\ref{sec:functional}--\ref{sec:lower} involve the elastic coefficients only through the constants $\bar\mu_i$; replacing them by $\bar\mu_{i,\eps}$ from \eqref{eq:drift_coeff} leaves every constant uniform in $\eps$, because $\bar\mu_{i,\eps}$ stays in a fixed compact subset of $(0,+\infty)$. In particular, the energy estimates of Section~\ref{sec:dissipativity}, the resolvent estimates of Theorem~\ref{thm:unif_exp}, and the compactness argument of Proposition~\ref{prop:quasi} apply mutatis mutandis, with constants uniform because $\bar\mu_{i,\eps}$ stays in a compact subset of $(0,+\infty)$.

For the resolvent limit, the computation now gives the coefficient
\[
a_{\eps,\mu}=\int_0^{+\infty}\tilde g_{1,\eps}(s)\frac{1-e^{-\mu s}}{\mu}\dd s
\le\frac{\overline\beta}{\beta_1}\,\eps\int_0^{+\infty}\!\zeta g_1(\zeta)\dd\zeta\to0,
\]
and similarly for the second component. The limiting elastic operator is therefore $\Lcal_{e,0}$ with coefficients $\bar\mu_{i,0}$; since $\bar\mu_{i,\eps}\to\bar\mu_{i,0}$ and the elliptic estimates are uniform, the finite-mode Duhamel argument of Proposition~\ref{prop:conv} goes through without change. In the memory part of that proof one uses the transport identity $\sigma=u(t)-u(t-s)$ and the tail bound above; the displacement remains bounded in $H_a^2$ uniformly because of the positive lower bound in \eqref{eq:drift_coeff}, so the finite-mode compactness argument of Proposition~\ref{prop:conv} applies without change. Hence finite-time convergence and upper semicontinuity follow as before.

Finally, lower semicontinuity follows from Proposition~\ref{prop:lower} because $\Acal_0$ is compact: the finite-mode argument in Proposition~\ref{prop:conv} uses only the uniform lower bound in \eqref{eq:drift_coeff} and the convergence $\beta_{i,\eps}\to\beta_{i,0}$. No kernel lower bound or parameter-H\"older estimate is needed.
\end{proof}

\begin{remark}\label{rem:drift_degenerate}
If $\beta_{i,0}\ge\mu_i$, then $\bar\mu_{i,0}\le0$ and the limiting elastic operator loses coercivity; the estimates above break down and the problem becomes genuinely singular. We return to this degenerate case in the last section.
\end{remark}

\subsection{The limit of infinitely fast exponential decay}\label{subsec:fast_decay}

The limit $\lambda_i\to+\infty$ for exponential kernels is not a different singular limit: it is exactly the rescaling \eqref{eq:rescaled} after a change of parameter.

\begin{proposition}\label{prop:fast_decay}
Let $g_i(s)=\beta_ie^{-s}$ and, for $\lambda>0$, set $g_i^{(\lambda)}(s)\coloneqq\beta_i\lambda e^{-\lambda s}$. Then $g_i^{(\lambda)}=g_{i,\eps}$ with $\eps=1/\lambda$. Consequently, as $\lambda\to+\infty$, all continuity conclusions of Theorem~\ref{thm:upper}, and of Proposition~\ref{prop:lower} when (S) holds, apply directly with $\bar\mu_i=\mu_i-\beta_i$; for a single damped mode the slow eigenvalues converge to the roots of $\lambda^2+\gamma_i\lambda+\bar\mu_i=0$.
\end{proposition}

\begin{proof}
Indeed, $g_{i,\eps}(s)=\eps^{-1}\beta_i e^{-s/\eps}=\beta_i\lambda e^{-\lambda s}$ with $\lambda=\eps^{-1}$, while $\int_0^{+\infty}g_i^{(\lambda)}=\beta_i$ and $(g_i^{(\lambda)})'+\lambda g_i^{(\lambda)}=0$, so \eqref{eq:kernel} holds uniformly. The characteristic polynomial of the corresponding damped mode has coefficients that converge to those of $\lambda^2+\gamma_i\lambda+\bar\mu_i=0$ as $\lambda\to+\infty$, so its slow roots converge to the roots of that limiting quadratic; this is confirmed numerically in Section~\ref{sec:numerics} below.
\end{proof}

\begin{remark}\label{rem:fast_general}
For a general family of kernels with $\lambda_i\to+\infty$ the identification with \eqref{eq:rescaled} requires that the family be a rescaling of a fixed profile; under that structural condition the same conclusion holds. Without it, additional uniform assumptions are needed and the question is not covered by the present analysis. The opposite limit $\lambda_i\to0^+$ spreads the mass over an unbounded time interval, the tail estimates of Section~\ref{sec:functional} fail, and the question remains open.
\end{remark}

\subsection{Towards quantitative rates}\label{subsec:rates}

Theorems~\ref{thm:upper} and Proposition~\ref{prop:lower} are qualitative because finite-time convergence is qualitative. The next statement is conditional: it isolates the two extra ingredients that would yield a rate.

\begin{proposition}[Conditional quantitative lower semicontinuity]\label{prop:quant}
Assume the hypotheses of Theorem~\ref{thm:upper}. Let $E_\eps:\Hcal_e\to\Hcal_\eps$ assign zero memory components. Assume, in addition, that
\begin{enumerate}
\item[(Q1)] there exist $M\ge1$ and $\omega>0$, independent of $\eps$, such that
$\dist_{\Hcal_\eps}(S_\eps(t)\Bcal_{0,\eps},\Acal_\eps)\le Me^{-\omega t}$ for all $t\ge0$;
\item[(Q2)] for some $\rho:(0,1]\times(0,+\infty)\to(0,+\infty)$,
\[
\sup_{U\in\Bcal_0}\ \sup_{t\in[0,T]}
\|S_\eps(t)E_\eps U-S_0(t)U\|_{\Hcal_\eps}\le\rho(\eps,T).
\]
\end{enumerate}
Then
\begin{equation}\label{eq:attractor_rate}
\dist_{\Hcal_\eps}(\Acal_0,\Acal_\eps)
\le\inf_{T>0}\big\{\rho(\eps,T)+Me^{-\omega T}\big\}.
\end{equation}
If $\rho(\eps,T)\le C\eps^\alpha e^{KT}$ with $\alpha,K>0$, then
$\dist_{\Hcal_\eps}(\Acal_0,\Acal_\eps)\le C'\eps^{\alpha\omega/(\omega+K)}$.
\end{proposition}

\begin{proof}
Let $U\in\Acal_0$ and fix $T>0$. Since $\Acal_0$ is invariant, $S_0(T)\Acal_0=\Acal_0$, so $U=S_0(T)V$ for some $V\in\Acal_0\subset\Bcal_0$. By (Q2),
$\|S_\eps(T)E_\eps V-U\|_{\Hcal_\eps}\le\rho(\eps,T)$. By (Q1), $\dist(S_\eps(T)E_\eps V,\Acal_\eps)\le Me^{-\omega T}$. The triangle inequality and the supremum over $U$ give \eqref{eq:attractor_rate}; the last assertion follows by balancing the two terms.
\end{proof}

\begin{remark}\label{rem:rate_obstruction}
Condition (Q1) is not a consequence of quasi-stability alone; it is an additional exponential-attraction hypothesis. Condition (Q2) is the genuine finite-time convergence rate. On compact sets of initial data the elastic part of that rate is $O(\eps)$, as $a_{\eps,\mu}=O(\eps)$ shows, while the memory part is controlled only by the modulus of continuity of the displacement in $H_a^2$. Compactness gives convergence of this modulus but no rate, because the energy estimates bound $u_t$ in $L^2$ and not in $H_a^2$. If one additionally assumes $\omega_\Phi(\delta)\le L\delta$ uniformly, the finite-mode argument gives $\rho(\eps,T)=C_T\eps|\log\eps|$ and hence a rate from Proposition~\ref{prop:quant}. Whether this regularity holds under the standing assumptions is open; the same argument gives no rate for $\dist(\Acal_\eps,\Acal_0)$, since trajectories on $\Acal_\eps$ carry non-vanishing memory.
\end{remark}

\section{Numerical Experiments}\label{sec:numerics}

The continuity results established in Sections~\ref{sec:functional}--\ref{sec:extensions} are qualitative: they assert convergence without providing a rate. We now report four numerical experiments that check the theory and probe the results that lie beyond it. Specifically, they test (i) that the dynamics selects the limiting coefficient $\bar\mu_i=\mu_i-\beta_i$, rather than the original $\mu_i$; (ii) that the finite-time convergence of Proposition~\ref{prop:conv} is already visible on moderate time intervals; (iii) that the identification of the fast-decay limit in Proposition~\ref{prop:fast_decay} is reflected in the slow poles; (iv) that the nondegenerate mass drift of Theorem~\ref{thm:drift} shifts the limit stiffness to $\mu_i-\beta_{i,0}$; and (v) that the empirical rate degrades exactly when the displacement loses its $H_a^2$-time regularity, as diagnosed in Remark~\ref{rem:rate_obstruction}. The outcome also motivates the questions formulated in the last section: the observed convergence rates are not proved by the theory, and the experiments do not reach the attractor scale. All experiments use the exponential kernels
\begin{equation}\label{eq:num_kernel}
g_i(s)=\beta_i e^{-s},\qquad g_{i,\eps}(s)=\frac{\beta_i}{\eps}e^{-s/\eps},
\end{equation}
which satisfy \eqref{eq:kernel} with $\lambda_i=1$. For such kernels the infinite history can be eliminated in favour of two auxiliary variables: setting
$z_1(x,t)=\int_0^{+\infty}g_{1,\eps}(s)u(x,t-s)\dd s$ and
$z_2(x,t)=\int_0^{+\infty}g_{2,\eps}(s)v(x,t-s)\dd s$, one has
\begin{equation}\label{eq:num_aux}
\eps\,\partial_t z_1+z_1=\beta_1u,\qquad
\eps\,\partial_t z_2+z_2=\beta_2v,
\end{equation}
and the convolution terms in \eqref{eq:Aeps} become $(\beta_1u-z_1)_{xxxx}$ and $(\beta_2v-z_2)_{xx}$, respectively. Thus the numerical scheme integrates a finite-dimensional system in history, with $z_i$ initialised from a zero past history.

\subsection{Test 1: a single mode}

\begin{figure}[htbp]
\centering
\includegraphics[width=0.97\textwidth]{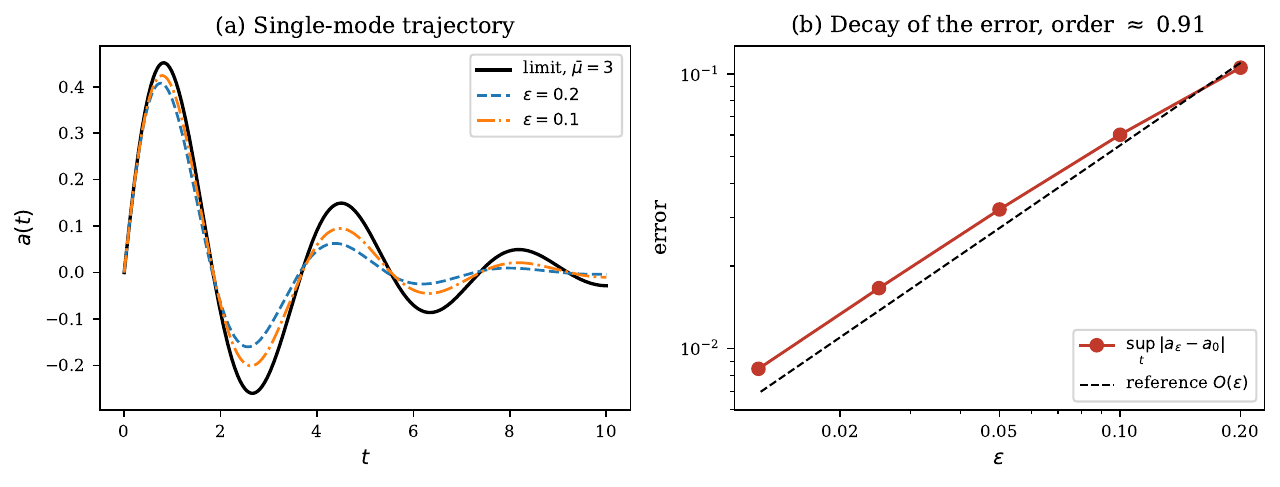}
\caption{Test~1: (a) mode amplitude for two values of $\eps$ and for the memoryless limit; (b) uniform error against $\eps$ (log--log), with an $O(\eps)$ reference line.}
\label{fig:mode}
\end{figure}

Project the deck equation onto the first sine mode, $u(x,t)=a(t)\sin x$, and drop the coupling, the lower-order nonlinearity and the forcing. With $\mu_1=5$, $\beta_1=2$ (hence $\bar\mu_1=3$) and $\gamma_1=0.6$, the amplitude and its auxiliary variable solve
\begin{equation}\label{eq:num_mode}
\begin{cases}
\ddot a+\gamma_1\dot a+\bar\mu_1 a+(\beta_1a-z)=0,\\
\eps\dot z+z=\beta_1a,
\end{cases}
\qquad a(0)=0,\quad \dot a(0)=1,\quad z(0)=0.
\end{equation}
Substituting $a=e^{\lambda t}$ into \eqref{eq:num_mode} gives, after eliminating $z$ by means of $z=\beta_1a/(1+\eps\lambda)$, the characteristic equation
\begin{equation}\label{eq:num_char}
(\lambda^2+\gamma_1\lambda+\bar\mu_1)(1+\eps\lambda)
+\beta_1\eps\lambda=0,
\end{equation}
that is, $\eps\lambda^3+(1+\gamma_1\eps)\lambda^2+[\gamma_1+(\bar\mu_1+\beta_1)\eps]\lambda+\bar\mu_1=0$.
Its two slow roots converge to $-0.3\pm1.7059\,i$, the roots of
$\lambda^2+0.6\lambda+3=0$; the third root is real, approximately $-1/\eps$, and represents the fast relaxation mode. Since \eqref{eq:num_kernel} is the exponential family of Proposition~\ref{prop:fast_decay} with $\lambda=1/\eps$, the rows $\eps=0.2,\dots,0.0125$ below correspond to $\lambda=5,\dots,80$ and confirm that proposition. Expanding the small roots of \eqref{eq:num_char} gives $\lambda_\eps=\lambda_0-\eps\beta_1\lambda_0/(2\lambda_0+\gamma_1)+O(\eps^2)$ for each root $\lambda_0$ of $\lambda^2+\gamma_1\lambda+\bar\mu_1=0$, which explains the pole convergence analytically; the numerical roots and trajectory errors are listed in Table~\ref{tab:num_poles}. Had the limit coefficient been incorrectly taken as $\mu_1=5$, the limiting roots would be $-0.3\pm2.2159\,i$, which are clearly not approached. Figure~\ref{fig:mode}(a) shows the trajectories for $\eps=0.2,0.1$ against the limit solution; Figure~\ref{fig:mode}(b) plots the error as a function of $\eps$, with an empirical slope $0.91$ on $[0.0125,0.2]$.
\begin{table}[htbp]
\centering
\caption{Test~1: slow roots of \eqref{eq:num_char} and uniform error on $[0,15]$; the limit roots are $-0.3\pm1.7059\,i$.}
\label{tab:num_poles}
\begin{tabular}{c|c|c}
\hline
$\varepsilon$ & slow roots & $\sup_{t\in[0,15]}|a_\eps(t)-a_0(t)|$\\
\hline
$0.2$    & $-0.5162\pm1.7371\,i$ & $1.06\times10^{-1}$\\
$0.1$    & $-0.4053\pm1.7030\,i$ & $6.01\times10^{-2}$\\
$0.05$   & $-0.3514\pm1.7006\,i$ & $3.21\times10^{-2}$\\
$0.025$  & $-0.3254\pm1.7023\,i$ & $1.66\times10^{-2}$\\
$0.0125$ & $-0.3126\pm1.7039\,i$ & $8.44\times10^{-3}$\\
\hline
\end{tabular}
\end{table}

\subsection{Test 2: the coupled system}

\begin{figure}[htbp]
\centering
\includegraphics[width=0.97\textwidth]{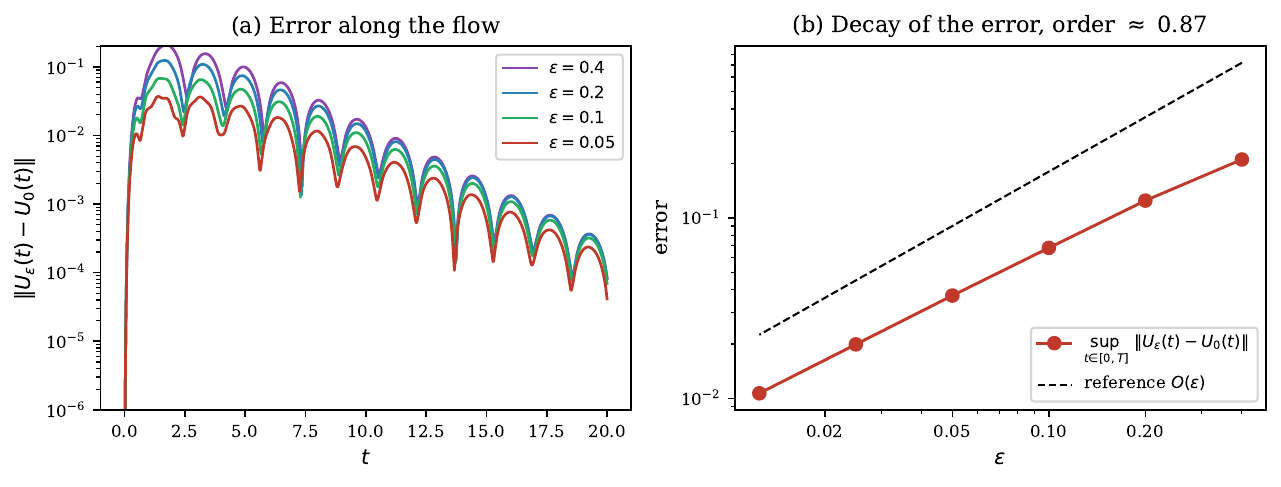}
\caption{Test~2: (a) $E_\eps(t)$ along the flow for four values of $\eps$; (b) $\sup_{t\in[0,20]}E_\eps(t)$ against $\eps$ (log--log) with an $O(\eps)$ reference line.}
\label{fig:bridge}
\end{figure}

We discretise \eqref{eq:Mukiawa} on $(0,\pi)$ by the Fourier--Galerkin expansion
$u=\sum_{n=1}^{6}a_n(t)\sin nx$, $v=\sum_{n=1}^{6}b_n(t)\sin nx$. The non-smooth coupling $k[u-v]^+$ is evaluated on a $128$-point grid and projected back onto the sine basis; the resulting stiff system, including the two auxiliary variables per mode in \eqref{eq:num_aux}, is integrated by an implicit Radau method with relative tolerance $10^{-8}$. We take
\[
\mu_1=\mu_2=5,\quad \beta_1=\beta_2=2,\quad \bar\mu_1=\bar\mu_2=3,\quad
\gamma_1=\gamma_2=0.8,\quad k=1,\quad h_i(s)=s,
\]
\[
f_1(x)=3\sin x+0.5\sin2x,\qquad f_2(x)=1.5\sin x,
\]
with all initial data and the past history equal to zero. The limit solution is computed from the same code with the convolutions removed and the coefficients $\bar\mu_i$ in place of $\mu_i$. For each value of $\eps$ we measure the discrete energy distance
\begin{equation}\label{eq:num_error}
E_\eps(t)=\Big(\sum_{n=1}^{6}n^4|a_n^\eps(t)-a_n^0(t)|^2
+\sum_{n=1}^{6}n^2|b_n^\eps(t)-b_n^0(t)|^2\Big)^{1/2}.
\end{equation}
For $\eps=0.4,0.2,0.1,0.05,0.025,0.0125$, the values of $\sup_{t\in[0,20]}E_\eps(t)$ are
$2.10\times10^{-1}$, $1.24\times10^{-1}$, $6.81\times10^{-2}$,
$3.71\times10^{-2}$, $1.99\times10^{-2}$ and $1.07\times10^{-2}$, respectively, corresponding to an empirical slope $0.87$ (the two-decimal value; the refinement below resolves it as $0.865$--$0.866$); the time profiles are shown in Figure~\ref{fig:bridge}. To check that the slope is not a discretisation artefact, we repeated the experiment with $N=4,6,8$ Fourier modes and grids of $64,128,256$ points for the projection of $[u-v]^+$; the fitted slopes all round to $0.866$ at three decimal places, as reported in Table~\ref{tab:grid}. Halving the Radau relative tolerance from $10^{-8}$ to $10^{-9}$ changed the six recorded suprema by at most $4.6\times10^{-5}$ in relative terms, so the displayed digits are tolerance-stable.
\begin{table}[htbp]
\centering
\caption{Test~2: fitted slope for $\sup_{[0,20]}E_\eps$ under mesh and tolerance refinement.}
\label{tab:grid}
\begin{tabular}{c|c|c}
\hline
modes $N$ & grid points & fitted slope\\
\hline
$4$  & $128$ & $0.866$\\
$6$  & $64$  & $0.866$\\
$6$  & $128$ & $0.866$\\
$6$  & $256$ & $0.866$\\
$8$  & $128$ & $0.866$\\
\hline
\end{tabular}
\end{table}
\subsection{Test 3: a nondegenerate mass drift}

\begin{figure}[htbp]
\centering
\includegraphics[width=0.97\textwidth]{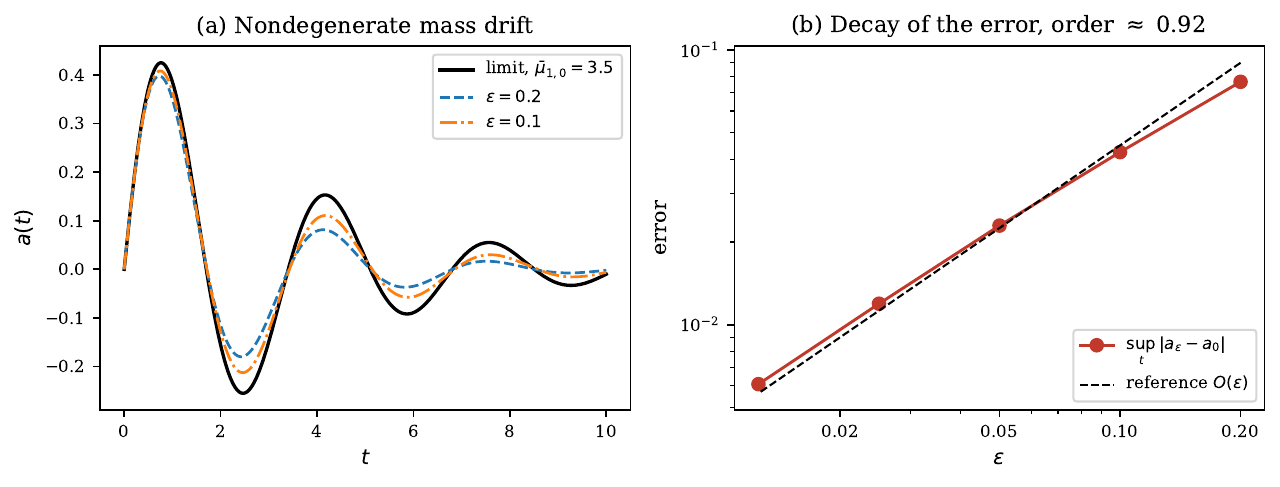}
\caption{Test~3: (a) mode amplitude for two values of $\eps$ and for the memoryless limit with stiffness $\mu_1-\beta_{1,0}=3.5$; (b) uniform error against $\eps$ (log--log), with an $O(\eps)$ reference line.}
\label{fig:drift}
\end{figure}

To check Theorem~\ref{thm:drift} we repeat Test~1 with a kernel whose mass is not preserved but converges. Take $\mu_1=5$, $\gamma_1=0.6$, and $\beta_{1,\eps}=1.5+0.5\eps$, so that $\beta_{1,\eps}\to\beta_{1,0}=1.5$ from above and stays in $(0,\mu_1)$; the rescaled kernel is $\tilde g_{1,\eps}(s)=\beta_{1,\eps}\eps^{-1}e^{-s/\eps}$. According to Theorem~\ref{thm:drift} the limit stiffness should be $\bar\mu_{1,0}=\mu_1-\beta_{1,0}=3.5$, so the slow poles must converge to the roots of $\lambda^2+0.6\lambda+3.5=0$, namely $-0.3\pm1.8466\,i$, and the errors must tend to zero. Repeating the computation of Test~1 with $\beta_1$ replaced by $\beta_{1,\eps}$ and the limit coefficient by $3.5$ gives
\begin{table}[htbp]
\centering
\caption{Test~3: slow roots and uniform error under nondegenerate mass drift; the limit roots are $-0.3\pm1.8466\,i$.}
\label{tab:drift}
\begin{tabular}{c|c|c}
\hline
$\varepsilon$ & slow roots & $\sup_{t\in[0,15]}|a_\eps-a_0|$\\
\hline
$0.2$    & $-0.4668\pm1.8507\,i$ & $7.67\times10^{-2}$\\
$0.1$    & $-0.3808\pm1.8335\,i$ & $4.25\times10^{-2}$\\
$0.05$   & $-0.3391\pm1.8367\,i$ & $2.30\times10^{-2}$\\
$0.025$  & $-0.3192\pm1.8409\,i$ & $1.19\times10^{-2}$\\
$0.0125$ & $-0.3095\pm1.8436\,i$ & $6.07\times10^{-3}$\\
\hline
\end{tabular}
\end{table}
with an empirical slope $0.92$. The poles approach $-0.3\pm1.8466\,i$, not the value $-0.3\pm2.2159\,i$ that would correspond to the original stiffness $\mu_1=5$. Thus the limiting coefficient is indeed $\mu_1-\beta_{1,0}$, as asserted by Theorem~\ref{thm:drift}. The trajectories for $\eps=0.2,0.1$ against the limiting solution and the decay of the uniform error are displayed in Figure~\ref{fig:drift}.

\subsection{Test 4: the rate bottleneck}

\begin{figure}[htbp]
\centering
\includegraphics[width=0.97\textwidth]{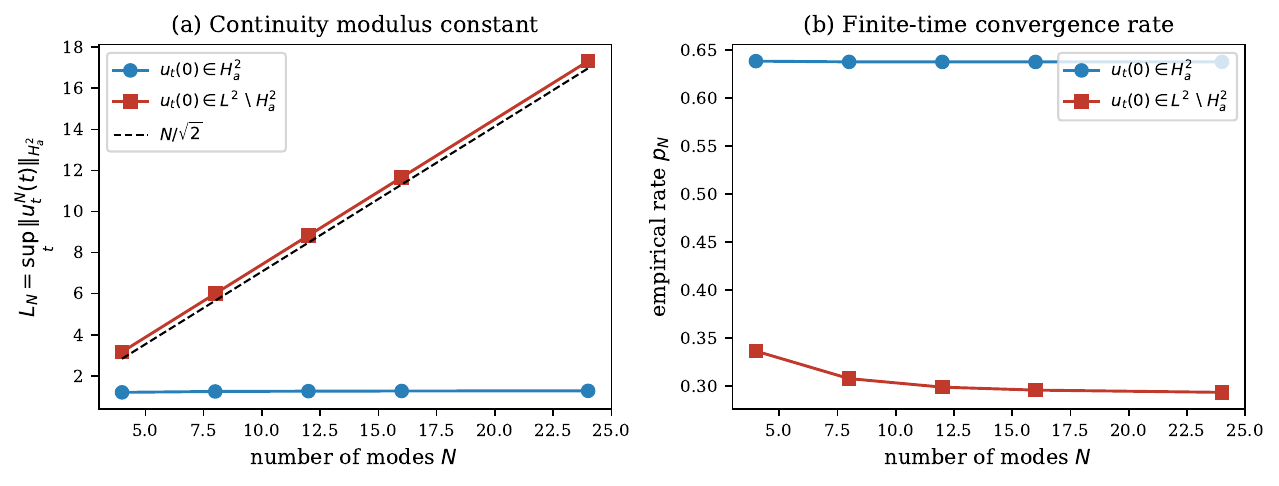}
\caption{Test~4: (a) the continuity-modulus constant $L_N$ against the number of modes; the dashed line is $N/\sqrt2$; (b) the empirical finite-time rate $p_N$ against $N$.}
\label{fig:rate}
\end{figure}

Remark~\ref{rem:rate_obstruction} traces the absence of a rate to the continuity modulus of the displacement in $H_a^2$: the elastic part of the finite-time error is $O(\eps)$ on compact time intervals, whereas the memory part is controlled only by $\omega_\Phi(\delta)$, since the energy estimates bound $u_t$ in $L^2$ and not in $H_a^2$. We probe this directly on the Galerkin truncation of \eqref{eq:num_mode} without forcing, with $\mu_1=5$, $\beta_1=2$, $\bar\mu_1=3$, $\gamma_1=0.6$, and with
\begin{equation}\label{eq:num_lowreg}
a_n(0)=0,\qquad \dot a_n(0)=n^{-s},\qquad z_n(0)=0,\qquad n=1,\dots,N .
\end{equation}
For $s=3$ the initial velocity belongs to $H_a^2(0,\pi)$; for $s=3/2$ it belongs to $L^2(0,\pi)\setminus H_a^2(0,\pi)$, which is the regularity actually supplied by the energy estimate. Let
\begin{equation}\label{eq:num_modulus}
L_N=\sup_{t\in[0,5]}\Big(\sum_{n=1}^N n^4|\dot a_n(t)|^2\Big)^{1/2},
\end{equation}
and let $p_N$ be the least-squares slope of $\log\sup_{t\in[0,5]}\|u_\eps(t)-u_0(t)\|_{H_a^2}$ against $\log\eps$ for $\eps=0.2,0.1,0.05,0.025$. The results are reported in Table~\ref{tab:num_rate}.
\begin{table}[htbp]
\centering
\caption{Test~4: the continuity-modulus constant $L_N$ in \eqref{eq:num_modulus} and the empirical finite-time rate $p_N$ for the two initial velocities in \eqref{eq:num_lowreg}.}
\label{tab:num_rate}
\begin{tabular}{c|cc|cc}
\hline
 & \multicolumn{2}{c|}{$s=3$: $\dot u(0)\in H_a^2$}
 & \multicolumn{2}{c}{$s=3/2$: $\dot u(0)\notin H_a^2$}\\
$N$ & $L_N$ & $p_N$ & $L_N$ & $p_N$\\
\hline
$4$  & $1.193$ & $0.639$ & $3.162$  & $0.337$\\
$8$  & $1.236$ & $0.638$ & $6.000$  & $0.308$\\
$12$ & $1.251$ & $0.638$ & $8.832$  & $0.299$\\
$16$ & $1.259$ & $0.638$ & $11.662$ & $0.296$\\
$24$ & $1.267$ & $0.638$ & $17.321$ & $0.294$\\
\hline
\end{tabular}
\end{table}
For $s=3$ the constant $L_N$ stabilises near $1.27$ and $p_N$ remains at $0.638$: the Lipschitz regularity invoked in Remark~\ref{rem:rate_obstruction} holds uniformly in $N$, and the finite-time error has a stable empirical rate. For $s=3/2$, $L_N$ grows like $N/\sqrt2$ (at $N=24$, $L_N=17.321$ against $24/\sqrt2=16.971$) while $p_N$ drops from $0.337$ to $0.294$: there is no $H_a^2$-Lipschitz bound uniform in $N$, and the finite-time error converges more slowly. Both trends are displayed in Figure~\ref{fig:rate}. The four-point log--log fits are stable: for $s=3$ the $R^2$ value is $0.998$ for every $N$ and the largest residual in logarithmic error is $0.033$; for $s=3/2$ the $R^2$ values increase from $0.997$ to $1.000$ and the largest residual falls from $0.018$ to $0.006$. Thus the difference between the two families of slopes is not an artefact of the fit: it reflects the growth of $L_N$ displayed in the same table.
This experiment neither proves nor disproves a rate for the attractors; it shows that the regularity gap identified in Remark~\ref{rem:rate_obstruction} is numerically detectable and that the slopes observed in Tests~1--2 cannot be expected for initial data at the energy regularity alone. This dependence of the finite-time rate on the continuity modulus is what motivates question~(a) in the last section. Table~\ref{tab:num_map} records which analytical statement each test addresses.
\begin{table}[htbp]
\centering
\caption{Correspondence between the numerical tests and the analytical results they check.}
\label{tab:num_map}
\begin{tabular}{p{0.8cm}|p{5.3cm}|p{7.4cm}}
\hline
Test & Analytical statement & Main outcome\\
\hline
$1$ & limiting coefficient $\bar\mu_i=\mu_i-\beta_i$; Proposition~\ref{prop:fast_decay}
    & poles $\to-0.3\pm1.7059\,i$; slope $0.91$\\
$2$ & finite-time convergence, Proposition~\ref{prop:conv}
    & slope $0.87$ for $\sup_{[0,20]}E_\eps$\\
$3$ & nondegenerate mass drift, Theorem~\ref{thm:drift}
    & poles $\to-0.3\pm1.8466\,i$, stiffness $\mu_1-\beta_{1,0}=3.5$; slope $0.92$\\
$4$ & rate bottleneck, Remark~\ref{rem:rate_obstruction}; Proposition~\ref{prop:quant}
    & $L_N$ bounded versus $L_N\sim N$; $p_N\approx0.64$ versus $0.29$\\
\hline
\end{tabular}
\end{table}

In all four tests the error is largest during the transient and subsequently decays: in the coupled example the memory contribution vanishes for time-independent states, so all trajectories approach the same equilibrium and the difference between $S_\eps(t)$ and $S_0(t)$ is a finite-time phenomenon. This is consistent with the mechanism of Proposition~\ref{prop:conv}, although it does not by itself probe the full attractor distance, since the invariant sets here contain no non-stationary component. The four experiments support the analytical conclusions, but the empirical slopes are not proofs and no convergence rate for the attractors is claimed; the question of such a rate is taken up again in the last section.

\section{Summary and Open Problems}\label{sec:conclusion}

We have developed a continuity theory for the global attractor of a dissipative evolution equation with memory under a mass-preserving singular limit of the memory kernel, and we have verified its ingredients for a coupled suspension bridge system with past history. The key elements are:
\begin{itemize}
\item[(1)] The positivity condition $\inf_\eps(\mu_i-\beta_{i,\eps})>0$, which keeps the effective elastic coefficients uniformly coercive; the mass-preserving rescaling is the cleanest way to guarantee it.
\item[(2)] The uniform exponential stability \eqref{eq:exp_stab} of the linear semigroup, established in \cite{Zhou2026} for fixed $\eps$ and extended to a uniform-in-$\eps$ version in Theorem~\ref{thm:unif_exp}.
\item[(3)] The uniform quasi-stability (Proposition~\ref{prop:quasi}), obtained from the uniform exponential stability and the sequential compactness of the integral remainder $K_{\eps,T}(\Bcal_0)$ in the common reference space.
\item[(4)] The finite-time convergence \eqref{eq:semigroup_conv} on compact sets of initial data, which together with (1)--(3) yields the upper semicontinuity.
\end{itemize}
Section~\ref{sec:extensions} shows that mass preservation can be relaxed to a nondegenerate mass drift, that the fast exponential decay limit is the same singular limit, and that a quantitative lower semicontinuity would follow from a finite-time rate together with an exponential-attraction hypothesis; the obstruction is identified in Remark~\ref{rem:rate_obstruction}. The numerical experiments of Section~\ref{sec:numerics} then support the ingredients that are most delicate in this chain (Table~\ref{tab:num_map}): the poles and trajectories of a damped mode converge to those of the system with stiffness $\bar\mu_i$, the energy-norm distance between perturbed and limit trajectories tends to zero on bounded time intervals, a nondegenerate mass drift shifts the limiting stiffness to $\mu_i-\beta_{i,0}$, and the finite-time rate degrades when the displacement lacks the $H_a^2$-time regularity required in Remark~\ref{rem:rate_obstruction}.

The results above settle several questions about the singular limit; the following remain open. Each is stated relative to what has just been proved:
\begin{itemize}
\item[(a)] \emph{Rate of convergence.} Can one obtain a quantitative rate for upper and lower semicontinuity under the standing assumptions alone? Proposition~\ref{prop:quant} reduces this to a rate for finite-time convergence and, under the additional Lipschitz regularity discussed in Remark~\ref{rem:rate_obstruction}, gives a rate for $\dist_{\Hcal_\eps}(\Acal_0,\Acal_\eps)$. Test~4 in Section~\ref{sec:numerics} supports this diagnosis numerically: the empirical finite-time rate falls from $0.638$ to $0.294$ when $\dot u(0)$ leaves $H_a^2$ and the truncated continuity-modulus constant grows like $N/\sqrt2$; it does not, however, supply a proof. A rate for $\dist_{\Hcal_\eps}(\Acal_\eps,\Acal_0)$, and the removal of the extra regularity assumption, remain open.
\item[(b)] \emph{Degenerate mass drift.} If the kernel mass varies with $\eps$, Theorem~\ref{thm:drift} settles the nondegenerate case $\beta_{i,0}<\mu_i$: continuity holds with limiting stiffness $\bar\mu_{i,0}=\mu_i-\beta_{i,0}$. What remains open is the genuinely singular case $\beta_{i,0}\ge\mu_i$, where $\bar\mu_{i,0}\le0$ and the elastic operator loses coercivity.
\item[(c)] \emph{Other singular limits.} The fast-decay limit $\lambda_i\to+\infty$ is settled for mass-preserving exponential families by Proposition~\ref{prop:fast_decay}, since it is exactly the rescaling \eqref{eq:rescaled} with $\eps=1/\lambda$; the general non-rescaling case and the slow-decay limit $\lambda_i\to0^+$ (where the mass escapes to infinity and the tail estimates fail) remain open.
\item[(d)] \emph{Lower semicontinuity without structural assumptions.} Is the lower semicontinuity true for the full range of nonlinearities $h_i$ satisfying \eqref{eq:h_assumptions}, without condition (S)? The set of equilibria is independent of $\eps$ (the memory contribution vanishes at equilibria), which removes the usual bifurcation mechanism; nevertheless, when an equilibrium of $\Acal_0$ is non-hyperbolic, continuity of its unstable manifold under perturbation is not guaranteed, and the standard abstract theorem does not apply. Example~\ref{ex:lower} shows that (S) is non-vacuous; removing it remains open.
\end{itemize}

\section*{Declarations}

\noindent\textbf{Conflict of interest:} The author has no relevant financial or non-financial interests to disclose.

\noindent\textbf{Data Availability Statement:} The numerical results reported in Section~\ref{sec:numerics} were produced by the author's code and are available from the author on reasonable request; no external data were used.


\begin{thebibliography}{99}

\bibitem{AhmedHarbi1998}
N.~U. Ahmed and H.~Harbi,
\emph{Mathematical analysis of dynamic models of suspension bridges},
SIAM J. Appl. Math. \textbf{58} (1998), no.~3, 853--874.

\bibitem{An2002}
Y.~An,
\emph{Nonlinear perturbations of a coupled system of steady state suspension bridge equations},
Nonlinear Anal. \textbf{51} (2002), 1285--1292.

\bibitem{BabinVishik1992}
A.~V. Babin and M.~I. Vishik,
\emph{Attractors of Evolution Equations},
North-Holland, Amsterdam, 1992.

\bibitem{BochicchioGiorgiVuk2010}
I.~Bochicchio, C.~Giorgi and E.~Vuk,
\emph{Long-term damped dynamics of the extensible suspension bridge},
Int. J. Differ. Equ. \textbf{2010} (2010), Article ID 383420.

\bibitem{ChueshovLasiecka2008}
I.~Chueshov and I.~Lasiecka,
\emph{Long-time behavior of second order evolution equations with nonlinear damping},
Mem. Amer. Math. Soc. \textbf{195} (2008), no.~912.

\bibitem{Gearhart1978}
L.~Gearhart,
\emph{Spectral theory for contraction semigroups on Hilbert space},
Trans. Amer. Math. Soc. \textbf{236} (1978), 385--394.

\bibitem{Pruss1984}
J.~Pr\"uss,
\emph{On the spectrum of $C_0$-semigroups},
Trans. Amer. Math. Soc. \textbf{284} (1984), 847--857.

\bibitem{ContiPataSquassina2006}
M.~Conti, V.~Pata and M.~Squassina,
\emph{Singular limit of differential systems with memory},
Indiana Univ. Math. J. \textbf{55} (2006), 169--215.

\bibitem{Dafermos1970}
C.~M. Dafermos,
\emph{An abstract Volterra equation with applications to linear viscoelasticity},
J. Differential Equations \textbf{7} (1970), 554--569.

\bibitem{DellOroGiorgiPata2015}
F.~Dell'Oro, C.~Giorgi and V.~Pata,
\emph{Asymptotic behavior of coupled linear systems modeling suspension bridges},
Z. Angew. Math. Phys. \textbf{66} (2015), 1095--1108.

\bibitem{EnyiMukiawa2020}
C.~D. Enyi and S.~E. Mukiawa,
\emph{Dynamics of a thermoelastic-laminated beam problem},
AIMS Math. \textbf{5} (2020), no.~5, 5261--5286.

\bibitem{GloverLazerMcKenna1989}
J.~Glover, A.~C. Lazer and P.~J. McKenna,
\emph{Existence and stability of large scale nonlinear oscillations in suspension bridges},
Z. Angew. Math. Phys. \textbf{40} (1989), 172--200.

\bibitem{Hale1988}
J.~K. Hale,
\emph{Asymptotic Behavior of Dissipative Systems},
Math. Surveys Monogr., vol.~25, Amer. Math. Soc., Providence, RI, 1988.

\bibitem{HolubovaMatas2003}
G.~Holubov\'a and A.~Matas,
\emph{Initial-boundary value problem for the nonlinear string-beam system},
J. Math. Anal. Appl. \textbf{288} (2003), 784--802.

\bibitem{Kang2016}
J.~R. Kang,
\emph{Asymptotic behavior of the thermoelastic suspension bridge equation with linear memory},
Bound. Value Probl. \textbf{2016} (2016), Article ID 206.

\bibitem{Kato1995}
T.~Kato,
\emph{Perturbation Theory for Linear Operators},
Classics in Mathematics, Springer, Berlin, 1995.

\bibitem{LazerMcKenna1987}
A.~C. Lazer and P.~J. McKenna,
\emph{Large scale oscillatory behavior in loaded asymmetric systems},
Ann. Inst. H. Poincar\'e Anal. Nonlinear \textbf{4} (1987), 243--274.

\bibitem{LazerMcKenna1990}
A.~C. Lazer and P.~J. McKenna,
\emph{Large-amplitude periodic oscillations in suspension bridges: some new connections with nonlinear analysis},
SIAM Rev. \textbf{32} (1990), no.~4, 537--578.

\bibitem{Litcanu2004}
G.~Li\c{t}canu,
\emph{A mathematical model of suspension bridges},
Appl. Math. \textbf{49} (2004), no.~1, 39--55.

\bibitem{MaZhong2005}
Q.~Z. Ma and C.~K. Zhong,
\emph{Existence of global attractors for the coupled system of suspension bridge equations},
J. Math. Anal. Appl. \textbf{308} (2005), 365--379.

\bibitem{MaZhong2009}
Q.~Ma and C.~Zhong,
\emph{Existence of strong solutions and global attractors for the coupled suspension bridge equations},
J. Differential Equations \textbf{246} (2009), 3755--3775.

\bibitem{Malik2004}
J.~Malik,
\emph{Mathematical modelling of cable-stayed bridges: existence, uniqueness, continuous dependence on data, homogenization of cable systems},
Appl. Math. \textbf{49} (2004), no.~1, 1--38.

\bibitem{McKenna1999}
P.~J. McKenna,
\emph{Torsional oscillations in suspension bridges revisited: fixing an old approximation},
Amer. Math. Monthly \textbf{106} (1999), 1--18.

\bibitem{McKennaTuama2001}
P.~J. McKenna and C.~Tuama,
\emph{Large torsional oscillations in suspension bridges revisited: vertical forcing creates torsional response},
Amer. Math. Monthly \textbf{108} (2001), 738--745.

\bibitem{McKennaWalter1987}
P.~J. McKenna and W.~Walter,
\emph{Nonlinear oscillations in a suspension bridge},
Arch. Ration. Mech. Anal. \textbf{98} (1987), no.~2, 167--177.

\bibitem{McKennaWalter1990}
P.~J. McKenna and W.~Walter,
\emph{Travelling waves in a suspension bridge},
SIAM J. Appl. Math. \textbf{50} (1990), no.~3, 703--715.

\bibitem{MessaoudiBonfohMukiawaEnyi2017}
S.~A. Messaoudi, A.~Bonfoh, S.~E. Mukiawa and C.~D. Enyi,
\emph{The global attractor for a suspension bridge with memory and partially hinged boundary conditions},
Z. Angew. Math. Mech. \textbf{97} (2017), no.~1, 1--14.

\bibitem{Mukiawa2025}
S.~E. Mukiawa and S.~A. Messaoudi,
\emph{On the well-posedness and long time dynamics for a coupled nonlinear bridge system with past history},
Appl. Math. Optim. \textbf{91} (2025), Article ID 53.

\bibitem{MukiawaEnyi2020}
S.~E. Mukiawa and C.~D. Enyi,
\emph{Well-posedness and attractors for a memory-type thermoelastic Timoshenko beam acting on shear force},
Taiwan. J. Math. \textbf{24} (2020), no.~5, 1179--1204.

\bibitem{ParkKang2011}
Y.~Park and J.~R. Kang,
\emph{Global attractors for the suspension bridge equations with nonlinear damping},
Q. Appl. Math. \textbf{69} (2011), 465--475.

\bibitem{Pazy1983}
A.~Pazy,
\emph{Semigroups of Linear Operators and Applications to Partial Differential Equations},
Applied Mathematical Sciences, vol.~44, Springer, New York, 1983.

\bibitem{Robinson2001}
J.~C. Robinson,
\emph{Infinite-Dimensional Dynamical Systems},
Cambridge University Press, Cambridge, 2001.

\bibitem{ZhongMaSun2007}
C.~K. Zhong, Q.~Z. Ma and C.~Y. Sun,
\emph{Existence of strong solutions and global attractors for the suspension bridge equations},
Nonlinear Anal. \textbf{67} (2007), 442--454.

\bibitem{Zhou2026}
J.~Zhou,
\emph{Exponential stability of the linear semigroup and the global attractor for a coupled suspension bridge system with past history},
preprint, arXiv:2610.11027v1, 2026.

\end{thebibliography}
\end{document}